\documentclass[10pt]{article}

\usepackage{mystyle,graphicx}
\usepackage{authblk}

\newcommand{\R}{\mathbb{R}}
\newcommand{\Z}{\mathbb{Z}}
\newcommand{\C}{\mathbb{C}}
\DeclareMathOperator{\symm}{Symm}

\newcommand{\X}{\mathcal{X}}
\newcommand{\Y}{\mathcal{Y}}
\newcommand{\G}{\mathcal{G}}
\newcommand{\F}{\mathcal{F}}
\newcommand{\bZ}{\mathbf{Z}}
\newcommand{\bm}{\mathbf{m}}
\newcommand{\pad}{\mathbf{M}}

\usepackage{enumitem}

\title{Computer-assisted proofs of Hopf bubbles and degenerate Hopf bifurcations}
\author[1]{Kevin Church\thanks{ kevin.church@umontreal.ca}}
\author[2]{Elena Queirolo\thanks{Corresponding author, elena.queirolo@tum.de}}
\affil[1]{Universit\'e de Montr\'eal, Montreal, Canada.}
\affil[2]{Technische Universit\"at M\"unchen, Munich, Germany.}
\date{\today}

\begin{document}
\maketitle

\begin{abstract}
We present a computer-assisted approach to prove the existence of Hopf bubbles and degenerate Hopf bifurcations in ordinary and delay differential equations. We apply the method to rigorously investigate these nonlocal bifurcation structures in the FitzHugh-Nagumo equation, the extended Lorenz-84 model and a time-delay SI model.
\end{abstract}


\section{Introduction}
The Hopf bifurcation is a classical mechanism leading to the birth of a periodic orbit in a dynamical system. In the simplest setting of an ordinary differential equation (ODE), a Hopf bifurcation \textit{generically} occurs when the linearization at a fixed point has a single pair of complex-conjugate imaginary eigenvalues. In this case, a perturbation by way of a scalar parameter would be expected, but not guaranteed, to result in a Hopf bifurcation. To prove the existence of the Hopf bifurcation, some non-degeneracy conditions must be checked. One such non-degeneracy condition pertains to the eigenvalues at the linearization: on variation of the distinguished scalar parameter, they must cross the imaginary axis transversally. We refer the reader to the papers \cite{Crandall1977,Faria1995,Jiang2020,Liu2011a,Rustichini1989} for some classical (and more recent) background concerning the Hopf bifurcation in the context of infinite-dimensional dynamical systems. A standard reference for the ODE case is the book of Marsden \& McCracken \cite{Marsden1976}.

\subsection{Hopf bubbles}
The Hopf bifurcation is \textit{generic} in one-parameter vector fields of ODEs. Roughly speaking, this means it is prevalent in a topological sense among all vector fields that possess an equilibrium point with a pair of complex-conjugate imaginary eigenvalues. In the world of two-parameter vector fields, the richness of possible bifurcations is much greater. In this paper, we are primarily interested in a type of degenerate Hopf bifurcation that results from a particular failure of the eigenvalue transversality condition. The failure we are interested in is where the curve of eigenvalues has a quadratic tangency with the imaginary axis, so the eigenvalues do not cross transversally and, importantly, do not cross the imaginary axis at all.

Before surveying some literature about this bifurcation pattern, let us construct a minimal working example. Consider the ODE system
\begin{align}
\label{ex.foldODE1}\dot x&=\beta x - y -x(x^2+y^2+\alpha^2)\\
\label{ex.foldODE2}\dot y&=x + \beta y -y(x^2+y^2+\alpha^2).
\end{align}
The reader familiar with the normal form of the Hopf bifurcation should find this familiar, but might be be unnerved by the $\alpha^2$ term, which is not present in the typical normal form. The linearization at the equilibrium $(0,0)$, produces the matrix 
$$\left(\begin{array}{cc}\beta-\alpha^2 & -1 \\ 1 & \beta-\alpha^2\end{array}\right),$$
which has the pair of complex-conjugate eigenvalues $\lambda = \beta - \alpha^2 \pm i$. Treating $\alpha$ as being a fixed constant, we have a supercritical Hopf bifurcation as $\beta$ passes through $\alpha^2$. Conversely, if $\beta>0$ is fixed and we interpret $\alpha$ as the parameter, there are \textit{two} supercritical Hopf bifurcations: as $\alpha$ passes through $\pm\sqrt{\beta}$. However, something problematic happens if $\beta=0$: the eigenvalue branches $\alpha\mapsto -\alpha^2 \pm i$ are tangent to the imaginary axis at $\alpha=0$ and do not cross at all, so the non-degeneracy condition of the Hopf bifurcation fails.

More information can be gleaned by transforming to polar coordinates. Setting $x=r\cos\theta$ and $y=r\sin\theta$ results in the equation for the radial component
$$\dot r = r(\beta - r^2 - \alpha^2),$$ while $\dot\theta=1$. There is a nontrivial periodic orbit (in fact, limit cycle) if $\beta-\alpha^2>0$, with amplitude $\sqrt{\beta-\alpha^2}$. Alternatively, there is a \emph{surface of periodic orbits} described by the graph of $\beta = r^2+\alpha^2$. For fixed $\beta>0$, the amplitude $r$ of the periodic orbit, as a function of $\alpha$, is $\alpha\mapsto\sqrt{\beta-\alpha^2}$, whose graph is half of an ellipse with semi-major axis $2\sqrt{\beta}$ and semi-minor axis $\sqrt{\beta}$. 

This non-local bifurcation structure, characterized by the connection of two Hopf bifurcations by a one-parameter branch of periodic orbits, has been given several names in different scientific fields. In intracellular calcium, it is frequently called a \textit{Hopf bubble} \cite{Dupont2016,Latulippe2018,Minicucci2021,Sneyd2017}. In infectious-disease modelling, where the bifurcation typically occurs at an endemic equilibrium, the accepted term is \textit{endemic bubble} \cite{Morshedy2021,LeBlanc2016,Liu2015d,Opoku-Sarkodie2022,Sherborne2018}. A more general definition of a (parametrically) non-local structure called \textit{bubbling} is given in \cite{Krisztin2011}. In the present work, we will refer to the structure as a \textit{Hopf bubble}, since that name is descriptive of the geometric picture, the bifurcation involved, and is sufficiently general to apply in different scenarios in a model-independent way.

Hopf bubbles can, in many instances, be understood as being generated by a codimension-two bifurcation; see Section \ref{sec-degenerate-intro}. However, this is not to say that they are rare. Aside from the applications in calcium dynamics and infectious-disease modelling in the previous paragraph, they have been observed numerically in models of neurons \cite{Ashwin2016}, condensed-phase combustion \cite{Margolis1988}, predator-prey models \cite{Braza2003}, enzyme-catalyzed reactions \cite{Hassard1993}, and a plant-water ecosystem model \cite{Wang2017a}. A recent computer-assisted proof also established the existence of a Hopf bubble in the Lorenz-84 model \cite{VandenBerg2021a}. We will later refer to the degenerate Hopf bifurcation that gives rise to Hopf bubbles as a \textit{bubble bifurcation}. The literature seems to not have an accepted name for this bifurcation, so we have elected to give it one here.

\subsection{Bautin bifurcation}
Bubble bifurcations are one type of dgenerate Hopf bifurcation. Another is the Bautin bifurcation, which occurs at a Hopf bifurcation whose first Lyapunov coefficient vanishes \cite{Kuznetsov}. Like the bubble bifurcation, the Bautin bifurcation is a codimension-two bifurcation of periodic orbits. To illustrate this bifurcation, consider the planar normal form
\begin{align}
\label{ex.Bautin1}\dot x&=\beta x - y -x(x^2+y^2)(\alpha - x^2-y^2)\\
\label{ex.Bautin2}\dot y&=x + \beta y -y(x^2+y^2)(\alpha - x^2-y^2).
\end{align}
Transforming to polar coordinates, the angular component decouples producing $\dot\theta=1$, while the radial component gives
$$\dot r = r(\beta + \alpha r^2 - r^4).$$
Nontrivial periodic orbits are therefore determined by the zero set of $r\mapsto\beta+\alpha r^2 - r^4$, which defines a two-dimensional smooth manifold. See later Figure \ref{fig:Bautin} for a triangulation of (part of) this manifold.

\subsection{Degenerate Hopf bifurcations and multi-parameter continuation}\label{sec-degenerate-intro}
The bubble bifurcation has been fully characterized by LeBlanc \cite{LeBlanc2016}, using the center manifold reduction and normal form theory for functional differential equations \cite{Faria1995} and a prior classification of degenerate Hopf bifurcations for ODEs by Golubitsky and Langford \cite{Golubitsky1981}. The study of the Bautin bifurcation goes back to 1949 with the work of Bautin \cite{Bautin}, and a more modern derivation based on normal form theory appears in the textbook of Kuznetsov \cite{Kuznetsov}.

Normal form theory and centre manifold reduction are incredibly powerful, providing both the direction of the bifurcation and criticality of bifurcating periodic orbits. The drawback is that they are computationally intricate and inherently local. In the present work, we advocate for an analysis of degenerate Hopf bifurcations in ODE and delay differential equations (DDE) by way of \textit{two-parameter continuation}, \textit{desingularization} and \textit{computer-assited proofs}. We will discuss the latter topics in the next section. 

\begin{figure}[!htb]
\centering\includegraphics[scale=0.45]{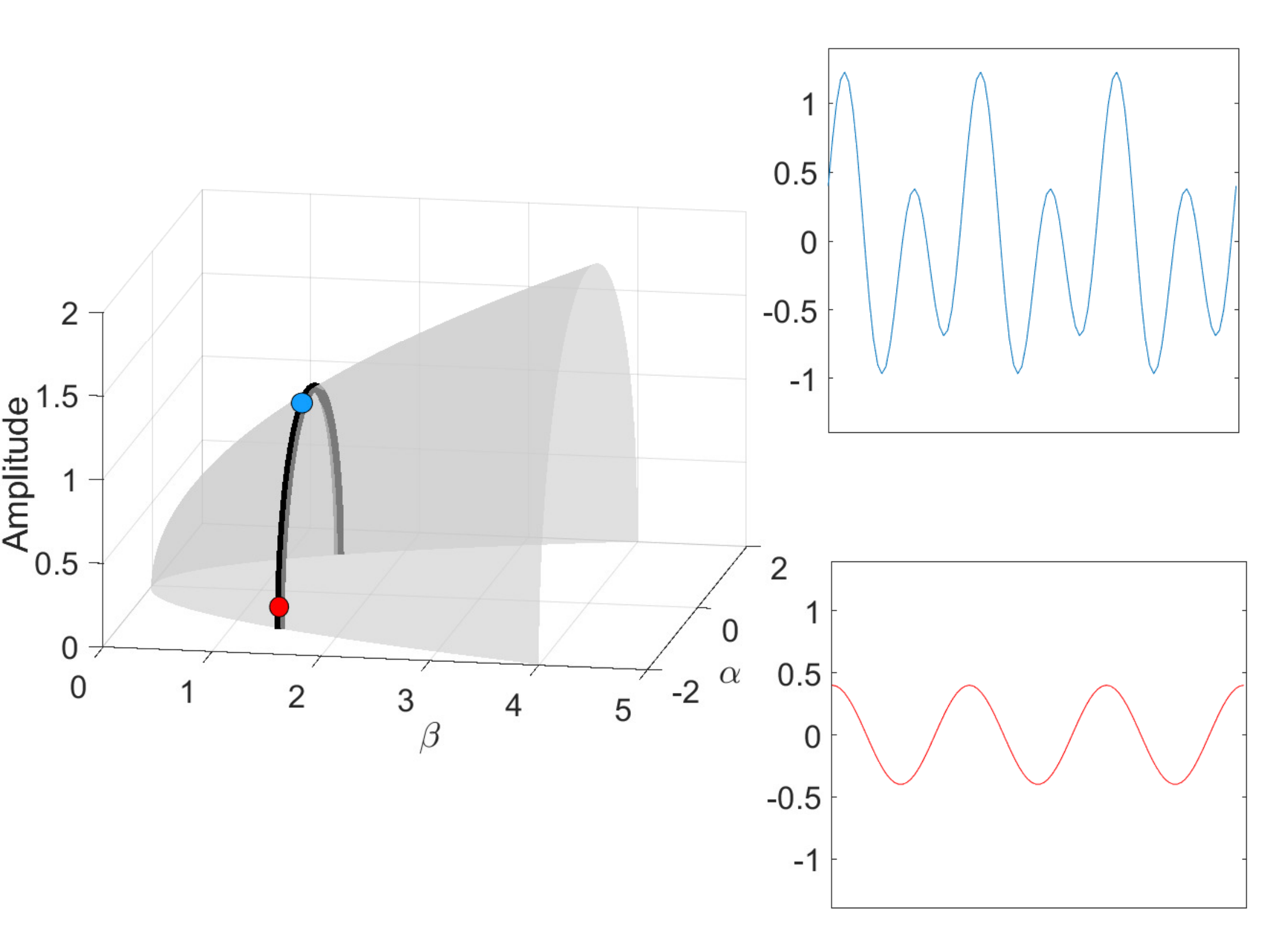}
\caption{Near a Hopf bifurcation, the bifurcating periodic orbit is close to a pure cosine (red solution) with small amplitude. Far from the bifurcation, the amplitude tends to grow and more Fourier modes are represented (blue solution). In a Hopf bubble, the amplitude grows from zero and then decreases to zero as a parameter ($\alpha$ in this figure) is varied montonically. Variation of a control parameter ($\beta$) can result in a smaller range of amplitudes and admissible interval over which the periodic orbit persists. If variation of the control parameter ultimately results in a collapse of the amplitude and admissibility ranges to a point, the result is a degenerate Hopf bifurcation.}
\label{cartoon-manifold}
\end{figure}

The intuition behind our continuation idea can be pictorially seen in Figure \ref{cartoon-manifold}, and understood analytically by way of our minimal working example, \eqref{ex.foldODE1}--\eqref{ex.foldODE2}, and its surface of periodic orbits described by the equation $\beta = \alpha^2+r^2$. There is an implicit relationship between two scalar parameters ($\alpha$ and $\beta$) and the set of periodic orbits of the dynamical system. At a Hopf bifurcation, one of these periodic orbits should retract onto a steady state. In other words, their amplitudes (relative to steady state) should become zero. Using a desingularization technique to topologically separate periodic orbits from steady states that they bifurcate from, we can use two-parameter continuation to numerically compute and continue periodic orbits as they pass through a degenerate Hopf bifurcation. With computer-assisted proof techniques, we can prove that the numerically-computed objects are close to true periodic orbits, with explicit error bounds. Further a posteriori analysis can then be used to prove the existence of a degenerate Hopf bifurcation based on the output of the computer-assisted proof of the continuation. 

\subsection{Rigorous numerics}\label{sec:newton-like}
Computer-assisted proofs of Hopf bifurcations have been completed in \cite{VandenBerg2021a} using a desingularization (sometimes called blow-up) approach, in conjunction with an \textit{a posteriori} Newton-Kantorvich-like theorem. We lean heavily on these ideas in the present work. The former desingularization idea, which we adapt to delay equations in Section \ref{sec-desingularization}, is used to resolve the fact that, at the level of a Fourier series, a periodic orbit that limits (as a parameter varies) to a fixed point at a Hopf bifurcation is, itself, indistinguishable from that fixed point. Our methods of computer-assisted proof are based on contraction mappings, and it is critical that the objects we prove are isolated as fixed points. The desingularization idea exploits the fact that an amplitude-based change of variables can be used to develop an equivalent problem where the representative of a periodic orbit consists of a Fourier series, a fixed point and a real number. The latter real number represents a signed amplitude. This reformulation results in fixed points being spatially islolated from periodic orbits, thereby allowing contraction-based computer-assisted proofs to succeed.

The other big point of inspiration in this work is \textit{validated multi-parameter continuation}. The technique was developed in \cite{Gameiro2016a} for continuation in general Banach spaces, and applied to some steady state problems for PDEs. We will overview this method in Section \ref{sec-continuation}. It is based on a combination of a two-dimensional analogue of pseudo-arclength continuation and \textit{a posteriori}, uniform validation of zeroes of nonlinear maps by a Newton-Kantorovich-like theorem applied to a Newton-like operator.

\subsection{Contributions and applications}
The main contribution of the present work is a general-purpose code for computing and proving two-parameter families of periodic orbits in polynomial delay differential equations. Equations of advanced or mixed-type can similarly be handled; there is no restriction whether delays are positive or negative. Ordinary differential equations can also be handled as a special case. Orbits can be proven in the vicinity of (and at) Hopf bifurcations, whether these are non-degenerate or degenerate. The first major release of the library \texttt{BiValVe} (\textbf{Bi}furcation \textbf{Val}idation \textbf{Ve}nture, \cite{CodeURL}) is being made in conjunction with the present work, and builds on an earlier version of the code associated to the work \cite{Vandenberg2021_0,VandenBerg2021a}. Some non-polynomial delay differential equations can be handled using the \textit{polynomial embedding technique}. The existence of Hopf bifurcation curves and degenerate Hopf bifurcations can then be completed by post-processing of the output of the computer-assisted proof.

To explore the applicability of our validated numerical methods, we explore Hopf bifurcations and degenerate Hopf bifurcations in
\begin{itemize}
\item the extended Lorenz-84 model (ODE)
\item a time-delay SI model (DDE)
\item the FitzHugh-Nagumo equation (ODE)
\item an ODE with complicated branches of periodic orbits (ODE)
\end{itemize}
The first two examples replicate and extend some of the analysis appearing in \cite{LeBlanc2016,VandenBerg2021a} using our computational scheme. The third examples provides, to our knowledge, the first analytically verified results on degenerate Hopf bifurcations and Hopf bubbles in that equation. The final example is a carefully designed ODE that exhibits the degenerate Hopf bifurcation in addition to folding and pinching in some projections of the manifold.

\subsection{Overview of the paper}
Section \ref{sec-continuation} serves as an overview of two-parameter continuation, both in the finite-dimensional and infinite-dimensional case. We introduce the continuation scheme for periodic orbits near Hopf bifurcations in Section \ref{sec3} in the general case of delay differential equations. Some technical bounds for the computer-assisted proofs in Section \ref{sec:CAP}. A specification to ordinary differential equations is presented in Section \ref{sec:ode}. In Section \ref{sec:analytical_proofs}, we connect the computer-assisted proofs of the manifold of periodic orbits to analytical proofs of Hopf bifurcation curves, Hopf bubbles, and degenerate Hopf bifurcations. Our examples are presented in Section \ref{sec:examples}, and we complete a discussion and comment on future research directions in Section \ref{sec:discussion}.

\subsection{Notation}\label{sec:notation}
Given $n\in\mathbb{N}$, denote $\C^\Z_n$ the vector space of $\Z$-indexed sequences of elements of $\C^n.$ Denote $\mbox{Symm}(\C^\Z_n)$ the proper subspace of $\C^\Z_n$ consisting of symmetric sequences; $z\in \mbox{Symm}(\C^\Z_n)$ if and only if $z\in\C^\Z_n$ and $z_k=\overline{z_{-k}}$ for all $k\in\Z$. For any subspace $U\subset\C^\Z_n$ closed under (componentwise) complex conjugation, denote $\symm(U) = U\cap\symm(\C^\Z_n)$. We will sometimes drop the subscript $n$ on $\C^\Z_n$ when the context is clear.

Given $\nu>1$, we denote $\ell_\nu^1(\C^n)$ the subspace of $\C_n^\Z$ whose elements $z$ satisfy $||z||_\nu\bydef \sum_{k\in\Z}|z_k|\nu^{|k|}<\infty$. $|\cdot|$ will always be taken to be the norm on $\C^n$ induced by the standard inner product. Denote $K_\nu(\C^n)$ the subspace of $\C_n^\Z$ whose elements $z$ satisfy $||z||_{\nu,K}\bydef \sum_{k\in\Z}(\nu^{|k|}/|k|)|z_k|<\infty.$ Introduce a bilinear form on $\langle\cdot,\cdot\rangle$ on $\ell_\nu^1(\C)$ as follows:
$$\langle v,w\rangle = \sum_{k\in\Z}v_k\overline{w_k}.$$
For $v,w\in\C^\Z_1$, define their convolution $v*w$ by 
\begin{align}\label{def:convolution}
(v*w)_k = \sum_{k_1+k_2=k}v_{k_1}w_{k_2}
\end{align}
whenever this series converges. Convolution is commutative and associative for sequences in $\ell_\nu^1(\C)$. In this case, we define multiple convolutions (e.g.\ triple convolutions $a*b*c$) inductively, by associativity.

A function $g:\ell_\nu^1(\C^n)\rightarrow\ell_\nu^1(\C)$ is a \emph{convlution polynomial of degree $q$} if
\begin{align*}
g(z)&=c_0 + \sum_{k=1}^q\sum_{p\in\mathcal{M}_k}c_p (z_{p_1}*\cdots*z_{p_k}),
\end{align*}
for $c_p\in\C$, where $\mathcal{M}_k$ denotes the set of $k$-element multisets of $\{1,\dots,n\}$, and each multiset $p$ is identified with the unique tuple $(p_1,\dots,p_k)\{1,\dots,n\}^k$ such that $p_i\leq p_{i+1}$ for all $i=1,\dots,k-1$. Analogously, $g:\ell_\nu^1(\C^n)\rightarrow\ell_\nu^1(\C^m)$ is a \emph{convolution polynomial of degree $q$} if $z\mapsto g(z)_{\cdot,j}\in\ell_\nu^1(\C)$ is a convolution polynomial of degree $q$, for $j=1\dots,m$.

If $X$ is a vector space, $0_X$ will denote the zero vector in that space. If $X$ is a metric space and $U\subset X$, we denote $U^\circ$ its interior, $\partial U$ its boundary, and $\overline{U}$ its closure.

An \emph{interval vector} in $\R^k$ for some $k\geq 1$ is a subset of the form $v = [a_1,b_1]\times\cdots\times[a_k,b_k]$ for real scalars $a_j,b_j\in\R$, $j=1,\dots,k$. If $\R^k$ is equipped with a norm $||\cdot||$ we define $||v||=\sup_{w\in v}||w||$. Similarly, an interval vector in $\C^k$ is a product $v = A_1\times\cdots\times A_k$, where each $A_j$ is a closed disc in $\C$. We define the norm of a complex interval vector as $||v|| = \sup_{w\in v}||w||$ whenever $||\cdot||$ is a norm on $\C^k$. In each case, be it real or comlex, the \emph{unit interval vector} is the unique interval vector $V$ with the property that $||V||=1$, and for any other interval vector $v$ with $||v||=1$, the inclusion $v\subseteq V$ is satisfied.

\section{Validated two-parameter continuation}\label{sec-continuation}
In this section, we review validated two-parameter continuation. Our presentation will loosely follow \cite{Gameiro2016a}. Some noteworthy changes compared to the references are that we work in a complexified (as opposed to strictly real) vector space, which causes some minor difficulties at the level of implementation.

We first review the continuation algorithm as it applies to finite-dimensional vector spaces in Section \ref{sec-continuation-findim}. We make comments concerning implementation in Section \ref{sec-technicals}. We then describe how it is extended to general Banach spaces in Section \ref{sec-continuation-infdim}. Validated continuation (i.e.\ computer-assisted proof) is discussed in Section \ref{sec-continuation-proofs}.

\subsection{Continuation in a finite-dimensional space}\label{sec-continuation-findim}
Let $\mathcal{X}$ and $\mathcal{Y}$ be finite-dimensional vector spaces over the field $\R$, with $\dim \X = \dim \Y + 2$, and consider a map $\G:\X\rightarrow \Y$. We are interested in the zero set of $\G$. Given the codimension of $\G$, we expect zeroes to be in a two-dimensional manifold in $\X$. In what follows, equality with zero should be understood in the sense of numerical precision. For example, if we write $\G(x)=0$, what we mean is that $||\G(x)||$ is machine precision small.

Let $\hat x_0\in \X$ satisfy $G(\hat x_0)=0$, and suppose $D\G(\hat x_0)$ has two-dimensional kernel. This property is generic. Let $\{\hat\Phi_1,\hat\Phi_2\}$ span the kernel. The tangent plane $\mathcal{T}_{\hat x_0}\mathcal{M}$ at $\hat x_0$ of the two-dimensional solution manifold $\mathcal{M}$ is therefore spanned by $\{\hat\Phi_1,\hat\Phi_2\}$. If $\epsilon_1,\epsilon_2$ are small, we have
$$\G(\hat x_0 + \epsilon_1\hat\Phi_1 + \epsilon_2\hat\Phi_2) = \G(\hat x_0) + \epsilon_1 D\G(\hat x_0)\hat\Phi_1 + \epsilon_2 D\G(\hat x_0)\hat\Phi_2 + o(|\epsilon|) = o(|\epsilon|)$$ by Taylor's theorem, provided $\G$ is differentiable at $\hat x_0$, and $|\epsilon|$ is any Euclidean norm of $\epsilon=(\epsilon_1,\epsilon_2)$. If $\G$ is twice continuously differentiable, the error is improved to $O(|\epsilon|^2)$. Therefore, new candidate zeroes of $\G$ can be computed using $\hat x_0$ and a basis for the tangent space. This idea is at the heart of the continuation.

\begin{figure}
\centering\includegraphics[scale=0.65]{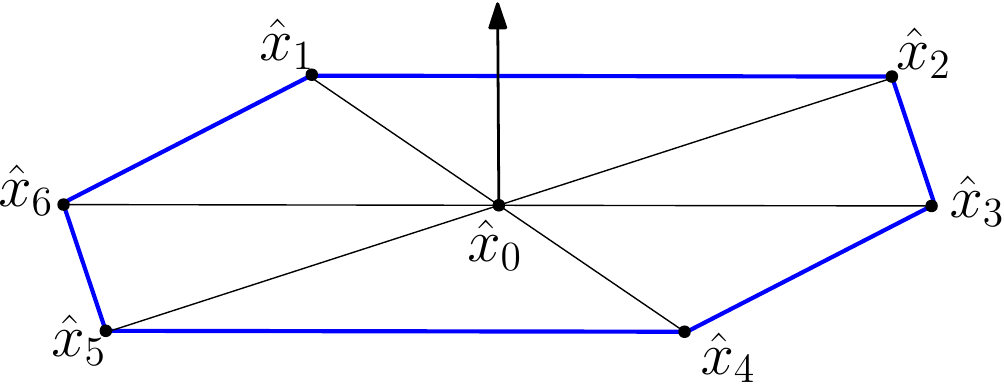}
\caption{The hexagonal patch with normal vector indicated by an arrow. In this initial configuration, every node except for $\hat x_0$ is a boundary node. The boundary edges are blue.}\label{fig-localhexagon}
\end{figure}

The continuation from $\hat x_0$ is done by way of iterative triangulation of the manifold $\mathcal{M}$. First, we compute an orthonormal basis of $\mathcal{T}_{\hat x_0}\mathcal{M}$ by applying the Gram-Schmidt process to $\{\hat\Phi_1,\hat\Phi_2\}$; see Section \ref{sec-technicals} for some technical details. Using this orthonormal basis to define a local coordinate system, six vertices of a regular hexagon are computed around $\hat x_0$ at a specified distance $\sigma$ from $\hat x_0$; see Figure \ref{fig-localhexagon}. Let these vertices be denoted $\hat x_1$ through $\hat x_6$, arranged in counterclockwise order (relative to the local coordinate system) around $\hat x_0$. Note that this means $$\hat x_i = \hat x_0 + \epsilon_{j,1}\hat\Phi_1 + \epsilon_{j,2}\hat\Phi_2$$ for some small $\epsilon_{j,1},\epsilon_{j,2}$, which means that $G(\hat x_j)\approx 0$ for $j=1,\dots,6$. Each of these candidate zeroes $\hat x_j$ are then refined by applying Newton's method to the map
\begin{align}\label{newton-G}x\mapsto \G_j(x) = \left(\begin{array}{c}\G(x) \\ \langle \hat\Phi_1, x - \hat x_j\rangle \\ \langle \hat\Phi_2, x - \hat x_j\rangle   \end{array}\right).\end{align}
The two added inner product equations ensure isolation of the solution (hence quadratic convergence of the Newton iterates) and that the Newton correction is perpendicular to the tangent plane. This map can similarly be used to refine the original zero $\hat x_0$.

This initial hexagonal ``patch" is itself formed by six triangles; see Figure \ref{fig-localhexagon}. The continuation algorithm proceeds by selecting one of the boundary vertices (i.e.\ one of the vertices $\hat x_1$ through $\hat x_6$) and attempting to ``grow" the manifold further. We describe this ``growth" phase below. However, first, some terminology. The vertices will now be referred to as \emph{nodes}. An \emph{edge} is a line segment connecting two nodes, and they will be denoted by pairs of nodes: $\{\hat x_i,\hat x_j\}$ for node $\hat x_i$ connected to $\hat x_j$. Two nodes are \emph{incident} if they are connected by an edge. A \emph{simplex} is the convex hull of three edges that form a triangle. Once the hexagonal patch is created, the data consists of:
\begin{itemize}
\item the nodes $\hat x_0,\dots,\hat x_6$;
\item 
the ``boundary edges" $\{\hat x_1,\hat x_2\}, \dots,\{\hat x_5,\hat x_6\},\{\hat x_6,\hat x_1\}$; 
\item the six simplices formed by triangles with $\hat x_0$ as one of the nodes.
\end{itemize}
Two simplices are \emph{adjacent} if they share an edge. An \emph{internal simplex} is a simplex that is adjacent to three other simplices, and it is a \emph{boundary simplex} otherwise. Therefore, the six simplices of the initial patch are all considered boundary since they are adjacent to exactly two others. Similarly, an edge of a simplex can be declared boundary or internal; internal edges are those that are shared with another simplex, and boundary edges are not. A \emph{boundary node} is any node on a boundary edge, a \emph{frontal node} is a boundary node on an edge shared by two simplices, and an \emph{internal node} is a node that is not a boundary node.

\begin{figure}
\centering
\includegraphics[scale=0.55]{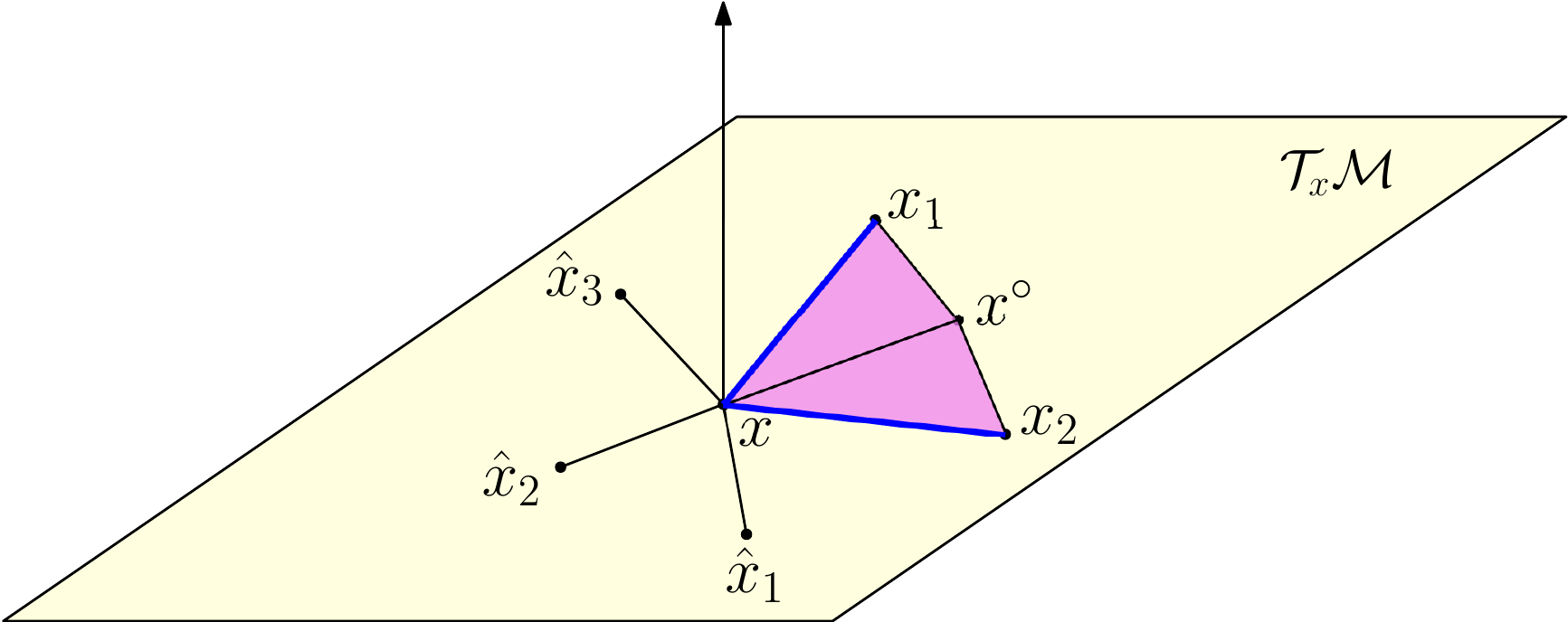}
\caption{Cartoon diagram of the simplex growth phase. For visualization purposes, we think of $x_1$ and $x_2$ as being in the tangent space $\mathcal{T}_x\mathcal{M}$; in reality, these should be close to the tangent space but not strictly contained in it. The boundary edges $\{x,x_1\}$ and $\{x,x_2\}$ are blue, and these form a simplices with the interior edge $\{x,x^\circ\}$ that generates the \textit{gap complement direction}. The predictor ``fan", after being mapped back to the tangent space, consists of the vertices $\hat x_1$, $\hat x_2$ and $\hat x_3$. The four new simplices that can be formed, namely $\{x,\hat x_1,\hat x_2\}$, $\{x,\hat x_2,\hat x_3\}$, $\{x,\hat x_1,x_2\}$ and $\{x,\hat x_3,x_1\}$, share the same angle at the $x$ vertex.}
\label{fig-growsimplex}
\end{figure}

Let $x$ be a chosen frontal node (see Section \ref{sec-technicals-nodes} for a general discussion on frontal node selection). By construction, $x$ is part of (at least) three distinct edges, two of which connect to boundary nodes, and at least one that connects to an internal node. The algorithm to \emph{grow a simplex} is as follows. See Figure \ref{fig-growsimplex} for a visualization.
\begin{enumerate}
\item Compute an orthonormal basis for the tangent space $\mathcal{T}_{x}\mathcal{M}$.
\item Compute the \emph{(average) gap complement direction}. Let $x^\circ_1,\dots,x^\circ_m$ denote the list of internal nodes nodes incident to $x$. The gap complement direction is $y_c = -x + \frac{1}{m}\sum_{i=1}^m x_i^\circ$.
\item Let $x_1$ and $x_2$ denote the boundary nodes incident to $x$; form the edge directions $y_1 = x_1 - x$ and $y_2 = x_2 - x$.
\item Orthogonally project $y_c$, $y_1$ and $y_2$ onto the tangent plane $\mathcal{T}_x\mathcal{M}$. Let these projections be denoted $Py_c$, $Py_1$ and $Py_2$. Introduce a two-dimensional coordinate system on $\mathcal{T}_x\mathcal{M}$ by way of a ``unitary" (see Section \ref{sec-technicals}) transformation to $\R^2$. Let $\tilde Py_c$, $\tilde Py_1$ and $\tilde Py_2$ denote the representatives in $\R^2$.
\item In the local two-dimensional coordinate system, compute the counter-clockwise (positive) angle $\theta_1$ required to complete a rotation from $\tilde Py_c$ to $\tilde Py_1$, and the counter-clockwise (positive) angle $\theta_2$ required for rotation from $\tilde Py_c$ to $\tilde Py_2$. The \emph{gap angle} $\gamma$ and \emph{orientation} $\rho$ is
$$\gamma = \max\{\theta_1-\theta_2,\theta_2-\theta_1\},\hspace{1cm}\rho = \left\{\begin{array}{ll}1,&\theta_2>\theta_1 \\ 2 & \theta_2<\theta_1 \end{array}\right.$$
\item If $\gamma<\frac{\pi}{6}$ then \emph{close the gap}: add the the triangle formed by the nodes $\{x,x_1,x_2\}$ to the list of simplices, flag it as a boundary simplex, flag $x$ as an internal node, and conclude the growth step. Otherwise, proceed to step 7.
\item Generate the predictor ``fan" in $\R^2$: let $k=\min\left\{1,\lfloor 3\theta/\pi\rfloor\right\}$ and define the predictors
$$\tilde P\tilde y_j = R(j\gamma/k)\tilde Py_\rho,\hspace{1cm}j=1,\dots,k,$$
for $R(\theta)$ the $2\times 2$ counterclockwise rotation matrix through angle $\theta$.
\item Invert the unitary transformation and map $\tilde P\tilde y_j$ into the tangent plane $\mathcal{T}_x\mathcal{M}$; let the result be the vectors $P\tilde y_j$, $j=1,\dots,k$.
\item Define predictors $\hat x_j = x + \sigma P\tilde y_j$ for $j=1,\dots,k$, where $\sigma$ is a user-specified step size. Refine them using Newton's method applied to \eqref{newton-G}, where $\hat\Phi_1$ and $\hat\Phi_2$ are now the orthonormal basis for $\mathcal{T}_x\mathcal{M}$.
\item Add the triangles formed by nodes $\{x,\hat x_1,\hat x_2\}$, $\{x,\hat x_2,\hat x_3\}$, $\dots$, $\{x,\hat x_{k-1},\hat x_k\}$. To the list of simplices, flag them as boundary simplices, and flag $x$ as an internal node. Do the same with the triangles formed by $\{x,\hat x_1,*\}$ and $\{x,\hat x_k,*\}$, where $*$ denotes ones of $x_1$ and $x_2$, depending on orientation $\rho$.
\end{enumerate}

At the end of the simplex growth algorithm, one or more boundary simplices is added to the dictionary, and one additional node will be flagged as internal. The structure of this algorithm ensures that each simplex added to the dictionary will be adjacent to exactly two others, indicating that our convention of internal vs.\ boundary simplices is effective. Once the growth algorithm is complete, another frontal node is selected and the growth phase is repeated. This continues until a sufficient portion of the manifold has been computed (i.e.\ a user-specified exit condition is reached), or until a Newton's method correction fails. 

\subsection{Comments on implementation}\label{sec-technicals}
Here we collect some remarks concerning implementation of the finite-dimensional continuation. These comments may be applicable for general continuation problems, but we will often emphasize our specific situation which is continuation of periodic orbits in delay differential equations.

\subsubsection{Complex rotations for the kernel of $D\G(x)$}
First, generating the kernel elements $\{\hat\Phi_1,\hat\Phi_2\}$ of $D\G(\hat x)$ for an approximate $\hat x$ must be handled in a way that is appropriate to the space $\X$. In our problem, $\X$ has additional structure: it is a subset of a complexified real vector space equipped with a discrete symmetry. However, in our implementation (that is, in the environment of MATLAB) we work on a generic complex vector space without this symmetry, so when we compute kernel elements using QR decomposition, the computed vectors are not necessarily in $\X$. We fix this by performing a complex rotation to put the kernel back into $\X$. This can always be done because $\hat\Phi_1$ and $\hat\Phi_2$, as computed using QR, are always $\C$-linearly independent.

\subsubsection{Orthonormal basis for the tangent space}\label{sec-Lmap-tangentspace}
The next point concerns the ``orthonormal basis" of $\mathcal{T}_{\hat x}\mathcal{M}$. Let us be a bit more precise. In our problem, $\mathcal{X}$ is a product of the form $\R^n\times \mathcal{V}$, where $\mathcal{V}\subset\C^{2M+1}$ is characterized by $v = (v_{-M},v_{-M+1},\dots,v_{M-1},v_M)\in\mathcal{V}$ satisfies $v=Sv$, where $Sv = (\overline{v_M},\overline{v_{M-1}},\dots,\overline{v_{-M+1}},\overline{v_{-M}})$. Consider the standard inner product $\langle\cdot,\cdot\rangle$ on $\C^{n+2M+1}$. Once a basis $\{\hat\Phi_1,\hat\Phi_2\}$ of $\mathcal{T}_{\hat x}\mathcal{M}$ has been computed, these basis vectors can be interpreted as being elements of $\C^{n+2M+1}$, and we say that they are orthogonal if $\langle\hat\Phi_1,\hat\Phi_2\rangle = 0$. It is straightforward to verify that the Gram-Schmidt process applied to this basis produces yet another basis of $\mathcal{T}_{\hat x}\mathcal{M}$ (that is, it does not break the symmetry), and that the new basis is orthonormal with respect to $\langle\cdot,\cdot\rangle$.

\subsubsection{Node prioritization for the simplex growth algorithm}\label{sec-technicals-nodes}
A suitable selection of a frontal node for simplex growth might not be obvious. First, nodes that are more re-entrant (i.e.\ have many edges incident to them) are generally given higher priority. This is because such nodes are more likely to have the simplex growth algorithm perform the \emph{close the gap} sub-routine. We want to avoid having thin simplices, so prioritizing the closing off of re-entrant nodes takes priority. After this, typically grow from the ``oldest" boundary node.

\subsubsection{Local coordinate system in $\R^2$ for the tangent space}
Finally, we must discuss the generation of a two-dimensional coordinate system for $\mathcal{T}_x\mathcal{M}$. If $\hat\Phi_1$ and $\hat\Phi_2$ are an orthonormal basis for the tangent space, then the map $$y\mapsto Ly = (\hat\Phi_1^*y,\hat\Phi_2^*y)$$ is invertible, where $\hat\Phi_i^*$ denotes the conjugate transpose of $\Phi_i$. Writing $Ly = (u,v)\in\C^2$ for $y\in\X$, one can verify that $(u,v)$ is the unique solution of
$$y = \hat\Phi_1 u + \hat\Phi_2 v.$$
In particular, the range of this map is $\R^2$ (that is, each of $u$ and $v$ is real) for our specific problem, where the space $\X$ is $\R^n\times\mathcal{V}$, with $\mathcal{V}$ having the symmetry described two paragraphs prior. The inverse map is $$(u,v)\mapsto L^{-1}(u,v) =  \hat\Phi_1 u + \hat\Phi_2 v.$$ Moreover, this map is unitary in the sense that $\langle Ly, Lw\rangle = \langle y,w\rangle$ for $y,w\in\X$, where on the left-hand side we take the standard inner product on $\R^2$. It is for these reasons that rotations on $\R^2$ performed after the action of $L$ are consistent with rotations in the tangent space relative to the gap complement direction.

\subsection{Continuation in a Banach space}\label{sec-continuation-infdim}
Two-parameter continuation can be introduced generally in a Banach space, but here we take the perspective that $\G:\X\rightarrow\Y$ of Section \ref{sec-continuation-findim} is a projection of a nonlinear map $G:X\rightarrow Y$, and we are \emph{actually} interested in computing the zeroes of $G$, rather than those of $\G$. That is to say, there exist projection operators $\pi_\X:X\rightarrow\X$, $\pi_\Y:Y\rightarrow\Y$, and associated embeddings $i_\X:\X\hookrightarrow X$, $i_\Y:\Y\hookrightarrow Y$, such that $\G = \pi_\Y G\circ i_\X$. 

While the zeroes of $G$ are what we want, we will still use the approximate zeroes of $\mathcal{G}$ to rigorously enclose them. Abstractly, let $\hat x_0,\hat x_1,\hat x_2\in\X$ be the three nodes of a simplex. Then we may introduce a coordinate system on this simplex as follows: write each element $\hat x_s$ of this simplex as a unique linear combination
\begin{align}\label{simplex-xhat}\hat x_s = \hat x_0 + s_1(\hat x_1 - \hat x_0) + s_2(\hat x_2-\hat x_0)\end{align}
for $s = (s_1,s_2)\in\Delta = \{(a,b)\in\R^2 : 0\leq a,b\leq 1,\hspace{1mm}0\leq a+b\leq 1\}$. Let $\{\hat\Phi_{j,1},\hat\Phi_{j,2}\}$ for $j=0,1,2$ denote an orthonormal basis for the kernel of $D\G(\hat x_j)$. We can then form the \emph{interpolated kernels}
$$\hat\Phi_{s,i} = \hat\Phi_{0,i} + s_1(\hat\Phi_{1,i} - \hat\Phi_{0,i}) + s_2(\hat\Phi_{2,i}-\hat\Phi_{0,i})$$
for $i=1,2$. We introduce a nonlinear map $G_s:X\rightarrow Y\times\R^2$, 
\begin{align}\label{Gmap_s}
G_s(x) = \left(\begin{array}{c} G(x) \\ \langle \hat\Phi_{s,1} ,\pi_\X x - \hat x_s \rangle \\ \langle \hat\Phi_{s,2}, \pi_\X x - \hat x_s \rangle \end{array}\right).
\end{align}
The objective is to prove that $G_s$ has a unique zero close to $\hat x_s$ for each $s\in\Delta$. If this can be proven, and in fact $G_s$ is $C^k$ for some $k\geq 1$, then one can prove \cite{Gameiro2016a} that the the zero set of $G$ is a $C^k$, two-dimensional manifold. If this same can be proven for a collection of simplices, then the zero set is globally (i.e.\ on the union of the cobordant simplicial patches) a $C^k$ manifold over all patches that can be proven; that is, the transition maps between patches are $C^k$.

\begin{remark}
There is a subtle point concerning the relative orientations of the individual kernel elements $\hat\Phi_{0,i}$, $\hat\Phi_{1,i}$ and $\hat\Phi_{2,i}$ that are used to define the interpolation $\hat\Phi_{s,i}$. The validated continuation is unstable (and can even fail outright) if the interpolated kernels $\hat\Phi_{s,i}$ vary too much over $s\in\Delta$. If these vectors all lived in the exact same two-dimensional tangent space, this would be very straightforward; we could simply (real) rotate and/or reflect each basis $\{\hat\Phi_{1,i},\hat\Phi_{2,i}\}$ for $i=1,2$ so that they matched $\{\hat\Phi_{1,0},\hat\Phi_{2,0}\}$ exactly. However, these vectors live in different tangent spaces, so it is not as easy. The more consistent differential-geometric way to solve the problem would be to parallel transport the tangent basis $\{\hat\Phi_{1,0},\hat\Phi_{2,0}\}$ to the other two nodes and use these as bases for the tangent spaces there. We do nothing so sophisticated. We merely perform (real) rotations or reflections of the orthonormal bases $\{\hat\Phi_{1,i},\hat\Phi_{2,i}\}$ in their respective tangent spaces (two-dimensional) for $i=1,2$ in such a way that, in norm, these are as close as possible to $\{\hat\Phi_{1,0},\hat\Phi_{2,0}\}$. In practice, this has the effect of promoting enough alignment of the bases that proofs are feasible. We emphasize that this alignment process is only needed for the validated continuation; using misaligned tangent bases is not a problem for the steps described in Section \ref{sec-continuation-findim}.
\end{remark}

\subsection{Validated continuation and the radii polynomial approach}\label{sec-continuation-proofs}
The \emph{radii polynomial approach} is essentially a Newton-Kantorovich theorem with an approximate derivative, approximate inverse, and domain parametrization. It will be used to connect the approximate zeroes of $\G$ with exact zeroes of $G$ from the previous section. We state it generally for a family of $s\in\Delta$ -dependent maps $F_s:X_1\rightarrow X_2$. We include a short proof for completeness.
\begin{theorem}\label{thm-rpa}
Suppose $F_s:X_1\rightarrow X_2$ is differentiable for each $s\in\Delta$, where $X_1$ and $X_2$ are Banach spaces. Let $\hat x_s\in X_1$ for all $s\in\Delta$. Suppose there exist for each $s\in\Delta$ a bounded linear operator $A_s^\dagger:X_1\rightarrow X_2$, a bounded and injective linear operator $A_s:X_2\rightarrow X_1$, and non-negative reals $Y_0$, $Z_0$, $Z_1$ and $Z_2=Z_2(r)$ such that
\begin{align}
\label{thm-rpa-Y0}||A_sF_s(\hat x_s)||&\leq Y_0\\
||I-A_sA_s^\dagger||_{B(X_1,X_1)}&\leq Z_0\\
\label{thm-rpa-Z1}||A_s(DF_s(\hat x_s)-A_s^\dagger)||_{B(X_1,X_1)}&\leq Z_1\\
\label{thm-rpa-Z2}||A_s(DF_s(\hat x_s+\delta) - DF_s(\hat x_s))||_{B(X_1,X_1)}&\leq Z_2(r),\hspace{5mm}\forall \delta\in \overline{B_r(\hat x_s)}\subset X,
\end{align}
for all $s\in\Delta$, where $B_r(\hat x_s)$ is the closed ball of radius $r$ centered at $\hat x_s$, and $||\cdot||_{B(X_1,X_1)}$ denotes the induced operator norm on $X_1$. Suppose there exists $r_0>0$ such that the \emph{radii polynomial}
$$p(r) = rZ_2(r) + (Z_1+Z_0-1)r + Y_0$$ satisfies $p(r_0)<0$. Then for each $s\in\Delta$, there is a unique $x_s\in B_{r_0}(\hat x_s)$ such that $F_s(x_s)=0$. If $(s,x)\mapsto F_s(x)$ and $s\mapsto A_s$ are $C^k$, then the same is true of $s\mapsto x_s$.
\end{theorem}

\begin{proof}
Define the Newton-like operator $T_s(x) = x - A_sF_s(x)$. We will show that $T_s$ is a contraction on $\overline{B_{r_0}(\hat x_s)}$, uniformly for $s\in\Delta$. First, write $x\in\overline{B_{r_0}(\hat x_s)}$ in the form $x = \hat x_s + \delta$ for some $||\delta||_\X\leq r$. Then
\begin{align*}
||DT_s(x)||&=||I-A_sDF_s(x)||\\
&\leq ||I - A_sA_s^\dagger|| + ||A_s(A_s^\dagger - DF_s(\hat x_s))|| + ||A_s(DF_s(\hat x_s) - DF_s(\hat x_s+\delta))||\\
&\leq Z_0+Z_1+Z_2(r),
\end{align*}
Now, using the triangle inequality and the mean-value inequality,
\begin{align*}
||T_s(x)-\hat x_s||&\leq||T_s(\hat x_s+\delta) - T_s(\hat x_s)|| + ||T_s(\hat x_s) - \hat x_s||\\
&\leq r\sup_{t\in[0,1]}||DT_s(\hat x_s + t\delta)||_{B(X_1,X_1)} + ||A_sF_s(\hat x_s)||\\
&\leq (Z_0 + Z_1 + Z_2(r))r + Y_0.
\end{align*}
Choosing $r=r_0$, the radii polynomial implies that $||T_s(x)-\hat x_s||<r$, so $T_s$ is a self-map on $\overline{B_{r_0}(\hat x_s)}$ with its range in the interior. Moreover, since $p(r_0)<0$, we get that $Z_0+Z_1+Z_2(r_0)<1$, which proves that $T_s:\overline{B_{r_0}(\hat x_s)}\rightarrow B_{r_0}(\hat x_s)$ is a contraction (uniformly in $s\in\Delta$). By the Banach fixed point theorem, $T_s$ has a unique fixed point $x_s\in B_{r_0}(\hat x_s)$ for each $s\in\Delta$, and $s\mapsto x_s$ is $C^k$ provided the same is true of $(s,x)\mapsto T_s(x)$. Since $A_s$ is injective, $F_s$ has a unique zero in $B_{r_0}(\hat x_s)$ if and only if $T_s$ has a unique fixed point there.
\end{proof}
For our problem, injectivity of $A_s$ will always follow from a successful verification of $p(r_0)<0$. See Lemma \ref{lemma-injectivity}.

Let $\hat x_s$ be the convex combination defined by \eqref{simplex-xhat} for the simplex nodes $\hat x_0$, $\hat x_1$ and $\hat x_2$. We will say this simplex has been \emph{validated} if we successfully find a radius $r_0$ such that the conditions of the radii polynomial theorem are successful for the nonlinear map $G_s:X\rightarrow Y\times\R^2$. 

In our implementation, we generate simplices ``offline" first. This allows the workload to be distributed across several computers, since the validation step can be restricted to only a subset of the computed simplices. We do not implement a typical refinement procedure, where simplices that fail to validated are split, with more nodes added and corrected with Newton. Rather, we implement an \textit{adaptive refinement} step, which can help with the validation if failure is primarily a result of interval over-estimation. See Section \ref{sec-adaptive}.

\begin{remark}
The operators $A_s$ and $A_s^\dagger$ have standard interpretations in terms of the nonlinear map $F_s$. The operator $A_s^\dagger:X_1\rightarrow X_2$ is expected to be an approximation of $DF_s(\hat x_s)$, which means that $Z_1$ is a measure of the quality of the approximation. Conversely, $A_s:X_2\rightarrow X_1$ is expected to be an approximation of the \emph{inverse} of $A_s^\dagger$, so that $Z_0$ measures the quality of this approximation. Indirectly, $A_s$ acts as an appoximation of $DF_s(\hat x_s)^{-1}$.
\end{remark}

\subsection{Globalizing the manifold}\label{sec:globalize}
Theorem \ref{thm-rpa} guarantees that the map from the standard simplex to the zero set of $G_s(\cdot)$ is $C^1$. There is then a natural question as to the smoothness of the manifold obtained by gluing together the images of the $C^1$ maps. This is answered in the affirmative in \cite{Gameiro2016a}, and is primarily a consequence of the numerical data being equal on cobordant simplices. 

\section{Continuation of periodic orbits through (degenerate) Hopf bifurcations}\label{sec3}
In this section, we construct a nonlinear map whose zeroes will encode periodic solutions of a delay differential equation
\begin{align}\label{eq:DDE}\dot y(t)=f(y(t+\mu_1),\dots,y(t+\mu_J),\alpha,\beta)\end{align} for $f:(\R^n)^{J+1}\times\R^2\rightarrow\R^n$, with some (positive, negative or zero) constant delays $\mu_1,\dots,\mu_J$, and distinguished systems parameters $\alpha,\beta$. We will briefly consider ordinary differential equations in Section \ref{sec:ode} as a special case. We assume that $f$ is sufficiently smooth to permit further partial derivative computations. 

Following \cite{VandenBerg2021a}, we use the desingularization approach to isolate periodic orbits from (potentially) nearby steady states. This approach allows us to put a large distance (in the sense of a suitable Banach space) between steady states and periodic orbits that arise from Hopf bifurcations. This is exposited in Section \ref{sec-desingularization}, where we also discuss some details concerning non-polynomial nonlinearities. 

The next Section \ref{sec-zfp} is devoted to the development of a nonlinear map whose zeroes encode periodic orbits of our delay differential equation.  In this map, periodic orbits are isolated from fixed points. We present the map abstractly at the level of a function space, and then with respect to a more concrete sequence space.

In Section \ref{sec-zfp-continuation}, we lift the map of the previous section into the scope of two-parameter continuation. We develop an abstract template for the map on a relevant Banach space, define an approximate Fr\'echet derivative near a candidate zero of this map, and investigate some properties of the Newton-like operator. Specifically, we verify that numerical data corresponding to an approximate \emph{real} periodic orbit, under conditional contraction of the Newton-like operator, will converge to a real periodic orbit.
 
\subsection{Desingularization, polynomial embedding and phase isolation}\label{sec-desingularization}
We begin by doing a ``blowup" around the periodic orbit. Write $y = x+az$, for $x$ a candidate equilibrium point. Then we get the rescaled vector field  
\begin{align*}
\tilde f(z_1,\dots,z_J,x,a,\alpha,\beta)=\left\{\begin{array}{ll}
a^{-1}(f(x+a z_1,\dots,x+a z_J,\alpha,\beta)-f(x,\dots,x,\alpha,\beta)),&a\neq 0\\
\sum_{j=0}^J d_{y_j}f(x,\dots,x,\alpha,\beta)z_j,&a=0.\end{array}\right.
\end{align*}
The interpretation is that $||z||=1$, so that $a$ behaves like the relative norm-amplitude of the periodic orbit. The vector field above is $C^{k-1}$ provided the original function $f$ is $C^k$ and $x$ is an equilibrium point; that is, $f(x,\dots,x,\alpha,\beta)=0$. At this stage we can summarize by saying that our goal is to find a pair $(x,z)$ such
\begin{align*}
f(x,\dots,x,\alpha,\beta)&=0\\
\dot z(t) &= \tilde f(z(t+\mu_1),\dots,z(t+\mu_J),x,a,\alpha,\beta),\\
||z||&=1
\end{align*}
where $z$ is $\omega$-periodic for an unknown period $\omega$; equivalently, the frequency of $z$ is $\psi = \frac{2\pi}{\omega}$.
\begin{remark}
It is a common strategy in rigorous numerics, especially for periodic orbits, that time is re-scaled so that the period appears as a parameter in the differential equation. We do not do that here, since this would have the effect of dividing every delay $\mu_j$ by the period. This causes its own set of problems.
\end{remark}

In what follows, it will be beneficial for the vector field $\tilde f$ to be polynomial. This is because we will make use of a Fourier spectral method, and polynomial nonlinearities in the function space translate directly to convolution-type nonlinearities on the sequence space of Fourier coefficients. While we can make due with non-polynomial nonlinearities, it is greatly simplifies the computer-assisted proof if they are polynomial. To fix this, we generally advocate the use of the \emph{polynomial embedding} technique. The idea is that many analytic, non-polynomial functions are themselves solutions of polynomial ordinary differential equations. The reader may consult [van den Berg, Groothedde, Lessard] for a brief survey of this idea in the context of delay differential equations. See Section \ref{sec:SI_model} for a specific example.

Applying the polynomial embedding procedure always introduces additional scalar differential equations. If we need to introduce $m$ extra scalar equations to get a polynomial vector field, this will also introduce $m$ \emph{natural boundary conditions} that fix the initial conditions of the new components. As a consequence, we need to bring in $m$ \emph{unfolding parameters} to balance the system. This is accomplished in a problem-specific way; see [van den Berg, Groothedde, Lessard] for some general guidelines and a discussion for the need of these extra unfolding parameters. By an abuse of notation, we assume $\tilde f$ is polynomial (that is, the embedding has already been performed), and we write it as $$\tilde f(z_1,\dots,z_J,x,a,\alpha,\beta,\eta),$$ where $\eta\in\mathbb{R}^m$ is a vector representing the unfolding parameter, and we now interpret $\tilde f:(\R^{n+m})^{J+1}\times\R^2\rightarrow\R^{n+m}$. We then write the natural boundary condition corresponding to the polynomial embedding as $$\theta_{BC}(z(0),x,a,\alpha,\beta,\eta)=0\in\mathbb{R}^m.$$

It can also be useful to eliminate non-polynomial parameter dependence from the vector field, especially if the latter has high-order polynomial terms with respect to the state variable $z$. This can often be accomplished by introducing extra scalar variables. For example, if $\tilde f$ is $$\alpha e^{-px}z_1 - \beta z_2^2  + \eta_1$$
and we want to eliminate the non-polynomial term $e^{-px}$ from the vector field, then we can introduce a new variable $\eta_2$ and impose the equality constraint $0=\eta_2-e^{-px}$. The result is that the vector field becomes $$\alpha \eta_2z_1 - \beta z_2^2 + \eta_1.$$
Since this operation introduces new variables and additional constraints, we incorporate the constraints as extra components in the natural boundary condition function $\theta_{BC}$ of the polynomial embedding. Since this type of operation will introduce an equal number of additional scalar variables \textit{and} boundary conditions, we will neglect them from the dimension counting. 

\begin{remark}\label{rem-embedding-dimension}
If we want to formalize the embedding process for parameters, we can introduce differential equations for them. Indeed, in the example above, we have $\dot\eta_2=0$, and this differential equation can be added to the list of differential equations that result from polynomial embeddings of the original state variable $z$. In this way, we can understand $m$ as the total embedding dimension. It should be remarked, however, that in numerical implementation, objects like $\eta_2$ really are treated as scalar quantities. 
\end{remark}

The final thing we need to take into account is that every periodic orbit is equivalent to a one-dimensional continuum by way of phase shifts. Since our computer-assisted approach to proving periodic orbits is based on Newton's method and contraction maps, we need to handle this lack of isolation. This can be done by including a \emph{phase condition}. In this paper we will make use of an \emph{anchor condition}; we select a periodic function $\hat z$ having the same period as $z$, and we require that $$\int\langle z(s),\hat z'(s)\rangle ds = 0.$$

\subsection{Zero-finding problem}\label{sec-zfp}
We are ready to write down a zero-finding problem for our rescaled periodic orbits. First, combining the work of the previous sections, we must simultaneously solve the equations
\begin{equation}\label{eq:zfp-DDE}
\begin{cases}
\dot z = \tilde f(z(t+\mu_1),\dots,z(t+\mu_J),x,a,\alpha,\beta,\eta),&\quad\text{(delay differential equations)}\\
\| z\| = 1, &\quad\text{(amplitude condition of scaled orbit)}\\
\int\langle z(s),\hat z'(s)\rangle ds=0, &\quad\text{(anchor condition)}\\
f(x,\dots,x,\alpha,\beta)=0, &\quad\text{($x$ is a steady state)}\\
\theta_{BC}(z(0),x,a,\alpha,\beta,\eta)=0.&\quad\text{(embedding boundary condition)}
\end{cases}
\end{equation}
At this stage, the period $\omega$ of the periodic orbit is implicit. In passing to the spectral representation, we will make it explicit. Define the frequency $\psi = \frac{2\pi}{\omega}$ and expand $z$ in Fourier series:
\begin{align}
\label{eq:fourier_DDE}
z(t)&=\sum_{k\in\Z}z_k e^{ik\psi t}.
\end{align}
Recall that for a real (as opposed to complex-valued) periodic orbit, $z_k\in\C^{n+m}$  will satisfy the symmetry $z_k = \overline{z_{-k}}$. To substitute \eqref{eq:fourier_DDE} into the differential equation in \eqref{eq:zfp-DDE}, we must examine how time delays transform under Fourier series. Observe 
$$z(t+\mu) = \sum_{k\in\Z}z_ke^{ik\psi\mu}e^{ik\psi t},$$ which means that at the level of Fourier coefficients, a delay of $\mu$ corresponds to a complex rotation 
\begin{align}\label{eq:fourier_DDE_rotation}
z_k\mapsto (\zeta_\mu(\psi)z)_k\bydef e^{ik\psi\mu}z_k.
\end{align} 
Note that this operator is linear on $C^\Z_{n+m}$ and bounded on $\ell_\nu^1(\C^{n+m})$. To densify our notation somewhat, we define $\zeta(\psi):\ell_\nu^1(\C^{n+m})\rightarrow \ell_\nu^1(\C^{n+m})^J$ by
$$\zeta(\psi)z = (\zeta_{\mu_1}(\psi)z,\dots,\zeta_{\mu_J}(\psi)z).$$

As for the derivative, we make use of the operator $(Kz)_k = kz_k$. Since $\tilde f$ is polynomial, subsituting \eqref{eq:fourier_DDE} into the first equation of \eqref{eq:zfp-DDE} will result in an equation of the form
$$\psi iKz = \mathbf{f}(\zeta(\psi)z,x,a,\alpha,\beta,\eta)$$
for a function $\mathbf{f}:(\C^\Z_{n+m})^{J}\times\R^n\times\R\times\R^2\times\R^m\rightarrow (\C^\Z_{n+m})$ being a (formal) vector polynomial with respect to Fourier convolution, in the arguments $(\C^\Z_{n+m})^{J}$. Observe that we have abused notation and identified the function $z$ in \eqref{eq:fourier_DDE} with its sequence of Fourier coefficients. As an example, the nonlinearity $z\mapsto z(t)^2z(t+\mu)$ is transformed in Fourier to the nonlinearity $$z\mapsto (z*z)*(\zeta_\mu(\psi)z).$$

Now, let $\hat z$ be an approximate periodic orbit. We can define new amplitude and phase conditions as functions of $z$ and the numerical data $\hat z$; see \cite{VandenBerg2021a}. Then, define a map $G:X\times\R^2\rightarrow U$, with $X = \ell_\nu^1(\C^{n+m})\times\R^n\times\R\times\R\times\R^m$, and $U=K_\nu(\C^{n+m})\times\R^n\times\C\times\C\times\C^m$, by
\begin{align}\label{eq:F-DDE}
G(z,x,a,\psi,\eta,(\alpha,\beta))&=\left(\begin{array}{c}
-\psi iKz + \mathbf{f}(\zeta(\psi)z,x,a,\alpha,\beta,\eta) \\ f(x,\dots,x,\alpha,\beta) \\
\langle z,K^2\hat z\rangle - 1 \\
\langle z,iK\hat z\rangle \\ 
\theta_{BC}(\sum_k z_k,x,a,\alpha,\beta,\eta)
\end{array}\right).
\end{align}

Let $X$ be equipped with the norm \begin{align}\label{Xnorm}||(z,x,a,\psi,\eta)||=\max\{||z||_\nu,|x|,|a|,|\psi|,|\eta| \},\end{align} where all Euclidean space norms are selected a priori \emph{and could be distinct}. That is, we allow for the possibility of a refined weighting\footnote{It becomes notationally cubmersome to include references to weights, or to explicitly label the different norms, so we will refrain from doing so unless absolutely necessary. In these cases, footnotes will be used to emphasize any particular subtleties. Weighted max norms can make it easier to obtain computer assisted proofs, especially in continuation problem, which is why we have included this option.} of the norms being used. Then $X$ is a Banach space, and the same is true for $U$ when equipped with an analogous norm, replacing $||\cdot||_\nu$ with $||\cdot||_{\nu,K}$. If the vector field $f$ is twice continuously differentiable, then $F$ is continuously differentiable. This is because the polynomial embedding results in polynomial ordinary differential equations, while the blow-up procedure only results in one derivaive being lost. However, $F$ is generally not twice continuously differentiable; see later Section \ref{sec-partials}.

It will sometimes be convenient to compute norms on the $\R^{n+m+4}$-projection of $X$. In this case, if $u=(z,y)\in X$ and $y\in\R^{n+m+4}$ is represented (ismorphically) as $y=(x,a,\psi,\eta)\in\R^n\times\R\times\R\times\R^m$, then we define $||y||=\max\{|x|,|a|,|\psi|,|\eta|\}$, where any weighting is, again, implicit. Then $||(z,y)||= \max\{||z||_\nu,||y||\}$.

Introduce $V=\mbox{Symm}(\ell_\nu^1(\C^{n+m}))\times\R^n\times\R\times\R\times\R^m$. Any zero of $F$ in the space $V$ uniquely corresponds to a \emph{real} periodic orbit of \eqref{eq:DDE} by way of the Fourier expansion \eqref{eq:fourier_DDE} and the blow-up coordinates $y = x+az$. Moreover, the restriction of $F$ to $V$ has range in $W=\symm(K_\nu(\C^{n+m}))\times\R^n\times\R\times\R\times\R^m$. Each of $V$ and $W$ are Banach spaces over the reals, and so from this point on we work with the restriction $F:V\rightarrow W$.

\subsection{Finite-dimensional projection}\label{sec-zfp-continuation}
To set up the rigorous numerics and the continuation, we need to define projections of $V$ and $W$ onto suitable finite-dimensional vector spaces. Let $\pi^M:\ell_\nu^1(\C^{n+m})\rightarrow\ell_\nu^1(\C^{n+m})$ denote the projection operator $$(\pi^Mz)_k=\left\{\begin{array}{ll}z_k,&|k|\leq M \\ 0 & |k|>M, \end{array}\right.$$
and $\pi^\infty = I_{\ell_\nu^1(\C^{n+m})}-\pi^M$ its complementary projector. Consider the finite-dimensional vector space $$V^M\bydef\pi^M(\symm(\C^\Z_{n+m}))\times\R^n\times\R\times\R\times\R^m,$$ and extend the projection to a map $\pi^M:V\rightarrow V^M$ as follows: $$\pi^M(z,x,a,\psi,\eta) = (\pi^M z,x,a,\psi,\eta).$$ Now define a map $\G:V^M\times\R^2\rightarrow V^M$
\begin{align}
\G(z,x,a,\psi,\eta,(\alpha,\beta))&=\pi^M G(z,x,a,\psi,\eta,(\alpha,\beta)).\label{eq:FM_dde}
\end{align}
$\G$ is well-defined and smooth, and we have $\G = \pi^M G\circ i_{V^M}$, where $i_{V^M}:V^M\hookrightarrow V$ is the natural inclusion map. Therefore, this definition of the projection of $G$ is consistent with the abstract set-up of Section \ref{sec-continuation-infdim}.

\subsection{A reformulation of the continuation map}
It will be convenient to identify $V\times\R^2$ and $V^M\times\R^2$ respectively with isomorphic spaces
\begin{align*}
V\times\R^2&\sim \pi^M(\symm(\ell_\nu^1(\C^{n+m})))\times\R^{n+m+4}\times\pi^\infty(\symm(\ell_\nu^1(\C^{n+m})))\bydef\Omega\\
V^M\times\R^2&\sim \pi^M(\symm(\ell_\nu^1(\C^{n+m})))\times\R^{n+m+4}\bydef\Omega^M
\end{align*}
There is then the natural embedding $(z^M,\rho)\mapsto(z^M,\rho,0)$ of $\Omega^M$ into $\Omega$. The isomorphism of $V$ with $\Omega$ is given by $$(z,x,a,\psi,\eta,(\alpha,\beta))\mapsto(\pi^M z, (x,a,\psi,\eta,\alpha,\beta),\pi^\infty z).$$ $\Omega^M$ is finite-dimensional, and as such our language will sometimes reinforce this by describing matrices whose columns are elements of $\Omega^M$. This should be understood ``up to isomorphism". The purpose of the isomorphism of $V\times\R^2$ with $\Omega$ is to symbolically group all of the finite-dimensional objects together.

Given $\hat u_j = (\hat z_j,\hat \rho_j)\in \Omega^M$ for $j=0,1,2$, let $\hat \Phi_j$ be a matrix whose columns are a basis for the kernel of $DF^M(\hat z_j,\hat \rho_j)$, and are therefore elements of $V^M\times\R^2\sim\Omega^M$. For $s\in\Delta$, let $\hat u_s$ and $\hat\Phi_s$ be the usual interpolations of the elements $\hat u_j$ and bases $\hat \Phi_j$ for $j=0,1,2$. 

The continuation map $G_s$ of \eqref{Gmap_s} could now be defined for our periodic orbit function $G$. However, it will be convenient in our subsequent discussions concerning the radii polynomial approach to re-interpret the codomain of $G_s$ as being $$\tilde\Omega\bydef\pi^M\mbox{Symm}(K_\nu(\C^{n+m}))\times\R^{n+m+4}\times\pi^\infty \mbox{Symm}(K_\nu(\C^{n+m})).$$
Specifically, this will make it a bit easier to define an approximate inverse of $DG_s(\hat u_s)$. The codomain of $G_s$ is
\begin{align*}
W\times\R^2 &= (\mbox{Symm}(K_\nu(\C^{n+m}))\times\R^n\times\R\times\R\times\R^m)\times\R^2\\
&\sim \mbox{Symm}(K_\nu(\C^{n+m}))\times\R^{n+m+4}\\
&\sim \tilde\Omega,
\end{align*}
where the isomorphisms can be realized by permuting the relevant components of $G_s$ and splitting the Fourier space into direct sums. For $(z^M,\rho,z^\infty)\in\Omega$, a suitable isomorphic representation of $G_s$ is given by $\F_s:\Omega\rightarrow\tilde\Omega$,
\begin{align}
\label{Fmap}\F_s(z^M,\rho,z^\infty)&=\left(\begin{array}{c}
-\psi i K z^M + \pi^M\mathbf{f}(\zeta(\psi)(z^M+z^\infty),\rho) \\ 
\mathbf{J}_s(z^M,\rho,z^\infty)\\
-\psi i Kz^\infty + \pi^\infty\mathbf{f}(\zeta(\psi)(z^M+z^\infty),\rho)\end{array}\right)\bydef\left(\begin{array}{c}{\mathcal{F}}_s^{(1)} \\ {\mathcal{F}}_s^{(2)}\\ {\mathcal{F}}_s^{(3)}\end{array}\right),\\
\label{def-J}\mathbf{J}_s(z^M,\rho,z^\infty)&=\left(\begin{array}{c}f(x,\dots,x,\alpha,\beta) \\ \langle z^M,K^2\hat u_s\rangle - 1 \\ 
\langle z^M , iK \hat u_s\rangle\\
\theta_{BC}(\sum_k (z^M_k+z^\infty_k),x,a,\alpha,\beta,\eta) \\ \langle \hat\Phi_{s,1}, u^M\rangle - \hat c_{s,1}\\ \langle \hat\Phi_{s,2}, u^M\rangle- \hat c_{s,2}\end{array}\right), \hspace{2mm} u^M = \left(\begin{array}{c}z^M \\ \rho \end{array}\right)\in\Omega^M,
\end{align}
where $\rho = (x,a,\psi,\eta,\alpha,\beta)$, and the $\hat c_{s,j}$ for $j=1,2$ are interpolations
\begin{align}\label{csj}\hat c_{s,j}=\langle\hat\Phi_{0,j},\hat u_0\rangle + s_1(\langle\hat\Phi_{1,j},\hat u_1\rangle - \langle\hat\Phi_{0,j},\hat u_0\rangle) + s_2(\langle\hat\Phi_{2,j},\hat u_2\rangle - \langle\hat\Phi_{0,j},\hat u_0\rangle).\end{align}
\begin{remark}
Strictly speaking, the final two components of $\mathbf{J}_s$ are not the ``standard" ones from \eqref{Gmap_s}. The inner product $\langle \Phi_{s,j},\pi_\mathcal{X}x - \hat x_s\rangle$ in the latter results in quadratic terms in $s$, while the ones in \eqref{def-J} are $s$-linear. The impact of this change is that, theoretically, the Newton correction is not strictly in the direction orthogonal to the interpolation of the tangent planes. This change has no theoretical bearing on the validated continuation, and is done only foe ease of computation: it is easier to compute derivatives of linear functions than nonlinear ones.
\end{remark}

\begin{remark}
We have abused notation somewhat, since now we interpret $\mathbf{f}$ as a map $$\mathbf{f}:\symm(\ell_\nu^1(\C^{n+m}))\times\R^{n+m+4}\rightarrow\symm(K_\nu(\C^{n+m})).$$ It acts trivially with respect to the variable $\psi$. Also, we emphasize that $G_s$ depends on the numerical interpolants $\hat u_s$ and $\hat\Phi_s$.
\end{remark}

Since Theorem \ref{thm-rpa} is stated with respect to general Banach spaces, the validated continuation approach applies equally to the representation $\F_s:\Omega\rightarrow\tilde\Omega$ of $G_s:V\rightarrow W$. The only thing we need to do is specify a compatible norm on $\Omega$. This is straightforward: for $(z^M,\rho,z^\infty)\in\Omega$ and $\rho  = (x,a,\psi,\eta,\alpha,\beta)$, a suitable norm is $$||(z^M,\rho,z^\infty)||_\Omega = \max\{||z^M+z^\infty||_\nu,|x|,|a|,|\psi|,|\eta|,|\alpha|,|\beta|\},$$ where the Euclidean norms on the components of $\rho$ are the same\footnote{Including any weighting, which we recall is implicit.} as the ones appearing in \eqref{Xnorm}. With this choice, $||\cdot||_\Omega$ is equivalent to the induced max norm on $X\times\R^2$, with $X$ equipped with the norm \eqref{Xnorm} and $\R^2$ the $\infty$-norm.

\subsection{Construction of $A_s^\dagger$ and $A_s$}\label{sec-AsAsDagger}
Write $\hat u_s = (\hat z_s,\hat x_s,\hat \psi_s,\hat a_s,\hat\alpha_s,\hat\beta_s)$, for $s\in\Delta$. Denote the three vertices of $\Delta$ as $s_0 = (0,0)$, $s_1 = (1,0)$ and $s_2 = (0,1)$.  Introduce an approximation of $\mathcal{F}_s$ as follows:
\begin{align}
\tilde{\mathcal{F}}_s(z^M,\rho,z^\infty)&=\left(\begin{array}{c}
-\psi i Kz^M + \pi^M\mathbf{f}(\zeta(\psi)z^M,\rho) \\ 
\mathbf{J}_s(z^M,\rho,z^\infty)\\ 
-\psi iKz^\infty \end{array}\right)\bydef\left(\begin{array}{c}\tilde{\mathcal{F}}_s^{(1)} \\ \tilde{\mathcal{F}}_s^{(2)}\\\tilde{\mathcal{F}}_s^{(3)}\end{array}\right).
\end{align}
Formally, $\tilde\F_s$ approximates $\F_s$ in the Fourier tail by neglecting the nonlinear terms, leaving only the part coming from the differentiation operator. 

\begin{proposition}
$D\tilde{\mathcal{F}}_s(\hat u_s)$ has the representation
\begin{align}
D\tilde{\mathcal{F}}_s(\hat u_s)&=\left(\begin{array}{ccc}
D_1\tilde{\mathcal{F}}_s^{(1)}(\hat u_s)&D_2\tilde{\mathcal{F}}_s^{(1)}(\hat u_s)&0 \\ D_1\tilde{\mathcal{F}}_s^{(2)}(\hat u_s) & D_2\tilde{\mathcal{F}}_s^{(2)}(\hat u_s) & D_3\tilde{\mathcal{F}}_s^{(2)}(\hat u_s) \\ 
0&0&D_3\tilde{\mathcal{F}}_s^{(3)}(\hat u_s)
\end{array}\right)\bydef A_s^\dagger.
\end{align}
\end{proposition}

\begin{proof}
Since $\F_s^{(1)}$ does not depend on $z^\infty$, the upper-right block of $D\F_S(\hat u_s)$ is the zero map. Similarly, $\F_s^{(3)}$ does not depend on either $z^M$ or $\rho$, so the finite-dimensional blocks in the $z^\infty$ (bottom) row are zero.
\end{proof}
The upper left $2\times 2$ block is equivalent to a finite-dimensional matrix operator. In particular, $D_2\tilde{\mathcal{F}}_s^{(2)}(\hat u_s)$ is real. Also,  $$D_3\tilde{\mathcal{F}}_s^{(3)}(\hat u_s) = -\hat\psi_s iK\pi^\infty$$ is invertible, with
$$D_3\tilde{\mathcal{F}}_s^{(3)}(\hat u_s)^{-1}=i(\hat\psi_s)^{-1}(K\pi^\infty)^{-1},$$ where $((K\pi^\infty)^{-1}z)_k \bydef \frac{1}{k}z_k$ for $|k|\geq M+1$. Suppose we can explicitly compute $S_j\in\R^{(m+n+4)\times(m+n+4)}$ and $P_j,Q_j,R_j$ complex matrices for $j=0,1,2$ such that
\begin{align}\label{eq:DDE-approximate-inverse}
\left(\begin{array}{cc}P_j&Q_j\\R_j&S_j \end{array}\right)\left(\begin{array}{cc}
D_1\tilde{\mathcal{F}}_{s_j}^{(1)}(\hat u_{s_j})&D_2\tilde{\mathcal{F}}_{s_j}^{(1)}(\hat u_{s_j})\\
D_1\tilde{\mathcal{F}}_{s_j}^{(2)}(\hat u_{s_j}) & D_2\tilde{\mathcal{F}}_{s_j}^{(2)}(\hat u_{s_j})
\end{array}\right)\approx I_{\Omega^M}.
\end{align}
We can then prove the following lemma.
\begin{lemma}\label{lem:dde_symmetry}
For $s=(s^{(1)},s^{(2)})\in \Delta$, define matrix interpolants $P_s=P_1 + s^{(1)}(P_2-P_1) + s^{(2)}(P_3-P_1)$, and analogously define interpolants $Q_s$, $R_s$ and $S_s$. Introduce a family of operators $A_s$ as follows:
\begin{align}\label{As}
A_s&=\left(\begin{array}{ccc}
P_s&Q_s&-Q_sD_3\tilde{\mathcal{F}}_s^{(2)}(\hat u_s)i(\hat\psi_s)^{-1}(K\pi^\infty)^{-1} \\ 
R_s&S_s&-S_sD_3\tilde{\mathcal{F}}_s^{(2)}(\hat u_s)i(\hat\psi_s)^{-1}(K\pi^\infty)^{-1} \\ 
0&0&i(\hat\psi_s)^{-1}(K\pi^\infty)^{-1}
\end{array}\right).
\end{align}
Suppose for $j=1,2,3$, $S_j$ is real and, as maps, $P_j$, $Q_j$ and $R_j$ are equivalent to
\begin{align*}
P_j&:\pi^M(\symm(\C^\Z_{n+m}))\rightarrow\pi^M(\symm(\C^\Z_{n+m}))\\
Q_j&:\R^{m+n+4}\rightarrow\pi^M(\symm(\C^\Z_{n+m}))\\
R_j&:\pi^M(\symm(\C^\Z_{n+m}))\rightarrow\R^{n+m+4}.
\end{align*}
Then $A_s:\tilde\Omega\rightarrow\Omega$ is well-defined. 
\end{lemma}
\begin{proof}
One can show $A_s:\tilde\Omega\rightarrow\Omega$ is well-defined using the fact that each of $P_s$, $Q_s$ and $R_s$ is a real convex combination of maps to/from appropriate symmetric spaces, and noticing that $iK\pi^\infty$ (and hence its inverse) satisfy the symmetry $(iK\pi^\infty z)_k = \overline{(iK\pi^\infty z)_{-k}}$ 
\end{proof}
The point here is that, due to \eqref{eq:DDE-approximate-inverse}, we have $$A_{s_j}\approx D\tilde{\mathcal{F}}_{s_j}(\hat u_{s_j})^{-1} \approx (A^\dagger_{s_j})^{-1}$$ for $j=0,1,2$, and if the interpolation points $\hat u_j$ are close together, we expect $D\tilde{\mathcal{F}}_s(\hat u_{s_j})$ to be invertible for all $s\in\Delta$, and $D\tilde{\mathcal{F}}_s(\hat u_{s_j})^{-1}\approx A_s$. 
\begin{remark}
Checking the conditions of Lemma \ref{lem:dde_symmetry} amounts to verifying conjugate symmetries of the matrices $P_j$, $Q_j$ and $R_j$. Numerical rounding makes this a nontrivial task, so we generally post-process the numerically computed matrices to impose these symmetry conditions.
\end{remark}

\section{Technical bounds for validated continuation of periodic orbits}\label{sec:CAP}
Based on the previous section, we define a Newton-like operator $T_s:\Omega\rightarrow\Omega$
\begin{align}\label{eq:newton-like-dde}
T_s(u) = u - A_s\mathcal{F}_s(u),
\end{align}
for $s\in\Delta$. As described in Section \ref{sec-continuation-proofs}, our goal is to prove that $T_s$ is a uniform (for $s\in\Delta$) contraction in a closed ball centered at $\hat u_s$ using Theorem \ref{thm-rpa}. If $A_s$ can be proven (uniformly in $s$) injective, this will prove the existence of a unique zero of $\mathcal{F}_s(\cdot)$ close to $\hat u_s$, thereby validating the simplex and proving the smooth patch of our manifold.

In Section \ref{sec:injective}, we will demonstrate how the bound $Z_0$ of the radii polynomial approach can be used to obtain a proof of uniform (in $s$) injectivity of the operator $A_s$. We then provide some detailed discussion concerning general-purpose implementation of the bounds $Y$ and $Z$ in Section \ref{sec-Y0} through to Section \ref{sec-Z2}.

\subsection{Injectivity of $A_s$}\label{sec:injective}
The injectivity of $A_s$ is a consequence of the successful identification of bounds $Z_0$ of Theorem \ref{thm-rpa} and the negativity of the radii polynomial. In particular,
\begin{lemma}\label{lemma-injectivity}
Suppose $||I-A_sA_s^\dagger||_{B(\Omega,\Omega)}\leq Z_0$ for all $s\in\Delta$, with the operators $A_s$ and $A_s^\dagger$ of Section \ref{sec-AsAsDagger}. If $Z_0<1$, then $A_s$ is injective for $s\in\Delta$.
\end{lemma}
\begin{proof}
First, observe $A_s$ has non-trivial kernel if and only if there exists $u\in\Omega^M$ such that $A_s u = 0$. This is a consequence of the structure of the operator and injectivity of $(K\pi^\infty)^{-1}$. Therefore, it suffices to verify the injectivity of the restriction $\mathbf{A}_s = A_s|_{\Omega^M}$. Define $\mathbf{A}_s^\dagger = A_s^\dagger|_{\Omega^M}$. By definition of the norm on $\Omega$, we have $||I-\mathbf{A}_s\mathbf{A}_s^\dagger||_{B(\Omega^M,\Omega^M)}\leq Z_0<1$ for all $s\in\Delta$. By Neumann series, it follows that $\mathbf{A}_s\mathbf{A}_s^\dagger$ is boundedly invertible, which implies $\mathbf{A}_s$ is surjective. Since $\Omega^M$ is finite-dimensional, $\mathbf{A}_s$ is also injective.
\end{proof}

\begin{corollary}
If the radii polynomial satisfies $p(r_0)<0$ for some $r_0>0$, then $A_s$ is injective for $s\in\Delta$.
\end{corollary}

\subsection{The bound $Y_0$}\label{sec-Y0}
To begin, expand the product $A_s\F_s(\hat u_s)$. We get
\begin{align*}
A_s\F_s(\hat u_s)&=\left(\begin{array}{c} P_s\F_s^{(1)}(\hat u_s) + Q_s\F_s^{(2)}(\hat u_s) - Q_sD_3\tilde\F_s^{(2)}(\hat u_s)i(\hat\psi_s)^{-1}(K\pi^\infty)^{-1}\F_s^{(3)}(\hat u_s) \\ R_s\F_s^{(1)}(\hat u_s) + S_s\F_s^{(2)}(\hat u_s) - S_sD_3\tilde\F_s^{(2)}(\hat u_s)i(\hat\psi_s)^{-1}(K\pi^\infty)^{-1}\F_s^{(3)}(\hat u_s)\\ i(\hat\psi_s)^{-1}(K\pi^\infty)^{-1}\F_s^{(3)}(\hat u_s) \\  \end{array}\right)
\end{align*}
Remark that $\F_s^{(3)}(\hat u_s)$ has range in a finite-dimensional subspace of $\Omega$; specifically, it will be in the part of $\Omega$ such that the components in $\C^\Z_{n+m}$ with index (in absolute value) greater than $Md+1$ are zero, where $d$ is the maximum degree of the (convolution) polynomial $\mathbf{f}$. As such, $A_s\F_s(\hat u_s)$ is explicitly computable. 


In practice, we must compute an enclosure of the norm $||A_s\mathcal{F}_s(\hat u_s)||$ for all $s\in\Delta$. This is slightly less trivial. We accomplish this using a first-order Taylor expansion with remainder. For the function $(s,u)\mapsto\mathcal{F}_s(u)$, denote by $\partial_s$ the Fr\'echet derivative with respect to $s$, and $D$ the derivative with respect to $u$. Given the interpolants $\hat u_s=(\hat v_s,\hat\rho_s)$ and $\hat\Phi_{s,j}$, let $\hat u'$ and $\hat\Phi_{j}'$ denote their Fr\'echet derivatives with respect to $s$. Also, denote $\hat c_{s,j}'= \partial_s \hat c_{s,j}$, for $\hat c_{s,j}$ defined in \eqref{csj}. Note that these derivatives are constant, since the interpolants are linear in $s$. Then
\begin{align*}
\F_s(\hat u_s)&= \F_{s_0}(\hat u_0) + \big(\partial_s\F_{s_0}(\hat u_0) + D\F_{s_0}(\hat u_0)\hat u'\big)s + \frac{1}{2}\mathcal{R},
\end{align*}
where the remainder term $\mathcal{R}$ will be elaborated upon momentarily. The product of the linear-order terms with $A_s$ results in a function that is componentwise quadratic in $s\in\Delta$, for which we can efficiently compute an upper bound on the norm. The derivative $D\F_s(\hat u_s)$ is implementable, and will be further discussed in Section \ref{sec-partials}. As for the $\partial_s$ term, 
\begin{align*}
\partial_s\F_s(u)&=\left(\begin{array}{c}0\\ \partial_s\mathbf{J}_s(u)\\0\end{array}\right),\hspace{2mm}\partial_s\mathbf{J}_s(u)=\left(\begin{array}{c}0 \\ \langle z,K^2\hat z'\rangle \\ \langle z,iK\hat z'\rangle \\ 0 \\ \langle\hat\Phi_1',u\rangle - \hat c_{s,1}'\\ \langle\hat\Phi_2',u\rangle - \hat c_{s,2}' \end{array}\right),
\end{align*}
where $u=(z,\rho)$ and $\hat u_s = (\hat z_s,\hat\rho_s)$.

As for the remainder $\mathcal{R}$, it is bounded by the norm of the second Fr\'echet derivative of $s\mapsto\F_s(\hat u_s)$, uniformly over $s\in\Delta$. For each $s\in\Delta$, let this second derivative be denoted $\mathbf{D}^2\F_s(\hat u_s)$. Then
\begin{align}
\label{boldD2F}\mathbf{D}^2\F_s(\hat u_s)&=\partial_s^2\F_s(\hat u_s)[e_1,e_2] + \partial_s D\F_s(\hat u_s)[\hat u'e_1,e_2] + D\partial_s\F_s(\hat u_s)[e_1,\hat u' e_2]\\
\nonumber&\phantom{=}+D^2\F_s(\hat u_s)[\hat u'e_1,\hat u'e_2].
\end{align}
where $e_1$ and $e_2$ denote the first and second (factor) projection maps on $\R^2\times\R^2$. The derivative $D^2\F_s$ will be discussed in Section \ref{sec-partials}. 
At this stage, we need only mention that it acts bilinearly on $\hat u$, and the latter is proportional in norm to the step size $\sigma$ of the continuation scheme, so the $D^2\F_s$ term will be order $\sigma^2$. As for the derivatives involving $\partial_s$, most of the components are zero as evidenced by the previous calculation of $\partial_s\F_s(z,\rho)$, and it suffices to compute the relevant derivatives of $\mathbf{J}_s$. We have
\begin{align*}
\partial_s^2\mathbf{J}_s(\hat u_s)[t_1,t_2]&=0 
,\hspace{2mm}D\partial_s\mathbf{J}_s(\hat u_s)[t,h]=\left(\begin{array}{c}0 \\ \langle h_z,K^2\hat z't\rangle \\ \langle h_z,iK\hat z't\rangle \\ 0 \\ \langle\hat\Phi_1't,h\rangle \\ \langle\hat\Phi_2't,h\rangle \end{array}\right),
\end{align*}
with $(t,h)\in\R^2\times\Omega$, $h = (h_z,h_\rho)$. In terms of the step size, $\partial_s^2\mathbf{J}_s(\hat u_s)$ is order $\sigma^2$, while $D\partial_s\F_s(\hat u_s)$ is order $\sigma$. However, in \eqref{boldD2F}, the latter term is multiplied by $\hat z' = O(\sigma)$. Therefore, as expected, the remainder $\mathcal{R}$ is quadratic with respect to step size. We therefore compute
$$Y_0\geq ||A_s(\F_s(\hat u_s) - 2^{-1}\mathcal{R})|| + \frac{1}{2}||A_s\mathcal{R}||.$$
Since $||R||=O(\sigma^2)$, the bound can be tempered quadratically by reducing the step size. The caveat is that if $||A_s||$ and/or $\F_s$ has large quadratic terms, it might still be necessary to take small steps. 
\begin{remark}
Directly computing the norm $||A_s\mathcal{R}||$ would require a general-purpose implementation of the second derivative $D^2\F_s(\hat u_s)$; see \eqref{boldD2F}. As we have stated previously, such an implementation would be rather complicated. Therefore, in practice, we perform another level of splitting; namely, we use the bound $$||A_s \mathcal{R}||\leq ||A_s\left(\mathbf{D}^2\F_s(\hat u_s) - D^2\F_s(\hat u_s)[\hat u'e_1,\hat u'e_2] \right)|| + ||A_sD^2\F_s(\hat u_s)[\hat u'e_1,\hat u'e_2]||.$$
The first term on the right of the inequality is explicitly implementable using only the derivatives of $\mathbf{J}_s$ and the finite blocks of $A_s$. For the second term, we use the bound $||A_s||\cdot||D^2\F_s(\hat u_s)[\hat u'e_1,\hat u'e_2]||$ which, while not optimal, is implementable and good enough for our purposes.

\subsubsection{Adaptive refinement}\label{sec-adaptive}
In continuation, the size of the $Y_0$ bound is severely limited by the step size. To distribute computations, we often want to compute the manifold first and then validate patches \textit{a posteriori}. However, once the manifold has been computed, adjusting and re-computing patches of the manifold with smaller step sizes becomes complicated due to the need to ensure that cobordant simplices have matched data, as discussed in Section \ref{sec:globalize}. When the $Y_0$ bound is too large due to interval arithmetic over-estimation, we can circumvent this by using \textit{adaptive refinement} on the relevant simplex. Formally, we subdivide the simplex into four, using the interpolated zeroes at the nodes of the original simplex to define the data at the nodes of the four new ones. The result is that cobordant data still matches, allowing for globalization of the manifold. The advantage of this approach is that it can be safely done in a distributed manner; adaptively refining one simplex does not require re-validating any of its neighbours. 

\end{remark}
\subsubsection{The bound $Z_0$}
The product $A_sA_s^\dagger$ is block diagonal. Indeed, $(I-A_sA_s^\dagger)|_{\pi^\infty\Omega} = 0$, whereas
\begin{align*}
(I-A_sA_s^\dagger)|_{\Omega^M} = I_{\Omega^M} - \left(\begin{array}{cc}P_s&Q_s\\R_s&S_s \end{array}\right)\left(\begin{array}{cc}
D_1\tilde{\mathcal{F}}_{s}^{(1)}(\hat u_{s})&D_2\tilde{\mathcal{F}}_{s}^{(1)}(\hat u_{s})\\
D_1\tilde{\mathcal{F}}_{s}^{(2)}(\hat u_{s}) & D_2\tilde{\mathcal{F}}_{s}^{(2)}(\hat u_{s})
\end{array}\right).
\end{align*}
Therefore, to compute $Z_0$ it suffices to find an upper bound, uniformly in $s$, for the norm of the above expression as a linear map from $\Omega^M\rightarrow\Omega^M$. Interpreted as matrix operators, $P_s$, $Q_s$, $R_s$ and $S_s$ are interpolants of other explicit matrices. However, the derivatives $D_i\tilde\F_s^{(j)}(\hat u_s)$, while evaluated at interpolants, are themselves ``nonlinear" in $s$.
\begin{remark}\label{rem:mesh}
The implementation of $||(I-A_sA_s^\dagger)|_{\Omega^M}||$ is influenced by the way in which the dependence on $s$ is handled. For example, $s$ can be treated as a vector interval and the norm can be computed ``in one step" using interval arithmeic, then we take the interval supremum to get a bound. This can result in some wrapping (over-estimation). One way to control the wrapping is to cover $\Delta$ in a mesh of balls, compute the norm $||(I-A_sA_s^\dagger)|_{\Omega^M}||$ for $s$ replaced with each of these interval balls, and take the maximum. Still another way is to carefully compute a Taylor expansion with respect to $s$, although this task has a few technical issues due to the fact that $\F_s(\cdot)$ is generally only $C^1$. We therefore only consider the ``in one step" 
approach in our implementation.
\end{remark}


\subsection{A detour: derivatives of convolution-type delay polynomials}\label{sec-partials}
In our implementation of the $Z_1$ and $Z_2$ bounds, we will need to represent various partial derivatives of the map $(z,\rho)\mapsto\mathbf{f}(\zeta(\psi)z,\rho)$, where $\rho = (x,a,\psi,\eta,\alpha,\beta)$. This can quickly become notationally cumbersome. Therefore, in this section we will elaborate on the structure of the derivatives of convolution terms of the following $K_\nu(\C)$-valued map:
\begin{align}
\Theta(z,\rho)\bydef\rho^{\bm}\prod_{p=1}^d\left(e^{i\psi(\rho)\mu_{j_p}K}z_{c_p}\right)\label{poly}
\end{align}
for $z\in\ell_\nu^1(\C^{n+m})$ and $\rho\in\R^{m+n+4}$. That is $\Theta:\ell_\nu^1(\C^{n+m})\times\R^{n+m+4}\rightarrow K_\nu(\C)$. Here, $\psi(\rho)$ denotes the frequency component of $\rho$. In this section, $z$ will be indexed with the convention $z=(z_1,\dots,z_{n+m})$, where $z_q\in\ell_\nu^1(\C)$ for $q=1,\dots,n+m$. 
The objects \eqref{poly} can be interpreted as individual terms of $\mathbf{f}_1$ through $\mathbf{f}_{n+m}$, which each have codomain $K_\nu(\C)$. The product symbol indicates iterated convolution, while $(e^{i\psi\mu K}z)_k\bydef e^{i\psi\mu k}z_k$. Here, $d\in\mathbb{N}$ is the (polynomial power) order of the term, the indices $c_p\in\{1,\dots,n+m\}$ specify which factors of $\ell_\nu^1(\C^{n+m})$ are involved in the multiplication, while $j_p\in\{1,\dots,J\}$ indicates which delays are associated to each of them. Finally, there is a multi-index $\bm$ for multi-index power $\rho^\bm$. Importantly, the multi-index $\bm$ is trivial in the frequency ($\psi$) component, and the latter only enters $\mathbf{f}$ in the form of the delay mapping $\zeta(\psi)$. The following lemma can then be proven by means of some rather tedious bookkeeping and the commutativity of the convolution product.

\begin{lemma}\label{lem-2nd-derivatives}
If $q\in\{1,\dots,n+m\}$, then for $z\in\ell_\nu^1(\C^{n+m})$ and $h\in\ell_\nu^1(\C)$,\small
\begin{align*}
\frac{d}{dz_q}\Theta(z,\rho)h&=\rho^\bm\sum_{\substack{r=1 \\ c_r=q}}^d(e^{i\psi(\rho)\mu_{j_r}K}h)*\prod_{\substack{p=1 \\ p\neq r}}^d e^{i\psi(\rho)\mu_{j_p}K}z_{c_p}\\
\frac{d}{d\psi(\rho)}\Theta(z,\rho)&=\rho^\bm\sum_{r=1}^d\left(\prod_{\substack{p=1 \\ p\neq r}}^d e^{i\psi(\rho)\mu_{j_p}K}z_{c_p}\right)*\left(iK\mu_{j_r}e^{i\psi(\rho)\mu_{j_r}K}z_{c_r}\right)\\
\frac{d}{dz_q}\left[\frac{d}{d\psi(\rho)}\Theta(z,\rho)\right]h&=\rho^\bm\sum_{\substack{r=1 \\ c_r=q}}^d(iK\mu_{j_r}e^{i\psi(\rho)\mu_{j_r}K}h)*\prod_{\substack{p=1 \\ p\neq r}}^d e^{i\psi(\rho)\mu_{j_p}K}z_{c_p}\\
&\phantom{=}+\rho^\bm\sum_{\substack{r=1 \\ c_r=q}}^d(e^{i\psi(\rho)\mu_{j_r}K}h)*\sum_{\substack{p=1 \\ p\neq r}}^d\left(\prod_{\substack{\xi=1 \\ \xi\neq r,p}}^de^{i\psi(\rho)\mu_{j_\xi}K}z_{c_\xi}\right)* (iK\mu_{j_p}e^{i\psi(\rho)\mu_{j_p}K}z_{c_p})
\end{align*}
\normalsize Also, $\frac{d}{d\psi(\rho)}\left[\frac{d}{dz_q}\Theta(z,p)h\right] = \frac{d}{dz_q}\left[\frac{d}{d\psi(\rho)}\Theta(z,\rho)\right]h$. If $z\in\C_\Z^{n+m}$ is band-limited, $\frac{d^2}{d\psi(\rho)^2}\Theta(z,\rho)$ exists and
\small
\begin{align*}
\frac{d^2}{d\psi(\rho)^2}\Theta(z,\rho)&=\rho^\bm\sum_{r=1}^d\left(\sum_{\substack{p=1\\p\neq r}}^d\left(\prod_{\substack{q=1 \\ q\neq p,r}}^d e^{i\psi(\rho)\mu_{j_q}K}z_{c_q} \right)*(iK\mu_{j_p}e^{i\psi(\rho)\mu_{j_p}K}z_{c_p})*(iK\mu_{j_r}e^{i\psi(\rho)\mu_{j_r}K}z_{c_r}) \right)\\
&\phantom{=}+\rho^\bm\sum_{r=1}^d\left(\prod_{\substack{p=1 \\ p\neq r}}^d e^{i\psi(\rho)\mu_{j_p}K}z_{c_p}\right)*\left(-K^2\mu_{j_r}^2e^{i\psi(\rho)\mu_{j_r}K}z_{c_r}\right).
\end{align*}
\normalsize
Finally, if $q_1,q_2\in\{1,\dots,n+m\}$ and $h_1,h_2\in\ell_\nu^1(\C)$, then
\begin{align*}
\frac{d^2}{dz_{q_2}dz_{q_1}}\Theta(z,p)[h_1,h_2]&=\rho^\bm\sum_{\substack{r=1 \\ c_r=q_1}}^d(e^{i\psi(\rho)\mu_{j_r}K}h_1)*\left(\sum_{\substack{p=1 \\ c_p = q_2 \\ p\neq r}}^d(e^{i\psi(\rho)\mu_{j_p}K}h_2)*\prod_{\substack{\xi=1 \\ \xi\neq r,p}}^de^{i\psi(\rho)\mu_{j_\xi}K}z_{c_\xi}\right)
\end{align*}
\end{lemma}
\begin{remark}\label{remark-D2F-DNE}
The requirement that $z$ be band-limited really is necessary for the existence of $\frac{d^2}{d\psi(\rho)^2}\Theta(z,p)$. Indeed, if $z\in\ell_\nu^1(\C^{n+m})$, then $Kz_j\in K_\nu(\C)$ but $K^2 z$ is not, and the latter terms appear in the second derivative. This is precisely why $\F_s$ is \emph{not twice differentiable}.
\end{remark}
\subsection{The bound $Z_1$}
To have a hope at deriving a $Z_1$ bound, we will first determine the structure of the operator $D\F_s(\hat u_s) - A_s^\dagger$. Represented as an ``infinite block matrix", most blocks are zero. One can verify that
\begin{align}\label{Z1-difference}
D\F_s(\hat u_s)-A_s^\dagger&=\left(\begin{array}{ccc}0&0&\bZ_1^{(1,3)} \\ 0&0&0\\ \bZ_1^{(3,1)}&\bZ_1^{(3,2)}&\bZ_1^{(3,3)} \end{array}\right),
\end{align}
with the individual terms $\bZ_1^{(i,j)}$ being the operators
\begin{align*}
\bZ_1^{(1,3)}&=\pi^MD_1\mathbf{f}(\zeta(\hat\psi_s)\hat z_s,\hat\rho_s)\zeta(\hat\psi_s)\pi^\infty\\
\bZ_1^{(3,1)}&=\pi^\infty D_1\mathbf{f}(\zeta(\hat\psi_s)\hat z_s,\hat\rho_s)\zeta(\hat\psi_s)\pi^M\\
\bZ_1^{(3,2)}&=\pi^\infty D_1\mathbf{f}(\zeta(\hat\psi_s)\hat z_s,\hat\rho_s)\zeta'(\hat\psi_s)\hat z_s\psi(\cdot) + \pi^\infty D_2\mathbf{f}(\zeta(\hat\psi_s)\hat z_s,\hat\rho_s)\\
\bZ_1^{(3,3)}&=\pi^\infty D_1\mathbf{f}(\zeta(\hat\psi_s)\hat z_s,\hat\rho_s)\zeta(\hat\psi_s)\pi^\infty,
\end{align*}
and $\psi:\R^{n+m+4}\rightarrow\R$ once again denoting the frequency component projection. Computing $Z_1$ requires precomposing \eqref{Z1-difference} with $A_s$, which has the structure \eqref{As}. 

\subsubsection{Virtual padding}
For $h\in \Omega$ with $||h||_\Omega\leq 1$, denote 
\begin{align}\label{g-Z1}g = A_s(D\mathcal{F}_s(\hat u_s) - A_s^\dagger)h,\hspace{2mm}g=(g^M,g^\rho,g^\infty).\end{align}
Computing $Z_1$ is equivalent to a uniform (in $h$) bound for $||g||_\Omega$. The tightness of this bound is determined by two levels of computation:
\begin{itemize}
\item some finite norm computations that are done on the computer;
\item theoretical bounds, which are inversely proportional to the dimension of the object on which the finite norm computations are done.
\end{itemize}
By default, the size of the finite computation is linear in $M$, the number of modes. This might suggest that explicitly increasing $M$ -- that is, padding our solution with extra zeros and re-computing everything with more modes -- is the only way to improve the bounds. Thankfully, this is not the case. 

\begin{figure}
\centering\includegraphics[scale=0.55]{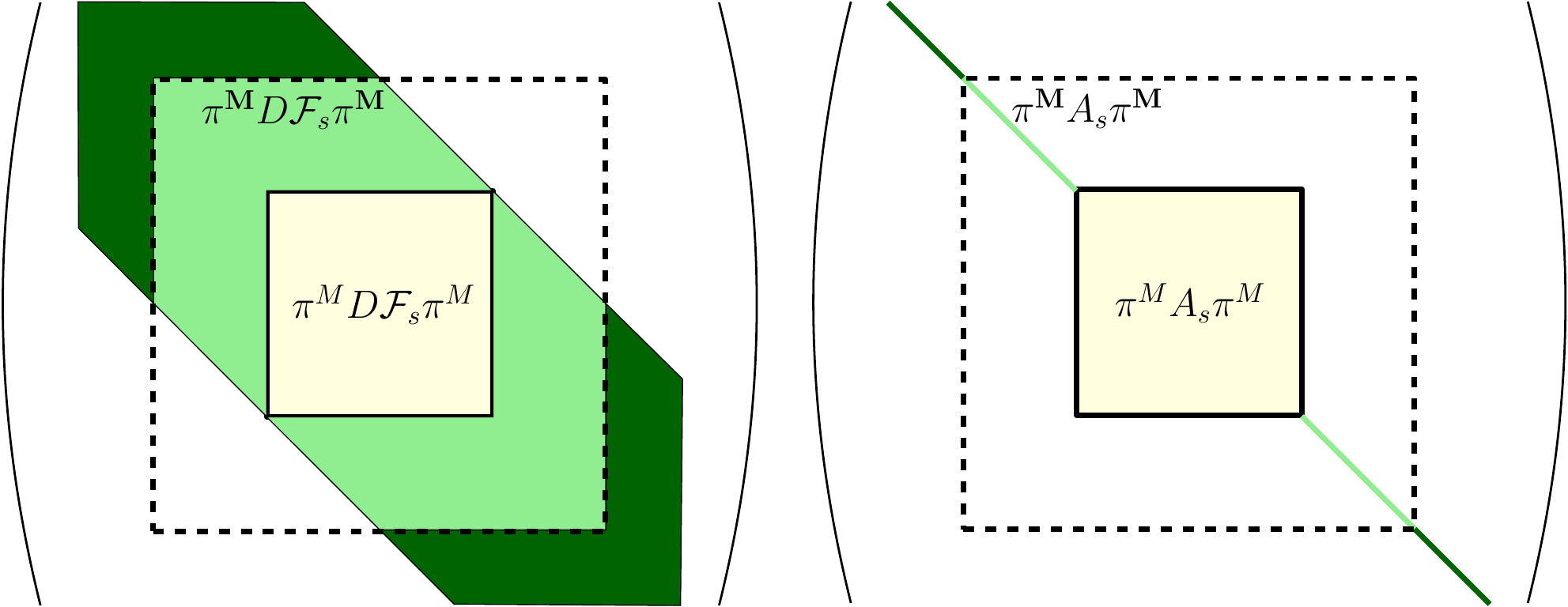}
\caption{Visual depiction of the virtual padding process as it applies to the operators $D\mathcal{F}_s$ and $A_s$, left and right. When we do not do virtual padding ($\mathbf{M}=M$), the size of the objects we store in memory on the computer for matrix operations are strictly determined by the dimension of the numerical zero. These objects are the inner boxes in the infinite matrices above, and are depicted in pale yellow. The part outside of this box is controlled analytically. When we use nontrivial virtual padding $(\mathbf{M}>M$), we store a larger amount of information in the computer for matrix operations; this is the outer dashed line box. Once again, the region outside the box is controlled analytically. Regions in white represent zeroes, while dark green represents the infinite part of the matrix, and light green a finite part of the matrix that is analytically controlled if there is no virtual padding.}
\label{fig-numerical-center}
\end{figure}

Intuitively, $A_s(D\F_s(\hat u_s-A_s^\dagger)$ is an ``infinite matrix", for which we have a canonical numerical center determined by the pre- and post- composition with the projection operator onto $\pi^M\symm(\ell_\nu^1(\C^{n+m}))\times\R^{n+m+4}$. This is determined by the number of modes $M$ we have specified in our numerical zero. However, there is no reason to only compute this much of the infinite matrix explicitly; we could instead choose $\mathbf{M}\geq M$ and compute the representation of this operator on $\pi^\pad\symm(\ell_\nu^1(\C^{n+m}))\times\R^{n+m+4}$. The result is that a larger portion of $A_s(D\F_s(\hat u_s-A_s^\dagger)$ is stored in the computer's memory. The advantage of doing this is that the explicit computations of norms are generally much tighter than theoretically guaranteed estimates, while the theoretical bounds related to the tail will be proportional to $\frac{1}{\mathbf{M}+1}$. See Figure \ref{fig-numerical-center} for a visual representation. 

To exploit these observations, we decompose $\ell_\nu^1(\C^{n+m})$ (and hence the symmetric subspace) as a direct sum $$\pi^\pad(\ell_\nu^1(\C^{n+m}))\oplus\pi^{\pad+1,\infty}(\ell_\nu^1(\C^{n+m})),$$
where $\pad\geq M$, and $\pi^{\pad+1,\infty}$ is now interpreted as the complementary projector to $\pi^\pad$. In the subsections that follow, $\pad$ will be the \textit{virtual padding dimension}. We therefore re-interpret $h=(h^\pad,h^\rho,h^\infty)$ and $g=(g^\pad,g^\rho,g^\infty)$ as being in the product space $$\pi^{\pad}\symm(\ell_\nu^1(\C^{n+m}))\times\R^{n+m+4}\times\pi^{\pad+1,\infty}\symm(\ell_\nu^1(\C^{n+m})),$$ which remains (isometrically) isomorphic to $\Omega$. Remark that the finite blocks $P_s,Q_s$ and $R_s$ coming from $A_s$ must now be re-interpreted as linear maps involving $\pi^\pad\symm(\ell_\nu^1(\C^{n+m}))$, as appropriate. Similarly, the $\mathbf{B}$ blocks can be replaced by
\begin{align*}
\bZ_1^{(1,3)}&=\pi^\pad D_1\mathbf{f}(\zeta(\hat\psi_s)\hat z_s,\hat\rho_s)\zeta(\hat\psi_s)\pi^{\pad+1,\infty}\\
\bZ_1^{(3,1)}&=\pi^{\pad+1,\infty} D_1\mathbf{f}(\zeta(\hat\psi_s)\hat z_s,\hat\rho_s)\zeta(\hat\psi_s)\pi^\pad \\
\bZ_1^{(3,2)}&=\pi^{\pad+1,\infty} D_1\mathbf{f}(\zeta(\hat\psi_s)\hat z_s,\hat\rho_s)\zeta'(\hat\psi_s)\hat z_s\psi(\cdot) + \pi^{\pad+1,\infty} D_2\mathbf{f}(\zeta(\hat\psi_s)\hat z_s,\hat\rho_s)\\
\bZ_1^{(3,3)}&=\pi^{\pad+1,\infty} D_1\mathbf{f}(\zeta(\hat\psi_s)\hat z_s,\hat\rho_s)\zeta(\hat\psi_s)\pi^{\pad+1,\infty},
\end{align*}
In our implementation, this virtual padding is implemented automatically on an as-needed basis.
\subsubsection{Some notation}
Let $\pi^{j_1,j_2}:\ell_\nu^1(\C^{n+m})\rightarrow \ell_\nu^1(\C^{n+m})$ be defined according to 
$$(\pi^{j_1,j_2}z)_k = \left\{\begin{array}{ll}z_k,&j_1\leq |k| \leq j_2 \\ 0 & \mbox{otherwise.}
\end{array}\right.$$ Let $\mathbf{1}_{j_1,j_2}\subset\ell_\nu^1(\mathbb{C}^{n+m})$ denote the interval sequence
\begin{align*}
(\mathbf{1}_{j_1,j_2})_k&=\left\{\begin{array}{ll} 
[\mathbf{1}_\C]^{n+m}\nu^{-|k|},&j_1\leq|k|\leq j_2 \\ 0 & \mbox{otherwise}, 
\end{array}\right.
\end{align*}
where $[\mathbf{1}_\C]^{n+m}\subset\C^{n+m}$ is the unit interval vector\footnote{Subject to any weighting of norms; that is, the components of $[\mathbf{1}_\C]^{n+m}$ should be scaled so that its norm is equal to 1 in $\C^{n+m}$.} in $\C^{n+m}$. Let $[\mathbf{1}_\R]^{n+m+4}$ denote the unit interval vector\footnote{In case the euclidean norms associated to the scalar variables are not homogeneous and different weighting is used, each component of $[\mathbf{1}_\R]^{n+m+4}$ should be scaled so that the norm of $[\mathbf{1}_\R]^{n+m+4}$ is equal to 1 in the induced max norm.} in $\R^{n+m+4}$. Finally, we will occasionally require the use of
\begin{align*}
T(u,w) & = \left( \begin{array}{c}0_{\R^m} \\ 0_\R \\ 0_\R \\ D_1\theta_{BC}(\sum_k z_k,x,a,\alpha,\beta,\eta)\sum_j w_j \\ 0_\R \\ 0_\R \end{array} \right),
\end{align*}
for $u=(z,(x,a,\alpha,\beta,\eta))\in\ell_\nu^1(\C^{n+m})\times\R^{n+m+4}$. Also, given $w\in\ell_\nu^1(\C^{n+m})$, define $\mathbf{S}(w)\in\ell_\nu^1(\C^{n+m})$ by $\mathbf{S}(w)_{j,0} = \sum_k w_{j,k}$, and $\mathbf{S}(w)_{\cdot,k}=0$ for $k\neq 0$. By construction, $T(u,w)=T(u,\mathbf{S}(w))$.
\subsubsection{The norm $||g^\pad+g^\infty||$}
In this section, the norm $||\cdot||$ on $\mbox{Symm}(\ell_\nu^1(\C^{n+m}))$ is identified with $(x^M,x^\infty)\mapsto||(x^M,0,x^\infty)||_\Omega$.
\begin{lemma}\label{lem-fourierZ1}
Let $z\mapsto\mathbf{f}(z,\cdot)$ be a convolution polynomial of degree $q$. Let\footnote{The absolute value of $\psi_\mathbf{1}$ is precisely the weight attributed to the frequency component, relative to the absolute value norm on $\R$.} $\psi_\mathbf{1} = [-1,1]\cdot ||(0,(0,0,1,0,0))||_{\Omega^M}$. Then
\begin{align*}
||g^\pad + g^\infty||&\leq ||P_s\pi^\pad D_1\mathbf{f}(\zeta(\hat\psi_s)\hat z_s,\hat\rho_s)\mathbf{1}_{\pad+1,q\pad}|| + \frac{1}{\pad+1}|\hat\psi_s|^{-1}||I-Q_sT(\hat u_s,\mathbf{v}_s + \mathbf{w}_s)||,
\end{align*}
where the finitely-supported interval sequences $\mathbf{v}_s,\mathbf{w}_s\subset\ell_\nu^1(\C^{n+m})$ are defined according to
\begin{align*}
\mathbf{w}_s&=\pi^{\pad+1,qM}\Big(D_1\mathbf{f}(\zeta(\hat\psi_s)\hat z_s,\hat\rho_s)\zeta'(\hat\psi_s)\hat z_s\psi_\mathbf{1}^{-1} +  D_2\mathbf{f}(\zeta(\hat\psi_s)\hat z_s,\hat\rho_s)[\mathbf{1}_\R]^{n+m+4}\Big)\\
(\mathbf{v}_s)_{j,k}& = \left\{\begin{array}{ll}
([\mathbf{1}_\C]^{n+m})_j||D_1\mathbf{f}_j(\zeta(\hat\psi_s)\hat z_s,\hat\rho_s)||_{L(\ell_\nu^1(\C^{n+m}),\ell_\nu^1(\C))},&\mbox{$j=1,\dots,m+n$, $k=0$}\\ 0 & \mbox{$j=1,\dots,m+n$, $k\neq 0$}\end{array}\right.
\end{align*}
\end{lemma}
\begin{proof}
Explicitly, the sum $g^\pad+g^\infty$ is
\begin{align*}
g^\pad+g^\infty&=P_s\bZ_1^{1,3}h^\infty + i\hat\psi_s^{-1}(I - Q_s D_3\tilde\F_s^{(2)}(\hat u_s))(K\pi^{\pad+1:\infty})^{-1}(\bZ_1^{3,1}h^\pad + \bZ_1^{3,2}h^\rho + \bZ_1^{3,3}h^\infty) \\
&=P_s\bZ_1^{1,3}h^\infty + i\hat\psi_s^{-1}(I - Q_s T(\hat u_s,\cdot))(K\pi^{\pad+1:\infty})^{-1}(\bZ_1^{3,1}h^\pad + \bZ_1^{3,2}h^\rho + \bZ_1^{3,3}h^\infty)
\end{align*}
Since $z\mapsto \mathbf{f}(z,\cdot)$ is a convolution polynomial of degree $q$, we have $(D_1\mathbf{f}(\zeta(\hat\psi_s)\hat z_s,\hat\rho_s)e_k)_j = 0$ for all $|j|\leq \pad$ whenever $|k|> q\pad$. It follows that the support of $w\mapsto \pi^\pad D_1\mathbf{f}(\zeta(\hat\psi_s)\hat z_s,\hat\rho_s)w$ is contained in the range of $\pi^{0,q\pad}$. Restricting to the range of $\pi^{\pad+1,\infty}$, $$||P_s\bZ_1^{1,3}h^\infty||\leq\left|\left|P_s\pi^\pad D_1\mathbf{f}(\zeta(\hat\psi_s)\hat z_s,\hat\rho_s)\mathbf{1}_{\pad+1,q\pad}\right|\right|.$$
Next, the range of $h^\rho\mapsto \bZ_1^{(3,2)}h^\rho$ is contained in that of $\pi^{0,qM}$ (recall that while we have used a different padding dimension $\pad$ for the computation of norms, our data $\hat z_s$ is still band-limited to $M$ modes). We have
\begin{align*}
||\bZ_1^{3,1}h^\pad + \bZ_1^{3,3}h^\infty||=||D_1\mathbf{f}(\zeta(\hat\psi_s)\hat z_s,\hat\rho_s)\zeta(\hat\psi_s)(h^\pad+h^\infty)|| &\leq ||\mathbf{v}_s||,\\
\bZ_1^{3,2}h^\rho=\pi^{\pad+1,\infty}\big(D_1\mathbf{f}(\zeta(\hat\psi_s)\hat z_s,\hat\rho_s)\zeta'(\hat\psi_s)\hat z_s\psi(\cdot) +  D_2\mathbf{f}(\zeta(\hat\psi_s)\hat z_s,\hat\rho_s)\big)h^\rho &\in \mathbf{w}_s
\end{align*}
Using these bounds/inclusions, the fact that $\mathbf{v}_s$ is supported on the zeroth Fourier index, and the properties of the map $T$, we get
\begin{align*}
||(I - Q_s T(\hat u_s,\cdot))(K\pi^{\pad+1:\infty})^{-1}(\bZ_1^{3,1}h^\pad + \bZ_1^{3,2}h^\rho + \bZ_1^{3,3}h^\infty)||&\leq\frac{||I-Q_sT(\hat u_s,\mathbf{v}_s + \mathbf{w}_s)||}{\pad +1}.
\end{align*}
Combining the previous bounds, we get the result claimed in the lemma.
\end{proof}

\begin{remark}
While perhaps symbolically intimidating, all of the quantities appearing in the statement of Lemma \ref{lem-fourierZ1} are explicitly machine-computable, so a rigorous uniform (in $s$ and $h$) bound for $||g^\pad + g^\infty||$ can indeed be computed.
\end{remark}


\subsubsection{The norm $||g^\rho||$}
In this section, the norm $||\cdot||$ on $\R^{n+m+4}$ is identified with $x\mapsto||0,x,0)||_\Omega$. The bound for $||g^\rho||$ bears a lot of similarity to the one for $||g^\pad + g^\infty||$, and its derivation is similar. We omit the proof of the following lemma.
\begin{lemma}\label{lem-grho-Z1}
With the same notation as in Lemma \ref{lem-fourierZ1}, we have the bound
\begin{align*}
||g^\rho||&\leq ||R_s\pi^\pad D_1\mathbf{f}(\zeta(\hat\psi_s)\hat z_s,\hat\rho_s)\mathbf{1}_{\pad+1,q\pad}|| + \frac{1}{\pad+1}|\hat\psi_s|^{-1}||S_sT(\hat u_s,\mathbf{v}_s + \mathbf{w}_s)||.
\end{align*}
\end{lemma}

\subsubsection{Summary of the $Z_1$ bound}
Combining the results of Lemma \ref{lem-fourierZ1} and Lemma \ref{lem-grho-Z1}, it follows that if we choose $Z_1$ such that
\begin{align*}
Z_1&\geq \max\left\{ ||P_s\pi^\pad D_1\mathbf{f}(\zeta(\hat\psi_s)\hat z_s,\hat\rho_s)\mathbf{1}_{\pad+1,q\pad}|| + \frac{1}{\pad+1}|\hat\psi_s|^{-1}||I-Q_sT(\hat u_s,\mathbf{v}_s + \mathbf{w}_s)||\hspace{1mm}, \right.\\
&\phantom{\geq\max\big\}\}}\left. ||R_s\pi^\pad D_1\mathbf{f}(\zeta(\hat\psi_s)\hat z_s,\hat\rho_s)\mathbf{1}_{\pad+1,q\pad}|| + \frac{1}{\pad+1}|\hat\psi_s|^{-1}||S_sT(\hat u_s,\mathbf{v}_s + \mathbf{w}_s)|| \right\}
\end{align*}
for all $s\in\Delta$, then $Z_1$ will satisfy the bound \eqref{thm-rpa-Z1} for our validated continuation problem. In the above bound, one must remember that the norms $||\cdot||$ are actually restrictions of the norm $||\cdot||_\Omega$ to various subspaces; see the preamble before Lemma \ref{lem-fourierZ1} and Lemma \ref{lem-grho-Z1}.

\subsection{The bound $Z_2$}\label{sec-Z2}
It is beneficial to perform a further splitting of the expression that defines the $Z_2$ bound. Recall that $Z_2(r)$ should satisfy
\begin{align}\label{Z2r-section}
||A_s(D\mathcal{F}_s(\hat u_s+\delta) - D\mathcal{F}_s(\hat u_s))||_{B(\Omega,\Omega)}\leq Z_2(r)
\end{align}
uniformly for $s\in\Delta$ and $\delta\in\overline{B_r(\hat u_s)}$. From Remark \ref{remark-D2F-DNE}, $D\mathcal{F}_s(\cdot)$ is not itself differentiable, so appealing to a ``typical" second derivative estimate for this $Z_2$ bound is not going to work. To handle this, we use the strategy of \cite{VandenBerg2020a} and split the norm to be computed using a triangle inequality. Given $\delta\in \overline{B_r(\hat u_s)}$, split it as $\delta = \delta_1+\delta_2$, where the components of $\delta_1$ are all zero \emph{except for the frequency component, $\psi$}, and $\delta_2$ has zero for its frequency component. This decomposition is uniquely defined. Then, consider two new bounds to be computed:
\begin{align*}
||A_s(D\mathcal{F}_s(\hat u_s + \delta_1) - D\mathcal{F}_s(\hat u_s))||_{B(\Omega,\Omega)}&\leq Z_{2,1}(r)\\
||A_s(D\mathcal{F}_s(\hat u_s + \delta_1+\delta_2) - D\mathcal{F}_s(\hat u_s + \delta_1))||_{B(\Omega,\Omega)}&\leq Z_{2,1}(r)
\end{align*}
uniformly for $s\in\Delta$ and $\delta\in\overline{B_r(\hat u_s)}$. If we define $Z_2(r)\bydef Z_{2,1}(r)+Z_{2,2}(r)$, then \eqref{Z2r-section} will be satisfied. It turns out that this decomposition is effective. We will elaborate on this now.

\subsubsection{The bound $Z_{2,1}$}\label{sec-Z21}
The partial derivative $\frac{d}{d\psi(\rho)}D\mathcal{F}_s(\hat z^M_s,\hat\rho_s,0)$ exists and is continuous. This is a consequence of Lemma \ref{lem-2nd-derivatives} (n.b.\ the inputs are band-limited) and the structure of $\F_s$; see \eqref{Fmap}. As $A_s$ is bounded, we can use the fundamental theorem of calculus (in Banach space) to obtain
\begin{align}\label{Z21}
||A_s(D\F_s(\hat u_s+\delta_1)-D\F_s(\hat u_s))h||_\Omega&\leq \sup_{t\in[0,1]}\left|\left|A_s\left(\frac{d}{d\psi(\rho)}D\F_s(\hat u_s + t\delta_1)h\right)\psi(\delta_1)\right|\right|_\Omega
\end{align}
for any $h\in\Omega$, $||h||_\Omega\leq 1$. Due to the structure of $\F_s$, the derivative term in the large parentheses will have 
\begin{itemize}
\item in the $\ell_\nu^1(\C^{n+m})$ components: convolution polynomials of the form \eqref{poly} involving Fourier components coming from the set $$\{\zeta(\hat\psi_s + \psi(\delta_1))\hat z_s\}\cup\{i\zeta(\hat\psi_s + \psi(\delta_1))K\hat z_s\}\cup\{-\zeta(\hat\psi_s + \psi(\delta_1))K^2\hat z_s\}\cup\{h\}$$ with coefficients in $\hat\rho_s$;
\item in the scalar components: finite-dimensional range functions of $\sum_{|k|\leq M}(\hat z_s)_k$ and $\delta(\psi_1)$.
\end{itemize}
To obtain the tightest possible upper bound for \eqref{Z21}, one would want to exploit as much structure of $D\F_s$ (and the directional derivative in the frequency direction) as possible. However, to get a general-purpose \emph{implementable} bound, we can apply the following operations to $\frac{d}{d\psi(\rho)}D\F_s(\hat u_s + t\delta_1)h\psi(\delta_1)$.
\begin{enumerate}
\item Replace all instances of $\zeta(\hat\psi_s + \psi(\delta_1))$ with the interval $[-1,1]$;
\item replace the component of $h$ in $\R^{n+m+4}$ with the unit interval vector\footnote{Note that the components of this ``unit" interval vector might be uneven due to weighting.} in $\times\R^{n+m+4}$;
\item replace all other variables\footnote{\label{footnote-ellnu}Note that if the norm on $\C^{n+m}$ associated to $\ell_\nu^1(\C^{n+m})$ is itself weighted, this must be taken into account.} described in the bulleted list above with bounds for their norms, multiplied by the zero-supported basis element of the relevant space (e.g.\ for $\ell_\nu^1$ elements, the identity $e_0$ of the Banach algebra; for $\ell_\nu^1(\C^{n+m})$, the vectorized version of $e_0$), multiplied by the interval $[-1,1]$;
\item compute operator norms of the blocks of $A_s$ in \eqref{As};
\item complete block multiplications, taking the norm of the result, followed by an interval supremum (note, $t\in[0,1]$ is also replaced by an interval).
\end{enumerate}
We emphasize in the algorithm above, every object is a finite-dimensional quantity, so operations can indeed be performed in a suitable complex vector space on the computer. This strategy produces a true bound for \eqref{Z21} largely because of the Banach algebra and the interval arithmetic. For example, consider the impact of this computation on the quantity
$$-\hat z_s*h + (e^{iK(\hat\psi_s + \psi(\delta_1))}\hat z_s)*(iKe^{iK(\hat\psi_s+\psi(\delta_1)}\hat z_s)*h$$
in $\ell_\nu^1$. From the Banach algebra, if $||h||_\nu\leq 1$, then a triangle inequality produces
\begin{align*}
||\hat z_s|| + ||e^{iK(\hat\psi_s + \psi(\delta_1)}\hat z_s||_\nu\cdot||iKe^{iK(\hat\psi_s + \psi(\delta_1)}\hat z_s||_\nu &= ||(ae_0)*(ce_0)||_\nu +  ||(ae_0)*(be_0)*(ce_0)||\\
&= ac + abc,
\end{align*}
where $a = ||\hat z_s||_\nu$, $b = ||iK\hat z_s||_\nu$, $c=||h||_\nu=1$, and $e_0$ is the identity element in the Banach algebra $\ell_\nu^1$. Now, define $\mathbf{1} = [-1,1]$ and compare the result to
\begin{align*}
-(\mathbf{1}ae_0)*(\mathbf{1}ce_0) + (\mathbf{1}a e_0)*(\mathbf{1}be_0)*(\mathbf{1}ce_0).
\end{align*}
The support of this sequence is the zeroth index, and we can explicitly compute the resulting interval. It is precisely
$\mathbf{1}\cdot(ac + abc)e_0$. The supremum of this interval is $ac + abc$, which matches the analytical triangle inequality / Banach algebra bound computed previously.

\begin{remark}
The bound obtained by applying the strategy above is incredibly crude; in effect, we use the triangle inequality for everything. However, producing fully general code to compute the second derivatives for this class of problems (polynomial delay differential equations with arbitrary numbers of delays) would be a messy programming task. Even with the present implementation, where we need only compute second derivatives evaluated at sequences that are supported on the zeroth Fourier mode, the implementation was far from trivial. Long term, it would be beneficial to implement second derivatives, as this would also allow for improvements to the $Y_0$ bound; see Section \ref{sec-Y0}. The good news is that, theoretically, the computed bound will be $O(r)$ for $r$ small enough, due to the linear multiplication of $\psi(\delta_1)=O(r)$ appearing in \eqref{Z21}.
\end{remark}

\subsubsection{The bound $Z_{2,2}$}
Let $d_{\delta_2}$ denote the Gateaux derivative of $D\F_s(\cdot)$ in the direction $\delta_2$. We claim $$d_{\delta_2}D\F_s(\hat u_s+\delta_1+t\delta_2)$$ exists and is continuous for $t\in[0,1]$. To see why, observe that that $D\F_s^{(2)}$ is continuously differentiable, so we need only worry about the part in $\ell_\nu^1(\C^{n+m})$, namely the components $D\F_s^{(1)}$ and $D\F_s^{(3)}$ that come from the vector field. The result now follows from Lemma \ref{lem-2nd-derivatives}. Indeed, a band-limited input is only required for the double derivative with respect to the frequency variable, and $\delta_2$ has zero for its frequency component. By the fundamental theorem of calculus, 
\begin{align}\label{Z22}
||A_s(D\F_s(\hat u_s+\delta_1+\delta_2)-D\F_s(\hat u_s+\delta_1))h||_\Omega&\leq \sup_{t\in[0,1]}\left|\left|A_s\left(d_{\delta_2}D\F_s(\hat u_s + \delta_1 +  t\delta_2)h\right)\right|\right|_\Omega
\end{align}
where we take $h\in\Omega$, $||h||_\Omega\leq 1$. 

To bound \eqref{Z22}, we make use of a very similar strategy to the one from Section \ref{sec-Z21} for $Z_{2,1}$. The difference here is that the convolution polynomials in the $\ell_\nu^1(\C^{n+m})$ part of the Gateaux derivative have a more difficult structure. The main problem is in the mixed derivatives with respect to frequency and Fourier space; see the derivatives $\frac{d}{dz_q}\frac{d}{d\psi(\rho)}$ in Lemma \ref{lem-2nd-derivatives}. These derivatives involve the action of the derivative operator $K$ on the sequence part of $h$, which is not necessarily band-limited. If $h=(h_z,h_\rho)$ for $h_z\in\ell_\nu^1$, and this sequence is not band-limited, then generally $Kh_z$ will not have a finite $\ell_\nu^1$-norm. To combat this problem, we remark that \emph{only one such term can appear in any given convolution polynomial}. This can be exploited as follows.

To begin, we factor $A_s$ as follows.
$$A_s = \left(\begin{array}{ccc}(M+1) P_s & Q_s & V_s^{(1)} \\ (M+1) R_s & S_s & V_s^{(2)} \\ 0&0&V_s^{(3)}\end{array}\right)\left(\begin{array}{ccc}I\frac{1}{M+1}&0&0\\0&I&0\\0&0&(K\pi^\infty)^{-1}\end{array}\right)\bydef \mathbf{A}_s \mathbf{K}^{-1}.$$ 
Note that $\mathbf{A}_s$ is obtained from $A_s$ by multiplying (on the right) the first ``column" by $(M+1)$, and the last ``column" by $K\pi^\infty$. The following lemma is a specific case of a lemma of van den Berg, Groothedde and Lessard \cite{VandenBerg2020a}.

\begin{lemma}
Define an operator $\tilde K^{-1} = \frac{1}{M+1}\pi^M + K^{-1}\pi^\infty$ on $\ell_\nu^1$. Let $u,v,\in\ell_\nu^1$. Then
\begin{align*}
||\tilde K^{-1}(Ku*v)||_\nu&\leq C_\nu||u||_\nu||v||_\nu,\hspace{1cm}
C_\nu = \left\{\begin{array}{ll}\frac{\nu^{2M+2}}{e\log\nu^{2M+2}},&\nu^{2M+2}<e \\ 1,&\nu^{2M+2}\geq e.
\end{array}\right.
\end{align*}
\end{lemma}

Now consider the product $\mathbf{K}^{-1}d_{\delta_2}D\F_s(\hat u_s + \delta_1 +  t\delta_2)h$. The part in $\ell_\nu^1(\C^{n+m})$ of the Gateaux derivative $d_{\delta_2}D\F_s(\hat u_s + \delta_1 +  t\delta_2)h$ will be multiplied by the diagonal operator $(\tilde K^{-1},\dots,\tilde K^{-1})$. So consider how we might bound a term of the form
\begin{align}\label{Kinv-Z2-thing1}\tilde K^{-1}\left(\alpha iKh_z*\prod_{j}e^{iK\psi(\hat\rho_s + \psi(\delta_1))\mu_j}(\hat z_s + tz(\delta_2))_j \right) \in\ell_\nu^1(\C)\end{align}
 in a given convolution polynomial that contains a factor $Kh_z$, for $h_z\in\ell_\nu^1$ and $||h_z||_\nu\leq 1$. Here, $\alpha$ is a scalar that could depend on $\hat\rho_s$, the scalar components of $t\delta_2$, and the delays, while $z(\delta_2)$ is the part of $\delta_2$ in $\ell_\nu^1(\C^{n+m})$. Note that by permutation invariance of the convolution, we may assume that $Kh_z$ appears as the first term on the left, as we have done here. By Lemma 9, the above is bounded above by
\begin{align*}
C_\nu\left|\left|\alpha\prod_{j}e^{iK\psi(\hat\rho_s + \psi(\delta_1))\mu_j}(\hat z_s + tz(\delta_2))_j\right|\right|_\nu.
\end{align*}
We can achieve the same effect by replacing $Kh_z$ in \eqref{Kinv-Z2-thing1} with $C_\nu e_0$, for $e_0$ being the identity in the convolution algebra, and taking the $\ell_\nu^1$-norm. When $\tilde K^{-1}$ is multiplied by a polynomial term that does \emph{not} contain a factor $Kh_z$, then a straightforward calculation shows that $||\tilde K^{-1}||_{B(\ell_\nu^1,\ell_\nu^1)}=(M+1)^{-1}$. The net result is that the impact on the norm of the polynomial term is a scaling by $(M+1)^{-1}$. Therefore, adapting the algorithm from the $Z_{2,1}$ bound calculation, it suffices to do perform the following operations to $d_{\delta_2}D\F_s(\hat u_s + \delta_1 +  t\delta_2)h$.
\begin{enumerate}
\item Multiply the part of $d_{\delta_2}D\F_s(\hat u_s + \delta_1 +  t\delta_2)h$ in $\ell_\nu^1(\C^{n+m})$ by $\frac{1}{M+1}$; 
\item replace all instances of $\zeta(\hat\psi_s + \psi(\delta_1))$ with the interval $[-1,1]$;
\item replace all (remaining) instances of $K$ with $(M+1)C_\nu e_0\cdot[-1,1]$; 
\item replace the component of $h$ in $\R^{n+m+4}$ with the unit interval vector\footnote{Note that the components of this ``unit" interval vector might be uneven due to weighting.} in $\R^{n+m+4}$; 
\item replace all other variables ($\hat z_s$, $\hat \rho_s$, $\delta_1$, $\delta_2$, the part of $h$ in $\ell_\nu^1(\C^{n+m})$, and their relevant projections\footnote{See footnote \ref{footnote-ellnu} concerning the projections of $h$ in $\ell_\nu^1(\C^{n+m})$.}) with bounds for their norms, multiplied by the zero-supported basis element of the relevant space, multiplied by the interval $[-1,1]$;
\item compute operator norms of the blocks of $\mathbf{A}_s$;
\item complete block multiplication of $\mathbf{A}_s$ on the left, take the $\Omega$-norm of the result, followed by an interval supremum (note: $t\in[0,1]$ is also replaced by an interval).
\end{enumerate}
The result is a (crude) upper bound of $\sup_{t\in[0,1]}||A_sd_{\delta_2}D\F_s(\hat u_s+\delta_1+t\delta_2)h||_\Omega$ for $||h||_\Omega\leq 1$. This bound is expected to be $O(r)$ for $r$ small, since the Gateaux derivative $d_{\delta_2}D\F_s(\hat u_s+\delta_1+t\delta_2)h$ is $O(||\delta_2||)$ for $||\delta_2||\leq r$ small. 
\begin{remark}
In step 3, the variable replacement negates (by multilinearity of the convolution) the multiplication of the relevant polynomial term by $(M+1)^{-1}$ that was done in step 1, while propagating forward the correct bound $C_\nu$ that results from the interaction between $\tilde K^{-1}$ and the $Kz_h*\prod_j(\cdots)$ polynomial term. This somewhat roundabout way of introducing the bound $C_\nu$ in the correct locations is done in order to make the process implementable in generality.
\end{remark}


\section{Specification to ordinary differential equations}\label{sec:ode}
In Section \ref{sec3}, we demonstrated how rigorous two-parameter continuation of periodic orbits in delay differential equations can be accomplished in such a way that the continuation can pass through degenerate Hopf bifurcations. As ordinary differential equations are a special case -- that is, where any delays are identically zero -- the theory of the previous sections naturally applies equally to them. However, the formulation of the map \eqref{eq:F-DDE} can be greatly simplified. Indeed, the reader familiar with numerical methods for periodic orbits has likely noticed that we have not performed the ``usual" scaling out by the frequency, so that the period can be abstractly considered as $2\pi$. With delay differential equations, this is not beneficial because it merely moves the frequency dependence into the delayed variables. Additionally, it is the presence of the delays that requires a more subtle analysis of the $Z_2$ bound; see Section \ref{sec-Z2}. 

To compare, with ordinary differential equations \emph{without delays}, the technicalities with the $Z_2$ bound are absent. At the level of implementation, computing second and even third derivatives of the vector field in the Fourier representation is also much easier. In this section, we will present an alternative version of the zero-finding problem of Section \ref{sec-zfp} and analogous map $G$ from \eqref{eq:F-DDE}. However, we will not discuss the general implementation of the technical bounds $Y$ and $Z$, since they are both simpler than the ones we have previously presented for the delay differential equations case, and can be obtained by fairly minor modifications of the bounds from \cite{VandenBerg2021a}. We have implemented them in general in the codes at \cite{CodeURL}.

\subsection{Desingularization, polynomial embedding and phase isolation}
The set-up now starts with an ordinary differential equation
\begin{align}\label{eq:ODE}\dot y(t)=f(y(t),\alpha,\beta)\end{align} 
again depending on two real parameters $\alpha$ and $\beta$. Performing the same blowup procecure as we did for the delay equations, we define 
\begin{align*}
\tilde f(z,x,a,\alpha,\beta)=\left\{\begin{array}{ll}
a^{-1}(f(x+a z,\alpha,\beta)-f(x,\alpha,\beta)),&a\neq 0\\
d_{y}f(x,\alpha,\beta)z,&a=0.\end{array}\right.
\end{align*}
The goal is therefore to find a pair $(x,z)$ such that
\begin{align*}
f(x,\dots,x,\alpha,\beta)&=0\\
\dot z(t) &= \tilde f(z(t),x,a,\alpha,\beta),\\
||z||&=1
\end{align*}
where $z$ is $\omega$-periodic for an unknown period $\omega$; equivalently, the frequency of $z$ is $\psi = \frac{2\pi}{\omega}$. Let $\mu=\psi^{-1}$ denote the reciprocal frequency. Define $\tilde z(t) = z(t\mu)$. Then $\tilde z(t)$ is $2\pi$-periodic. Substituting this into the differential equation above and dropping the tildes, we obtain the modified vector field
\begin{align*}
\dot z(t) &= \mu\tilde f(z(t),x,a,\alpha,\beta),
\end{align*}
where now the scope is that $z$ should be $2\pi$-periodic.

If a polynomial embedding is necessary to eliminate non-polynomial nonlinearities, we allow the inclusion of $m$ additional scalar equations that must be simultaneously solved, where we introduce an appropriate number $m$ of unfolding parameters, $\eta\in\R^m$, to balance them. We again use the symbol $\theta_{BC}$ to function that defines these boundary conditions, with the equation being $\theta_{BC}(z(0),x,a,\alpha,\beta,\eta)=0$. We use the same anchor condition to handle the lack of isolation from phase shifts.
\subsection{Zero-finding problem}
Abusing notation and assuming now that $\tilde f$ is polynomial after the embedding has been taken into account, the zero-finding problem is 
\begin{equation}\label{eq:zfp-ODE}
\begin{cases}
\dot z = \mu \tilde f(z(t),x,a,\alpha,\beta,\eta),&\quad\text{(differential equations)}\\
\| z\| = 1, &\quad\text{(amplitude condition of scaled orbit)}\\
\int\langle z(s),\hat z'(s)\rangle ds=0, &\quad\text{(anchor condition)}\\
f(x,\alpha,\beta)=0, &\quad\text{($x$ is a steady state)}\\
\theta_{BC}(z(0),x,a,\alpha,\beta,\eta)=0.&\quad\text{(embedding boundary condition)}
\end{cases}
\end{equation}
where $\hat z$ is understood to be a candidate numerical solution. Expanding $z$ as a Fourier series with period $2\pi$, we have $z(t) = \sum_{k\in\Z}z_ke^{ikt}$ for some complex vectors $z_k$ obeying the symmetry $z_k = \overline{z_{-k}}$. If we now allow $\mathbf{f}(z,x,a,\alpha,\beta,\eta)$ to be the representation of $\tilde f$ as a convolution polynomial in the coefficients $z=(z_k)_{k\in\Z}$, then an analogous derivativation to the delay differential equations case produces the map
\begin{align*}
G(z,x,a,\psi,\eta,(\alpha,\beta))&=\left(\begin{array}{c}
-iKz + \mu\mathbf{f}(z,x,a,\alpha,\beta,\eta) \\ f(x,\alpha,\beta) \\
\langle z,K^2\hat z\rangle - 1 \\
\langle z,iK\hat z\rangle \\ 
\theta_{BC}(\sum_k z_k,x,a,\alpha,\beta,\eta)
\end{array}\right)
\end{align*}
defined on the same Banach spaces, with the same codomain. However, this time, one can show that if $f$ is $C^k$, then $G$ really is $C^{k-1}$.

\section{Proving bifurcations and bubbles}\label{sec:analytical_proofs}
Modulo non-resonance conditions, we would generically expect Hopf bifurcations to occur on the level curve at amplitude zero of the computed 2-manifolds of Section \ref{sec3}. Hopf bubbles are equivalent to curves in the manifold that intersect the amplitude zero surface at two distinct points. 
As for bubble bifurcations, we can describe these in terms of a the existence of a relative local minimum of a projection of the computed 2-manifold, represented as the graph over amplitude and a distinguished parameter. We develop these points in this section, demonstrating how post-processing of data from the computer-assisted proofs can be used to prove the existence of Hopf bifurcations, bubbles, and degenerate bifurcations.

We should emphasize that to prove the existence of a \textit{single} Hopf bubble, it suffices to identify a curve connecting two points at amplitude zero in the projection of the proven manifold in $\alpha\times\beta\times\mbox{amplitude}$ space, for one of $\alpha$ or $\beta$ being fixed. While this does require verifying two Hopf bifurcations, the rest of the proof of an isolated bubble is fairly trivial, requiring only determining a sequence of simplices that enclose the curve. Therefore, we will only discuss how to prove bubble \textit{bifurcations} in this section, since this can require a bit more analysis. We will not analyze Bautin bifurcations.

To simplify the presentation, we will assume throughout this section that the norm on $\Omega$ is such that $||(0,(0,a,0,0,0),0)||=|a|$. That is, the amplitude component is unit weighted relative to the absolute value.

\subsection{Hopf bifurcation curves}
In delay differential equations, the sufficient conditions for the existence of a Hopf bifurcation include the non-resonance check, which involves counting all eigenvalues all eigenvalues of the lienarization on the imaginary axis. This is somewhat beyond the scope of our work here, although there is literature on how this can be done using computer-assisted proofs \cite{Church2022,P.Lessard2020}. In this paper, we will consider the following related notion.
\begin{definition}\label{def:hopfcurve}
The delay differential equation \eqref{sec3} has a \emph{Hopf bifurcation curve $\Theta_H\in\Delta$} if there exists a $C^1$ parametrization $\Delta\ni s\mapsto (\alpha(s),\beta(s),x(s),y(s))$ such that $x(s)$ is a steady state solution and $y(s)$ is a periodic orbit at parameters $(\alpha(s),\beta(s))$, 
such that:
\begin{itemize}
\item $\Theta_H(0)$ and $\Theta_H(1)$ are on the boundary of $\Delta$, and $\Theta_H(t)$ is in the interior of $\Delta$ for $t\in(0,1)$.
\item $x\circ\Theta_H=y\circ\Theta_H$; in other words, $y$ is a steady state on restriction to the image of $\Theta_H$.
\item For $s\in\Delta\setminus\{\Theta_H(t):t\in[0,1]\}$, $y(s)$ is not a steady state.
\end{itemize}
\end{definition}
Contrary to the usual Hopf bifurcation, we do not reference the \emph{direction} of the bifurcation, even along one-dimensional curves in the manifold of periodic orbits. 

\begin{figure}\centering
\includegraphics[scale=0.5]{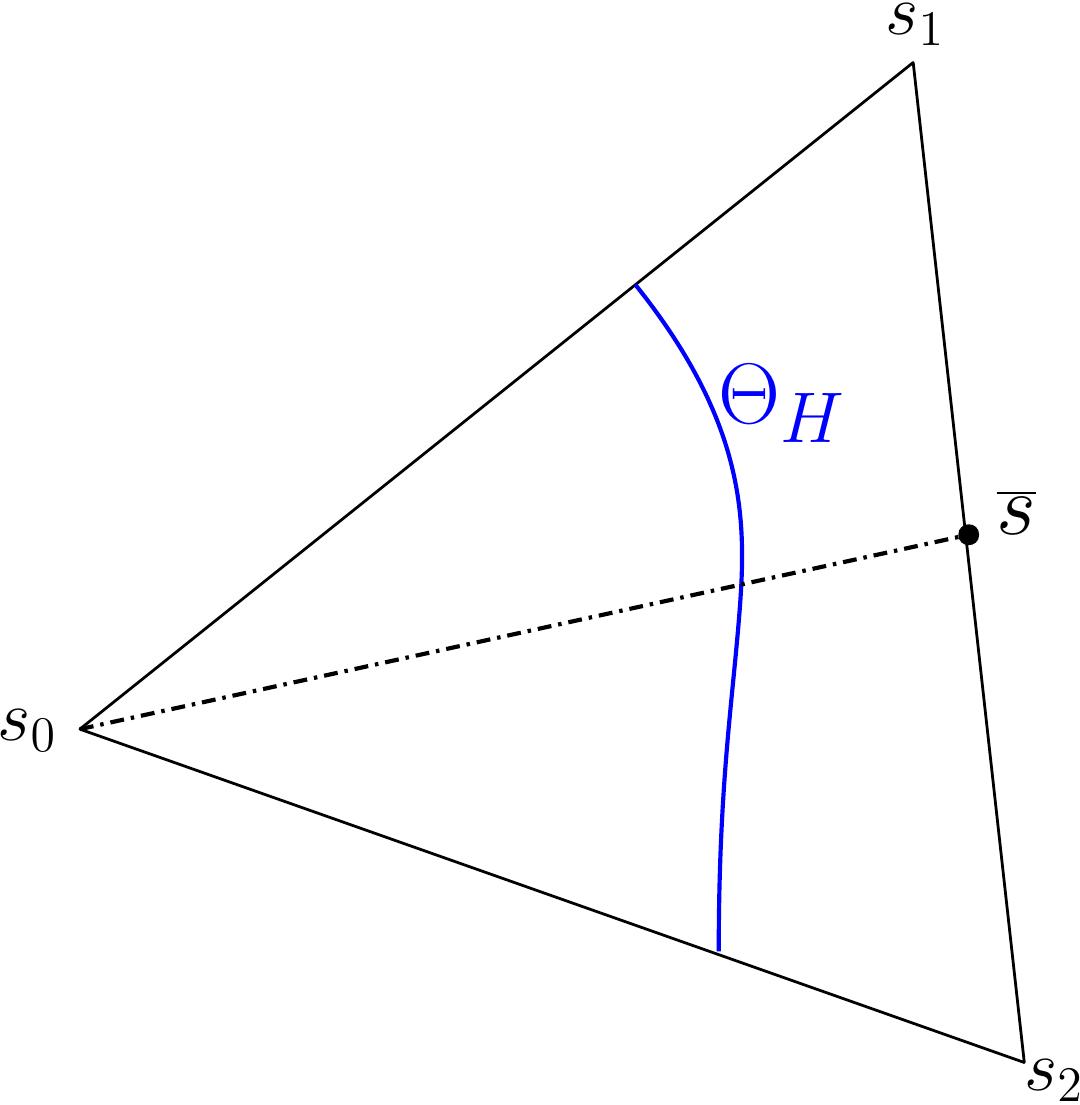}
\caption{Briefly, the partial derivatives of the amplitude of the family of periodic orbits, with respect to the simplex parametrization variables, have the same sign. From the assumptions of the proposition, the sign of the amplitude along the edge connecting $s_1$ to $s_2$ is opposite to that of the amplitude at $s_0$. Forming a ray (dashed-dotted line) connecting $s_0$ to some point $\overline s$ on the edge connecting $s_1$ to $s_2$, the mean-value theorem guarantees the existence of a unique point along the ray at which the amplitude is zero. Parameterizing over all $\overline s$ on the edge $s_1$ to $s_2$, the Hopf curve $\Theta_H$ (blue) is constructed.}
\label{fig-assist-proof-prop10}
\end{figure}

The existence of a Hopf bifurcation curve can be proven using post-processing of a validated contiuation on given simplex on which one \emph{expects} a Hopf bifurcation to occur. The idea is that on a simplex that encloses a Hopf bifurcation, the amplitude in the blown-up variables (see Section \ref{sec-desingularization}) is expected to cross through zero. Since we are working in two-parameter continuation, such a crossing point should persist as a curve. A geometric construction based on the partial derivatives of the amplitude with respect to the simplex paramaterization can be used to construct this curve. See Figure \ref{fig-assist-proof-prop10} for a visualization.

To establish a more constructive proof, we need to introduce a few projection maps. Let $u=(z^M,\rho)\in\Omega^M$, $\rho = (x,a,\psi,\eta,\alpha,\beta)$. Denote $\pi^a u = a$ the projection to the amplitude component, and $\pi^{(\alpha,\beta)}u = (\alpha,\beta)$ the projection onto the parameters. Similarly, write $\pi^a\rho = a$. Denote the three vertices of the standard simplex $\Delta$ as $s_0 = (0,0)$, $s_1 = (1,0)$ and $s_2 = (0,1)$.

\begin{proposition}\label{prop-Hopf}
Suppose a simplex containing the interpolated numerical data $\hat u_s$ has been validated, for $s\in\Delta.$ Let $r>0$ be an associated validation radius from the radii polynomial. Let $c = Z_0+Z_1 + Z_2(r)$. Suppose the following are satisfied.
\begin{itemize}
\item For $k=1,2$, $$[\pi^a(\hat u_{s_0})-r,\pi^a(\hat u_{s_0}) +r]\cdot[\pi^a(\hat u_{s_k})-r,\pi^a(\hat u_{s_k}) +r]<0,$$ where the multiplication is in the sense of interval arithmetic.
\item The sign of the interval $[\pi^a(\hat u_s) - r, \pi^a (\hat u_s)+r]$ is constant for $s$ on the edge of $\Delta$ not incident with $s_0$; that is, on the edge connecting $s_1$ with $s_2$. 
\item Let $\partial_s$ denote the Fr\'echet derivative operator with respect to the variable $s\in\Delta$. For $h\in\R^2$, define 
\begin{align}
\label{derivatives-tilde-au}\tilde\partial_s a(s)h&\bydef \pi^a\left(-S_s\partial_s[\mathbf{J}_s](\hat u_s + \delta_r)h\right),\hspace{1cm}\tilde\partial_s u(s)\bydef -A_s\partial_s [\F_s](\hat u_s + \delta_r), \\  \gamma &= \frac{c}{1-c}\sup_{s\in\Delta}||\tilde\partial_s u(s)|| \nonumber,\end{align}
 where $\delta_r=\overline{B_r(0)}\subset\Omega^M$.  That is, $\gamma$ is a real interval and $\tilde\partial_s a(s)$ is an inteval in $(\R^2)^*$. Let $\hat\gamma = \sup\gamma$. 
The components of $\tilde\partial_s a(s) + \hat\gamma[-1,1]^2$ have the same sign for $s\in\Delta$. 
\end{itemize}
Then there is a \textit{unique} (up to parameterizaztion) Hopf bifurcation curve $\Theta_H:[0,1]\rightarrow\Delta$, and each of $\Theta_H(0)$ and $\Theta_H(1)$ lie on the edges of $\Delta$ that are incident with $s_0$. 
\end{proposition}

\begin{proof}
Suppose the radii polynomial proves the existence of $u_s\in B_r(\hat u_s)$ such that $\F_s(u_s)=0$. The $C^1$ parametrization of the steady state ($x(s)$), periodic orbit ($y(s)$) and amplitude parametrization is a consequence of the computer-assisted proof. $\partial_s u_s$ exists, and $\partial_s \F_s(u_s) + D\F_s(u_s)\partial_s u_s = 0$. We claim $D\F_s(y)$ is boundedly invertible for $y\in \overline{B_r(\hat u_s)}$. First, 
\begin{align*}
I-A_sD\F_s(y) &= (I-A_sA_s^\dagger) + A_s(A_s^\dagger - D\F_s(\hat u_s)) + A_s(D\F_s(\hat u_s) - D\F_s(y))
\end{align*}
whose norm is bounded above by $Z_0 + Z_1 + Z_2(r)=c$. Since the radii polynomial has validated $\hat u_s$, it follows that $c<1$. By the Neumann series, $A_sD\F_s(y)$ is boundedly invertible, with norm bounded by $(1-c)^{-1}$. Since $A_s$ is boundedly invertible by its construction, the same is true of $D\F_s(y)$.

It follows that $\partial_s u_s = -D\F_s(u_s)^{-1}\partial_s\F_s(u_s)$. Observe that
\begin{align*}
(A_s-D\F_s(u_s)^{-1})\partial_s\F_s(\hat u_s + \delta_r)&=(AD\F_s(u_s))^{-1}(AD\F_s(u_s)-I)A_s\partial_s\F_s(\hat u_s + \delta_r)
\end{align*}
By previous estimates, the above is bounded in norm by $\gamma$. Since $\partial_s[\F_s](u_s)\in\partial_s[\F_s](\hat u_s + \delta_r)$ --- see Section \ref{sec-Y0} --- we have 
\begin{align}
\label{derivative-us}\partial_s u_s \in -A_s\partial_s\F(\hat u_s + \delta_r) + \overline{B_\gamma(0)}.
\end{align}
Note that $\tilde\partial_s a(t)h$ is precisely the amplitude component of $\partial_s (u_s h)|_{s=t}$.

Consider the line $[0,1]\mapsto t\mapsto s_1 + t(\overline s - s_1)$, where $\overline s$ is any point on the edge connecting $s_1$ to $s_2$ in $\Delta$. Consider the function $g(t) = \pi^a u_{s_1 + t(\overline s - s_1)}$. By the assumptions of the proposition and the inclusion \eqref{derivative-us}, we have $g(0)g(1)<0$ and $\frac{d}{dt}g(t)\neq 0$ for $t\in(0,1)$. By the mean-value theorem, there is a unique $t^*\in(0,1)$ such that $g(t^*)=0$. Moreover, $t^*=t^*(\overline s)$ depends continuously on $\overline s$. Letting $p:[0,1]\rightarrow\Delta$ be a $C^1$ parametrization of the edge connecting $s_1$ to $s_2$, it follows by the definition of the map $\F_s$ that $\Theta_H=t^*\circ p$ is a Hopf bifurcation curve. 
\end{proof}

\begin{remark}\label{rem-Hopf}
Proposition \ref{prop-Hopf} is stated in such a way that we assume $s_0$ and the edge connecting $s_1$ to $s_2$ lie on the opposite sides of Hopf curve $\Theta_H$. This is done \emph{almost} without loss of generality. First, if the Hopf curve should bisect the simplex in such a way that $s_1$ or $s_2$ is separated from its opposing edge, a suitable re-labeling of the nodes transforms the problem into the form of the proposition. Second, a problem can arise if the Hopf curve intersects one of the vertices $s_j$, $j=0,1,2$. Similarly, numerical difficulties can arise with the method of computer-assisted proof if the intersections of $\Theta_H$ with the edges of $\Delta$ occur very close to the vertices. These problems can be alleviated by a careful selection of numerical data.
\end{remark}

\begin{remark}
We emphasize that most of the sufficient conditions of Proposition \ref{prop-Hopf} can be checked using the output of the validation of a simplex. We also mention specifically that in the second point, the sign of $[\pi^a(\hat u_s) - r, \pi^a (\hat u_s)+r]$ need only be verified to match at $s=s_1$ and $s=s_2$, due to linearity of $\hat u_s$. The exception is the computation of $\tilde\partial_s a$ and $\gamma$. The former is a finite computation, and the latter can be bounded in a straightforward manner; see Section \ref{sec-Y0} for details concerning the structure of $\partial_s\F_s$.
\end{remark}

\begin{remark}
A Hopf bifurcation curve necessarily generates continua of Hopf bifurcations. These can be constructed by selecting curves in $\Delta$ that are transversal to $\Theta_H$. 
\end{remark}

\subsection{Degenerate Hopf bifurcations}
In this section, we consider how one might prove the existence of a degenerate Hopf bifurcation, be it a bubble bifurcation or otherwise. For a first definition, we consider the bubble bifurcation.
\begin{definition}\label{def-degenerate}
A \emph{bubble bifurcation (with quadratic fold)} occurs at $(\alpha^*,\beta^*)$ in \eqref{eq:DDE} if there exists a Hopf curve $\Theta_H$ with $(\pi^\alpha u_{\Theta_H(t^*)},\pi^\beta u_{\Theta_H(t^*)})=(\alpha^*,\beta^*)$ for some $t^*\in(0,1)$, such that
\begin{itemize}
\item the projection of $\Theta_H$ into the $(\alpha,\beta)$ plane can be realized as a $C^2$ graph $\beta = \beta(\alpha)$.
\item $\beta(\alpha^*)=\beta^*$, $\beta'(\alpha^*)=0$, and $\beta''(\alpha^*)\neq 0$.
\item there is a $C^2$ diffeomorphism $h:D\rightarrow h(D)$ defined on a neighbourhood $D$ of $\Theta_H(t^*)\in\Delta$, such that $h(\Theta_H(t^*))=(\alpha^*,0)$, and the periodic orbit $y$ (see Definition \ref{def:hopfcurve}) is a steady state if and only if $\pi_2 h(s)=0$, where $\pi_2$ is the projection onto the second factor. Also, the projection $\pi_1$ onto the first factor satisfies $\pi_1 h(s) = \pi^a u_{s}$ for all $s\in D$.
\item $(\alpha^*,0)$ is a strict local extremum of the map $(\alpha,a)\mapsto \pi^\beta u_{h^{-1}(\alpha,a)}$.
\end{itemize}
\end{definition}

Our perspective is that a bubble bifurcation is characterized by the existence of a manifold of periodic orbits, parameterized in terms of amplitude and a parameter, such that at an isolated critical point, the projection of the manifold onto the other parameter has a local extremum. The parametrization of $\beta$ near $\alpha^*$ and the local minimum condition reflects the observation that a fold in the curve of Hopf bifurcations is sufficient condition for the birth of the bubble. That is, we have a family of bubbles (loops of Hopf bifurcations) that can be parameterized by $\beta$ in an interval of the form $(\beta^*-\delta,\beta^*]$ or $[\beta^*,\beta^*+\delta)$, for $\delta>0$ small. This is a consequence of the $(\alpha,a)$-parametrization of the simplex.

\begin{remark}
We include the adjective ``quadratic" to describe the fold in the Hopf curve to contrast with a related definition in Section \ref{sec-degenerate2}, where the condition $\beta''(\alpha^*)\neq 0$ will be weakened.
\end{remark}

\subsubsection{Preparations}
It will facilitate the presentation of the bubble bifurcation if we are able to construct the first and second derivatives of $s\mapsto u_s$, the zeroes of \eqref{Fmap}, with respect to the simplex parameter $s$. This was partially done in the proof of Proposition \ref{prop-Hopf}, but we will elaborate further here. 

To compute the derivatives, the idea is to formally differentiate $s\mapsto \F_s(u_s)$ with respect to $s$ twice, allowing for the introduction of $\partial_s u_s$ and $\partial_s^2 u_s$. The result is a system of three nonlinear equations, and solving for $(u_s,\partial_s u_s,\partial_s^2 u_s)$ amounts to computing zeroes of a map 
\begin{align}\label{non-deg-hopf-map}\Omega^3\ni (u_s,\partial_s u_s) \mapsto (\F_s(u_s),\F_s^{[1]}(u_s,\partial_s u_s),\F_s^{[2]}(u_s,\partial_s u_s,\partial_s^2 u_s))
\end{align}

\begin{remark}\label{rem-2derivatives}
It may be unclear whether the second derivative of $s\mapsto\F_s(u_s)$ can be given meaning. Indeed, as we have previously explained in Remark \ref{remark-D2F-DNE}, the second Fr\'echet derivative of $u\mapsto \F_s(u)$ does not, in general, exist. The (formal) second derivative of both sides of $\F_s(u_s)=0$ produces
$$\partial_s^2\F_s(u_s) + 2\partial_s [D\F_s](u_s)\partial_s u_s + D^2\F_s(u_s)[\partial_s u_s,\partial_s u_s]=0,$$
and we can see that it is in fact possible to interpret $D^2\F_s(u_s)[\partial_s u_s,\partial_s u_s]$ as the Gateaux derivative of $w\mapsto D\F_s(w)\partial_s u_s$, at $w=u_s$ and in the direction $\partial_s u_s$. However, even this may not exist as an element of $K_\nu(\C^{n+m})$, since the Fourier components of $u_s$ are generally not band-limited. It is therefore necessary to specify the codomain of \eqref{non-deg-hopf-map} a bit more carefully. We will not elaborate on this subtlety here; the ramifications of this will be the topic of some of our future work. However, in the case of \textit{ordinary differential equations}, there is no major complication. When there are no delays (or when they are all zero), $\F_s$ is twice continuously differentiable provided the same is true of $f$.
\end{remark}

Regardless how the partial derivatives are enclosed, the following equivalent version of Proposition \ref{prop-Hopf} is available.

\begin{proposition}\label{prop-hopf2}
Assume the existence of a family of zeroes of \eqref{non-deg-hopf-map} parameterized by $s\in\Delta$ and close to a numerical interpolant $\hat{\mathbf{u}}_s=(\hat u_s,\partial_s\hat u_s,\partial_s^2\hat u_s)$ has been proven, in the sense that we have identified $r>0$ such that there is a unique zero of \eqref{non-deg-hopf-map}, for each $s\in\Delta$, in the ball $B_r(\hat{\mathbf{u}}_s)$. Assume the topolgy on this ball is such that the components $\hat a_s=\pi^a \hat u_s$ satisfy $$||(\hat a_s^{[1]} - a_s^{[1]})h||\leq r_a^{(1)}|h|,$$ for $h\in\R^2$, and $||\hat a_s^{[0]} - a_s^{[0]}||\leq r_a^{(0)}$. Suppose
\begin{itemize}
\item For $j=1,2$, $[\hat a^{[0]}_{s_0}-r_a^{(0)},\hat a^{[0]}_{s_0}+r_a^{(0)}]\cdot[\hat a^{[0]}_{s_j}-r_a^{(0)},\hat a^{[0]}_{s_j}+r_a^{(0)}]<0$, where the multiplication is in the sense of interval arithmetic.
\item The sign of $[\hat a^{[0]}_{s}-r_a^{(0)},\hat a^{[0]}_{s}+r_a^{(0)}]$ is constant for $s$ on the edge of $\Delta$ not incident with $s_0$; that is, on the edge connecting $s_1$ to $s_2$.
\item The components of the interval vector $\hat a_s^{[1]} + r_a^{(1)}[-1,1]^2$ in $\R^{2*}$ has the same sign for $s\in\Delta$.
\end{itemize}
Then there exists a Hopf bifurcation curve $\Theta_H:[0,1]\rightarrow\Delta$, each of $\Theta_H(0)$ and $\Theta_H(1)$ lie on the edges of $\Delta$ that are incident with $s_0$, and $\Theta_H$ is $C^2$.
\end{proposition}

Comments analogous to those appearing in Remark \ref{rem-Hopf} apply here as well. Note that Proposition \ref{prop-hopf2} requires only the first derivatives, which can be enclosed using \eqref{derivative-us}. However, if the existence of the second derivatives are in question, then $\Theta_H$ can at best be proven to be $C^1$.

\begin{definition}
A triple of line segments $(\overline v_1, \overline v_*, \overline v_2)$ in $\Delta$ is \textit{$s_0$-oriented} if there exist points $t_1,t_*,t_2$ on the edge $(s_1,s_2)$, such that
\begin{itemize}
\item $\overline v_1$ is a subset of the line connecting $s_0$ to $t_1$;
\item $\overline v_*$ is a subset of the line connecting $s_0$ to $t_*$;
\item $\overline v_2$ is a subset of the line connecting $s_0$ to $t_2$;
\item under the total order on the edge $(s_1,s_2)$ defined by $a\leq b \Leftrightarrow ||s_1-a||_2\leq||s_1-b||_2$, we have $\overline v_1\leq\overline v_*\leq\overline v_2$.
\end{itemize} 
In this case the associated \textit{wedge cover} is the simplex in $\Delta$ with vertices $s_0,t_1$ and $t_2$, and it is denoted $W(\overline v_1,\overline v_*,\overline v_2)$.
\end{definition}
These $s_0$-oriented line segments will be used to enclose potential bubble bifurcation points. Some finer control can be specified with the following.

\begin{definition}
Let $(\overline v_1,\overline v_*,\overline v_2)$ be an $s_0$-oriented triple in $\Delta$. Let $\Theta_H$ be a Hopf curve in $\Delta$ associated to a family $u_s$ of zeroes of $\F_s$. A \emph{Hopf-bounding trapezoid in $W(\overline v_1,\overline v_*,\overline v_2)$} is a trapezoid contained in $W(\overline v_1,\overline v_*,\overline v_2)$ whose edges include both $\overline v_1$ and $\overline v_2$, and such that the other edges, denoted $\overline w_+$ and $\overline w_-$, satisfy $\pi^a u_s <0$ for $s\in \overline w_-$, and $\pi^a u_s>0$ for $s\in \overline w_+$. 
\end{definition}

It is straightforward to verify that, under the assumptions of Proposition \ref{prop-hopf2} 
, a trapezoid in $W(\overline v_1,\overline v_*,\overline v_2)$ is Hopf-bounding provided two of its edges match $\overline v_1$ and $\overline v_2$, while
\begin{align*}
\hat a_s^{[0]}+r_a^{(0)}&<0,\quad s\in\overline w_-\\
\hat a_s^{[0]}-r_a^{(0)}&>0,\quad s\in\overline w_+.
\end{align*} 

\subsubsection{Enclosure of a bubble bifurcation}
The following proposition provides verifiable conditions for the existence of a bubble bifurcation. Some are more explicit than others.

\begin{proposition}\label{prop-degenerate}
Let the hypotheses of Proposition \ref{prop-hopf2} be satisfied. Let $(\overline v_1,\overline v_*,\overline v_2)$ be an $s_0$-oriented triple of line segments in $\Delta$. Suppose the topology on $B_r(\hat{\mathbf{u}}_s)$ is such that the components $\pi^\alpha\hat u_s^{[k]} = \hat \alpha_s^{[k]}$,  $\pi^\beta\hat u_s^{[k]} = \hat \beta_s^{[k]}$, $\pi^a\hat u_s^{[k]} = \hat a_s^{[k]}$ satisfy
\begin{align*}
||(\hat\alpha^{[k]}_s - \alpha^{[k]}_s)[h_1,\dots,h_k]||&\leq r_\alpha^{(k)}|h_1|\cdots|h_k|\\
||(\hat\beta^{[k]}_s - \beta^{[k]}_s)[h_1,\dots,h_k]||&\leq r_\beta^{(k)}|h_1|\cdots|h_k|\\
||(\hat a^{[k]}_s - a^{[k]}_s)[h_1,\dots,h_k]||&\leq r_a^{(k)}|h_1|\cdots|h_k|
\end{align*}
for $k$-tuples $h_1,\dots,h_k\in\R^2$, and $k=0,1,2$. With the empty tuple ($k=0$), the norm reduces to absolute value on $\R$. Let $\mathbf{r}_\kappa$, for $\kappa\in\{\alpha,\beta,a\}$ be interval vectors in $\R^2$ such that $(\mathbf{r}_\kappa)_j = [-1,1]r_\kappa^{(1)}$, for $j=1,2$.  Finally, given $s\in\Delta$, denote by $V_s$ the set
$$V_s = \{v\in\R^2 : (\hat a_s^{[1]}+\mathbf{r}_a)\cdot v=0,\hspace{2mm}||v||_2=1\}.$$
Suppose the following are satisfied.
\begin{enumerate}
\item $(\hat\alpha_s^{[1]} + \mathbf{r}_\alpha)\cdot v_s$ is bounded away from zero, for all $s\in \Delta$ and all $v_s\in V_s$.
\item $\inf \hat a_{\overline v_i}^{[0]}+r_a^{(0)}<0<\sup \hat a_{\overline v_i}^{[0]}-r_a^{(0)}$ for $i=1,2$.
\item $\inf \hat a_{\overline v_*}^{[0]}+r_a^{(0)}<0<\sup \hat a_{\overline v_*}^{[0]}-r_a^{(0)}$.
\item $\hat \beta^{[0]}_{\overline v_*}+r_\beta^{(0)} < \hat \beta^{[0]}_{\overline v_i}-r_\beta^{(0)}$ for $i=1,2$.
\item The components of $\hat a_s^{[1]}+\mathbf{r}_a$ are bounded away from zero for all $s\in \Delta$.
\item The determinant of the $2\times 2$ interval matrix
\begin{align}\label{invertible-matrix}\left(\begin{array}{cc}(\hat \alpha^{[1]}_s)_1 + [-1,1]r_\alpha^{(1)} & (\hat \alpha^{[1]}_s)_2 + [-1,1]r_\alpha^{(1)} \\ (\hat a^{[1]}_s)_1 + [-1,1]r_a^{(1)} & (\hat a^{[1]}_s)_2 + [-1,1]r_a^{(1)}  \end{array}\right)\end{align} is bounded away from zero for $s\in \Delta$.
\item Defining
\begin{align}\label{cs}
c_s = \frac{1}{||\hat a^{[1]}_s + \mathbf{r}_a||^2}(\hat \beta^{[1]}_s+\mathbf{r}_\beta)\cdot(\hat a^{[1]}_s + \mathbf{r}_a),
\end{align}
it holds that the interval $\hat\beta^{[2]}_s[v_s,v_s] - c_s\hat a^{[2]}_s[v_s,v_s] + (r_\beta^{(2)} + ||c_s||_2r_a^{(2)})[-1,1]$ is bounded away from zero for all $s\in \Delta$ and $v_s\in V_s$.
\item The matrix
\begin{align}\label{matrix-Hess}
\left(\begin{array}{cc}
\beta_s^{[2]}[\partial_\alpha s,\partial_\alpha s] + \beta_s^{[1]}\partial_\alpha^2 s & \beta_s^{[2]}[\partial_a s,\partial_\alpha s] + \beta_s^{[1]}\partial_a\partial_\alpha s \\ 
\beta_s^{[2]}[\partial_a s,\partial_\alpha s] + \beta_s^{[1]}\partial_a\partial_\alpha s & \beta_s^{[2]}[\partial_a s,\partial_a s] + \beta_s^{[1]}\partial_a^2 s
\end{array}\right)
\end{align}
is uniformly (for $s\in\Delta$) definite, where defining $h:\Delta\rightarrow\R^2$ by $h(s)=(\alpha_s,a_s)$, the partial derivatives in the matrix above are
\begin{align*}
\left[\begin{array}{cc}\partial_\alpha s & \partial_a s \end{array}\right]&=Dh(s)^{-1}\\
\partial_\alpha^2 s&=-Dh(s)^{-1}D^2 h(s)[\partial_\alpha s,\partial_\alpha s]\\
\partial_a\partial_\alpha s&=-Dh(s)^{-1}D^2 h(s)[\partial_a s,\partial_\alpha s]\\
\partial_a^2 s&=-Dh(s)^{-1}D^2 h(s)[\partial_a s,\partial_a s].
\end{align*}
\end{enumerate}
Then, there exists a bubble bifurcation with quadratic fold at some $(\alpha^*,\beta^*)$ in the projection of $\Theta_H\cap W(\overline v_1,\overline v_*,\overline v_2)$ onto the $(\alpha,\beta)$ plane.
\end{proposition}

\begin{proof}
By a suitable reparametrization, we may assume that $||\Theta_H'(t)||=1$ for all $t\in(0,1)$. Denote $a(s)=\pi^a u_s$. The definition of the Hopf curve is that $a(\Theta_H(t))=0$ or all $t\in[0,1]$. As consequence, $(\partial_s a)\Theta_H'=0$, so that $\Theta_H'$ is dual to the orthogonal complement of $\partial_s a$. That is, $\Theta_H'\in \cup_s V_s$. Now,
$$\frac{d}{dt}\alpha^{[0]}_{\Theta_H(t)} = \alpha_{\Theta_H(t)}^{[1]}\Theta_H'\in\bigcup_{s\in \Delta}\left( \bigcup_{v_s\in V_s} (\hat\alpha_s^{[1]}+\mathbf{r}_\alpha)\cdot v_s \right),$$ 

which is bounded away from zero by assumption (condition 1 of the proposition). It follows that $t\mapsto\alpha_{\Theta_H(t)}^{[0]}$ is monotone, so $\Theta_H$ can be parameterized by $\alpha$, for $\alpha$ in the monotone range of $t\mapsto\alpha_{\Theta_H(t)}^{[0]}$. Consequently, the projection of the Hopf curve $\Theta_H$ in the $(\alpha,\beta)$ plane can be represented as a graph $\beta=\beta(\alpha)$.

Conditions 2--4 of the proposition guarantee that each of the segments $\overline s_1$, $\overline s_2$ and $\overline t$ enjoy the following properties:
\begin{itemize}
\item as one-dimensional manifolds, they enclose an intersection with the Hopf curve;
\item the $\beta$-components of $\overline v_1$ and $\overline v_2$ are strictly greater than those of $\overline v_*$.
\end{itemize}
As consequence, $\beta=\beta(\alpha)$ possesses an internal (to its domain, relative to the previously-computed range) critical point which is a global minimizer. Let this point be $\alpha^*$, so $\beta(\alpha^*)\bydef\beta^*$. Let the associated point on the simplex be $\Theta_H(t^*)$. 

We wish to show that $\beta^*$ is a strict, isolated minimum of $\beta$. We will do that by proving $\beta''(\alpha^*)\neq 0$. It is enough to prove that $\frac{d^2}{dt^2}\beta^{[0]}_{\Theta_H(t)}\neq 0$ whenever $\frac{d}{dt}\beta^{[0]}_{\Theta_H(t)}=0$. If the latter is satisfied, then we have simultaneously
\begin{align*}
a_{\Theta_H(t)}^{[1]}\Theta_H'&=0\\
\beta_{\Theta_H(t)}^{[1]}\Theta_H'&=0.
\end{align*}
Since $\Theta_H'\neq 0$, it must be the case that $\alpha^{[1]}_{\Theta_H(t)}$ and $\beta_{\Theta_H(t)}^{[1]}$ are colinear. By assumption 5, $\alpha^{[1]}_{\Theta_H(t)}$ is bounded away from zero, so the quantity $c_s$ of \eqref{cs} is well-defined and \begin{align}\label{beta-thetaH}\beta_{\Theta_H(t)}^{[1]}\in c_{\Theta_H(t)} a_{\Theta_H(t)}^{[1]}.\end{align}
On its own, \eqref{beta-thetaH} might not seem particular useful. However, consider that
\begin{align*}
\frac{d^2}{dt^2}\beta^{[0]}_{\Theta_H(t)}&=\beta_{\Theta_H(t)}^{[2]}[\Theta_H',\Theta_H'] + \beta_{\Theta_H(t)}^{[1]}\Theta_H''\\
0&=a_{\Theta_H(t)}^{[2]}[\Theta_H',\Theta_H'] + a_{\Theta_H(t)}^{[1]}\Theta_H'',
\end{align*}
where the second equation comes from taking a second derivative of the definition $a(\Theta_H(t))=0$ of the Hopf curve. Using the second equation together with \eqref{beta-thetaH}, we can get the inclusion
$$\beta_{\Theta_H(t)}^{[1]}\Theta_H'' \in -c_{\Theta_H(t)}a_{\Theta_H(t)}^{[2]}[\Theta_H',\Theta_H'],$$ thereby removing the dependence on the second derivative $\Theta_H''$ of the Hopf curve. Substituting into the expression for the second derivative of $t\mapsto\beta_{\Theta_H(t)}$, we get
\begin{align*}
\frac{d^2}{dt^2}\beta^{[0]}_{\Theta_H(t)}&\in \beta_{\Theta_H(t)}^{[2]}[\Theta_H',\Theta_H'] -c_{\Theta_H(t)}a_{\Theta_H(t)}^{[2]}[\Theta_H',\Theta_H']\\
&\subset \bigcup_{s\in \Delta}\left(\bigcup_{v_s\in V_s} \hat\beta^{[2]}_s[v_s,v_s] - c_s\hat a^{[2]}_s[v_s,v_s] + (r_\beta^{(2)} + ||c_s||_2r_a^{(2)})[-1,1]\right),
\end{align*}
which is bounded away from zero by condition 7. Therefore, $\beta''(\alpha^*)\neq 0$.

Next, we need to verify the local parametrization of the simplex near $\Theta_H(t^*)$ in terms of $(\alpha,a)$. This is a fairly direct consequence of the inverse function theorem, using the condition 6 of the proposition. This shows that the function $h$ defined in condition 8 of the proposition defines a local diffeomorphism near $\Theta_H(t^*)$. The Hessian of $\Gamma:(\alpha,a)\mapsto \beta_{h^{-1}(\alpha,a)}$ can be computed by implicit differentiation; if $h(s)=(\alpha,a)$, then the Hessian is precisely the $2\times 2$ matrix of condition 8. Since this matrix is uniformly definite on the simplex $\Delta$, every critical point must be a local extremum. We already know that $\partial_\alpha\Gamma(\alpha^*,0)=\beta'(\alpha^*)=0$. The other partial derivative $\partial_a\Gamma(\alpha^*,0)$ is zero due to the amplitude symmetry of periodic orbits. Therefore, $(\alpha^*,0)$ is a strict local extremum of the map $(\alpha,a)\mapsto \pi^\beta u_{h^{-1}(\alpha,a)}$.
\end{proof}

\begin{remark}
Checking conditions 2--6 of the proposition is fairly routine. Conditions 1 and 7, however, deserve some extra attention. If the step size is small, we should expect the \textit{derivatives} of the solution $s\mapsto u_s$ to be close to constant. In this way, the set $V_s$ should not vary too much (in a Hausdorff sense). $V_s$ geometrically consists of two arcs on the unit circle, and we expect the angles defining these arcs to be nearly constant over the simplex. Consequently, the interval computations of conditions 1 and 7 are, indeed, implementable, but the \textit{feasibility} of the checks --- that these intervals are bounded away from zero --- will be strongly influenced by the size of the radius, $r$, and any weighting in the norm. As for condition 8, we have not included all of the associated radii, but since all of the intermediary derivatives appearing in the matrix \eqref{matrix-Hess} admit (by assumption) rigorous enclosures, the matrix is implementable. Therefore, condition 8 can be checked using a suitable package to compute eigenvalues of interval matrices.
\end{remark}

\begin{remark}
Condition 4 is formulated in such a way that $\beta^*$ is a local \textit{minimum} of the curve $\beta=\beta(\alpha)$. This condition can be reformulated in a straightforward way to allow it instead to be a local maximum. Also, comments analogous to those appearing in Remark \ref{rem-Hopf} apply here as well. 
\end{remark}

\begin{corollary}
It suffices to verify the conditions 6,7,8 of Proposition \ref{prop-degenerate} for all $s$ in a given Hopf-bounding trapezoid in $W(\overline v_1,\overline v_*,\overline v_2)$.
\end{corollary}


\subsubsection{A degenerate Hopf bifurcation without second derivatives}\label{sec-degenerate2}
In Remark \ref{rem-2derivatives}, we pointed out that, unfortunately, the lack of second-differentiability of the map $\F_s$ is a serious obstruction to computing second derivatives $\partial_s^2 u_s$ and, consequently, checking all the conditions of Proposition \ref{prop-degenerate}. While it is not a problem when all of the delays are zero (i.e.\ an ODE), we would like to provide a constructive result for delay equations. Along these lines, let us slightly weaken Definition \ref{def-degenerate}.

\begin{definition}\label{def-generateHopf2}
A \emph{degenerate Hopf bifurcation} occurs at $(\alpha^*,\beta^*)$ in \eqref{eq:DDE} if there exists a Hopf curve $\Theta_H$ with $(\pi^\alpha u_{\Theta_H(t^*)},\pi^\beta u_{\Theta_H(t^*)})=(\alpha^*,\beta^*)$ for some $t^*\in(0,1)$, such that
\begin{itemize}
\item the projection of $\Theta_H$ into the $(\alpha,\beta)$ plane can be realized as a $C^1$ graph $\beta = \beta(\alpha)$;
\item $\beta(\alpha^*)=\beta^*$ and $\beta'(\alpha^*)=0$;
\item there is a $C^1$ diffeomorphism $h:D\rightarrow h(D)$ defined on a neighbourhood $D$ of $\Theta_H(t^*)\in\Delta$, such that $h(\Theta_H(t^*))=(\alpha^*,0)$, and the periodic orbit $z$ (see Definition \ref{def:hopfcurve}) is a steady state if and only if $\pi_2 h(s)=0$, where $\pi_2$ is the projection onto the second factor. Also, the projection $\pi_1$ onto the first factor satisfies and $\pi_1 h(s) = \pi^a u_{s}$ for all $s\in D$.
\end{itemize}
\end{definition}


The main difference between the above and Definition \ref{def-degenerate} is we no longer require that $\beta^*$ is an extremum of $\beta$. We also do not make any specifications concerning the geometry of the implicit map $(\alpha,a)\mapsto\beta$ near the Hopf bifurcation curve. Clearly, a bubble bifurcation with quadratic fold satisfies the conditions of the above definition. However, the new definition permits other types of degenerate Hopf bifurcations, including Bautin bifurcations. 
The following proposition can now be proven using the same ideas as Proposition \ref{prop-degenerate}.


\begin{proposition}\label{prop-degenerate2}
Assume a family of zeroes $u_s$ of the map $\F_s$ has been proven, in addition to the first derivatives $\partial_s u_s$, close to a numerical interpolant $\hat{\mathbf{u}}_s=(\hat u_s,\partial_s\hat u_s)$, in the sense that we have identified $r>0$ such that there is a unique zero of \eqref{non-deg-hopf-map}, for each $s\in\Delta$, in the ball $B_r(\hat{\mathbf{u}}_s)$. Suppose the topology on $B_r(\hat{\mathbf{u}}_s)$ is such that the components $\pi^\alpha\hat u_s^{[k]} = \hat \alpha_s^{[k]}$,  $\pi^\beta\hat u_s^{[k]} = \hat \beta_s^{[k]}$, $\pi^a\hat u_s^{[k]} = \hat a_s^{[k]}$ satisfy
\begin{align*}
||(\hat\alpha^{[k]}_s - \alpha^{[k]}_s)[h_1,\dots,h_k]||&\leq r_\alpha^{(k)}|h_1|\cdots|h_k|\\
||(\hat\beta^{[k]}_s - \beta^{[k]}_s)[h_1,\dots,h_k]||&\leq r_\beta^{(k)}|h_1|\cdots|h_k|\\
||(\hat a^{[k]}_s - a^{[k]}_s)[h_1,\dots,h_k]||&\leq r_a^{(k)}|h_1|\cdots|h_k|
\end{align*}
for $k$-tuples $h_1,\dots,h_k\in\R^2$, and $k=0,1$. With the empty tuple ($k=0$), the norm reduces to absolute value on $\R$. Let $(\overline v_1,\overline v_*,\overline v_2)$ be an $s_0$-oriented triple of line segments in $\Delta$. Suppose conditions 1--3, 5 and 6 of Proposition \ref{prop-degenerate} are satisfied. Then, there exists a degenerate Hopf bifurcation at some $(\alpha^*,\beta^*)$ in the projection of $\Theta_H\cap W(\overline v_1,\overline v_*,\overline v_2)$ onto the $(\alpha,\beta)$ plane. 
\end{proposition}


\begin{corollary}
It suffices to verify condition 6 of Proposition \ref{prop-degenerate2} for all $s$ in a Hopf-bounding trapezoid in $W(\overline v_1,\overline v_*,\overline v_2)$.
\end{corollary}

\section{Examples}\label{sec:examples}
\subsection{Extended Lorenz-84 system}
The extended Lorenz-84 system is the following system of four coupled ODEs:
\begin{align*}
\dot u_1&=-u_2^2 - u_3^2 - au_1 - af - bu_4^2\\
\dot u_2&= u_1u_2 - cu_1u_3 - u_2 + d\\
\dot u_3&= cu_1u_2 + u_1u_3 - u_3\\
\dot u_4&= -eu_4 + bu_4u_1 + \mu
\end{align*}
We consider the parameters $a = 0.25$, $b=0.987$, $d=0.25$, $e=1.04$, $f=2$ 
to be fixed, while $\mu$ and $c$ are treated as parameters.

\begin{figure}[!hbt]
  \begin{minipage}[c]{0.67\textwidth}
    \includegraphics[width=0.9\textwidth]{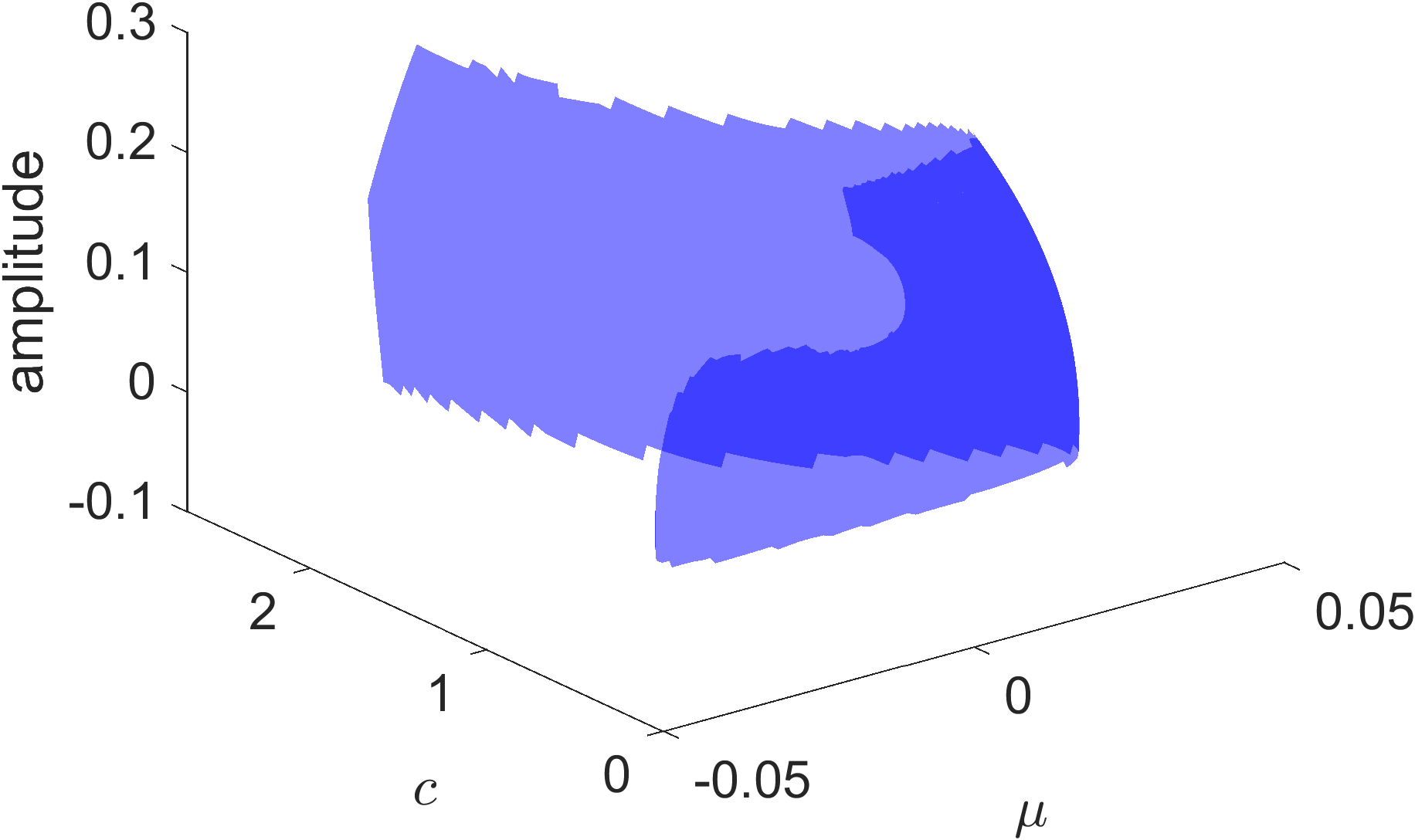}
  \end{minipage}\hfill
  \begin{minipage}[c]{0.33\textwidth}
    \caption{
       The manifold of (proven) periodic orbits in the extended Lorenz-84 system, in the projection of amplitude and the parameters $c,\mu$. Note the square-root curvature of amplitude near $\mu\approx 0.04$ (right side of plot), indicative of Hopf bubbles. There are over 10,000 simplices, so to provide a clean figure, we have not plotted the edges. 
    } \label{fig:Lorenz1}
  \end{minipage}
\end{figure}
\begin{figure}[!hbt]
  \begin{minipage}[c]{0.6\textwidth}
    \includegraphics[width=0.85\textwidth]{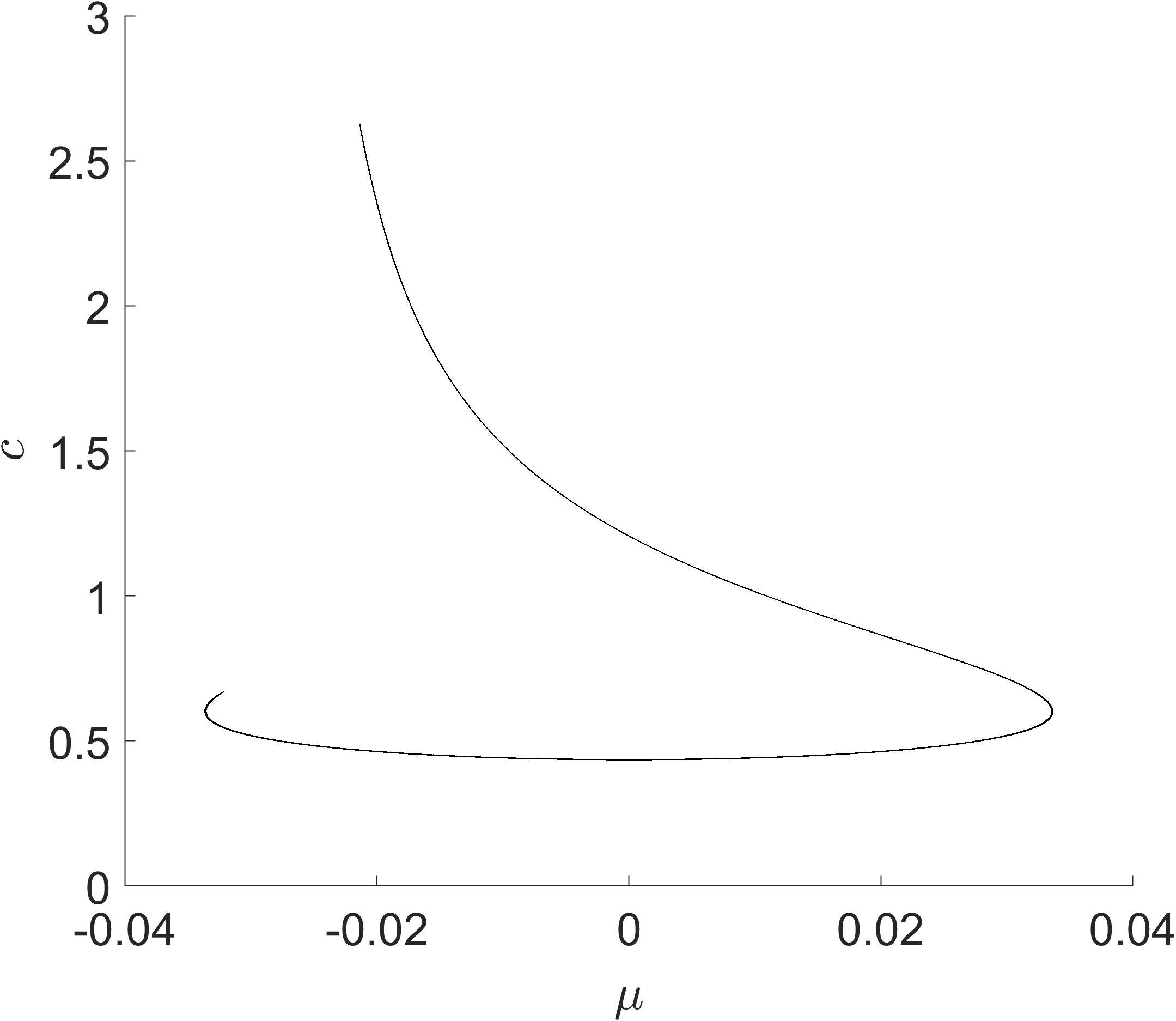}
  \end{minipage}\hfill
  \begin{minipage}[c]{0.4\textwidth}
    \caption{
       Intersection of the surface in Figure \ref{fig:Lorenz1} with the amplitude zero plane, which corresponds to the Hopf bifurcation curve. Hopf bubbles are present in the former figure near $\mu\approx 0.04$, so the fold in the present figure likely corresponds to a bubble bifurcation with quadratic fold. There is a symmetric fold near $\mu\approx-0.04$, and with respect to the other variable near $\mu\approx 0$, for a total of three bubble bifurcations.
    } \label{fig:Lorenz2}
  \end{minipage}
\end{figure}

We started the continuation at $c=1$, close to a Hopf bifurcation. Using a step size of $0.02$, we generated an approximate triangulation of the manifold with 11928 simplices (including simplices created by adaptive refinement needed for proofs). We used $N=5$ Fourier modes. To capture a ``wider" section of the manifold, we restricted the simplex growth to amplitude in the interval $[-0.1,0.3]$. Figure \ref{fig:Lorenz1} is a plot of the triangulation, projected into amplitude and parameter space, while we restricted to the zero amplitude plane in Figure \ref{fig:Lorenz2} to generate a plot of the Hopf bifurcation curve. The former figure allows for visualization of the traditional square-root amplitude curvature near the Hopf bifurcation curve. Interesting, far from the bifurcation curve, there appears to be a near-circular ``hole" in the manifold. We have not studied its structure in detail, and have no insight into its significance. The restriction of the amplitude to $[-0.1,0.3]$ results in the top and bottom edges appearing ``ragged", since simplices can not organically grow to produce hexagon patches. The latter Figure  \ref{fig:Lorenz2} indicates the presence of three bubble bifurcations. 

\subsection{Time-delay SI model}\label{sec:SI_model}
We consider the time-delay SI model
\begin{align*}
\dot y(t)&=-y(t) + R_0e^{-py(t-\tau)}y(t)(1-y(t)),
\end{align*}
for which there is an analytical proof of a bubble bifurcation \cite{LeBlanc2016} near $(R_0,p)\approx(2.1474, 1.6617)$, for delay $\tau=10$. We will replicate this analysis using our rigorous two-parameter continuation. 

\subsubsection{Set-up}
We begin by desingularizing the vector field. Writing $y = x + a\tilde y$, for $x$ being a steady state, we get, dropping the tilde,
\begin{align*}
\dot y= -y + g(a,y_\tau p)R_0e^{-px}x(1-x) + R_0e^{-px}e^{-apy_\tau}(-a y^2 + y(1-2x)),
\end{align*}
where $y_\tau = y(t-\tau)$, and $g$ is defined by
\begin{align*}
g(a,y)&=\left\{\begin{array}{ll}a^{-1}(e^{-ay }-1),&a\neq 0 \\ -y,&a=0.\end{array}\right.
\end{align*}
Observe, $\partial_y g(a,y)=-e^{-ay}$ and $$g(a,y)=-y\sum_{k=1}^\infty \frac{1}{k!}(-ay)^{k-1}.$$
$g$ is indeed analytic. Whenever $g$ (or its derivatives) must be rigorously evaluated, we construct Taylor polynomials of sufficiently high order and propagate error from the remainder accordingly. 

Now we polynomialize. Define $z_2=e^{-a py}$ and $z_1=g(a,yp)$. Then
\begin{align*}
\dot z_1&
= -pz_2\big(-y + R_0e^{-px}z_1(t-\tau)x(1-x) + R_0e^{-px}z_2(t-\tau)(-a y^2 + y(1-2x)) \big)\\
\dot z_2&
= -a p z_2\big(-y + R_0e^{-px}z_1(t-\tau)x(1-x) + R_0e^{-px}z_2(t-\tau)(-a y^2 + y(1-2x)) \big)
\end{align*}
They can be more compactly written as $\dot z_1 = -z_2\dot z_0$ and $\dot z_2=-apz_2\dot z_0$. We also have the implied boundary conditions
\begin{align*}
z_1(0)&=g(a,z_0(0)p)\\
z_2(0)&=e^{-ap z_0(0)},
\end{align*}
where $z_0=y$. We need two unfolding parameters to compensate for the two extra boundary conditions.
\begin{lemma}
If $z$ is a periodic solution of 
\begin{align*}
\dot z_0(t)&=-z_0(t) + R_0\mu z_1(t-\tau)x(1-x) + R_0\mu z_2(t-\tau)(-az_0(t)^2 + z_0(t)(1-2x))\\
\dot z_1(t)&=-pz_2(t)\dot z_0(t) + \eta_1\\
\dot z_2(t)&=-apz_2(t)\dot z_0(t) + \eta_2
\end{align*}
for some $\eta_1,\eta_2\in\R$, and $z$ satisfies $z_1(0)=g(a,z_0(0)p)$ and $z_2(0)=e^{-apz_0(0)}$, then $\eta_1=\eta_2=0$.
\end{lemma}

\begin{proof}
First, suppose $\eta_2\neq 0$. Then necessarily, $z_2$ has constant sign because the differential equation for $z_2$ is affine-linear and $\eta_2\neq 0$. Since $z_2(0)>0$, we must have $z_2>0$. But this means that $$\frac{\dot z_2(t)}{z_2(t)} + ap\dot z_0(t) = \frac{\eta_2}{z_2(t)},$$ a contradiction, since the left side is periodic and $\eta_2\neq 0$. We may therefore assume that $\eta_2=0$. Then $\dot z_2(t) = apz_2(t)\dot z_0(t)$, and it follows again that $z_2>0$. But this means $$\frac{\dot z_1(t)}{z_2(t)} + p\dot z_0(t) = \frac{\eta_1}{z_2(t)},$$ and since the left-hand side is periodic, it follows that $\eta_1=0$.
\end{proof}

To complete the polynomial embedding, we further polynomialize the parameters. This is done to ensure compatibility with the numerical implementation. We define $\mu = R_0e^{-px}$, so that the complete polynomialized vector field is
\begin{align*}
\dot z_0(t)&=-z_0(t) + \mu z_1(t-\tau)x(1-x) + \mu z_2(t-\tau)(-az_0(t)^2 + z_0(t)(1-2x))\\
\dot z_1(t)&=-pz_2(t)\dot z_0(t) + \eta_1\\
\dot z_2(t)&=-apz_2(t)\dot z_0(t) + \eta_2.
\end{align*}
The complete set of boundary conditions is
\begin{align*}
0&=z_1(0) - g(a,z_0(0)p)\\
0&=z_2(0) - e^{-apz_0(0)}\\
0&=\mu - R_0e^{-px}
\end{align*}
In the terminology of Remark \ref{rem-embedding-dimension}, the embedding dimension is $m=3$. The steady-state equation is scalar, and is given by
$$0=-x + \mu x(1-x).$$
\subsubsection{Results}
\begin{figure}
  \begin{minipage}[c]{0.6\textwidth}
    \includegraphics[width=0.9\textwidth]{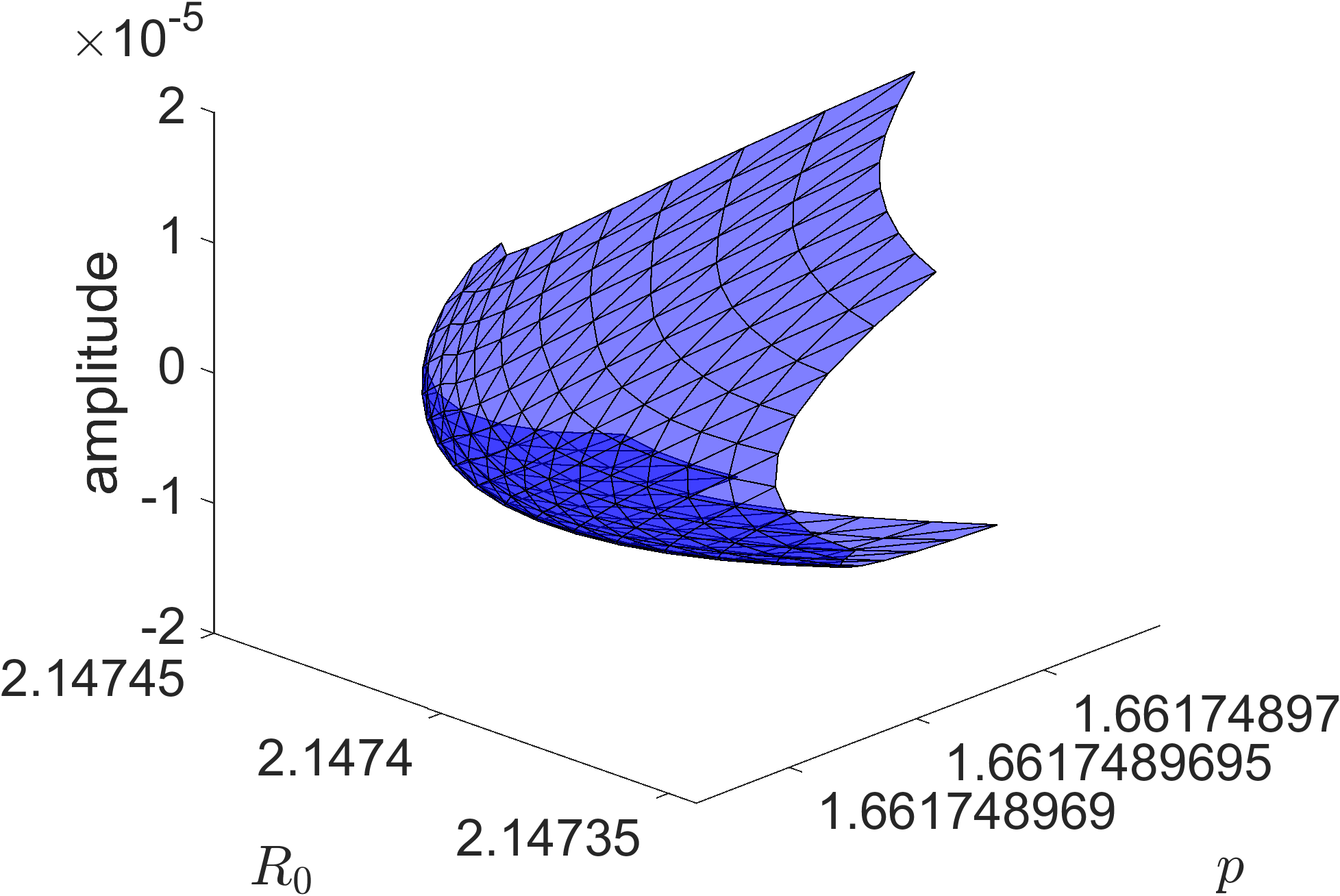}
  \end{minipage}\hfill
  \begin{minipage}[c]{0.4\textwidth}
    \caption{
       The manifold of (proven) periodic orbits in the time delay SI model. A very small step size was necessary to get proofs without too many adaptive refinement steps. The $Y_0$ bound is the clear bottleneck.
    } \label{fig:SI-manifold}
  \end{minipage}
\end{figure}

\begin{figure}[!htb]
  \begin{minipage}[c]{0.6\textwidth}
    \includegraphics[width=0.9\textwidth]{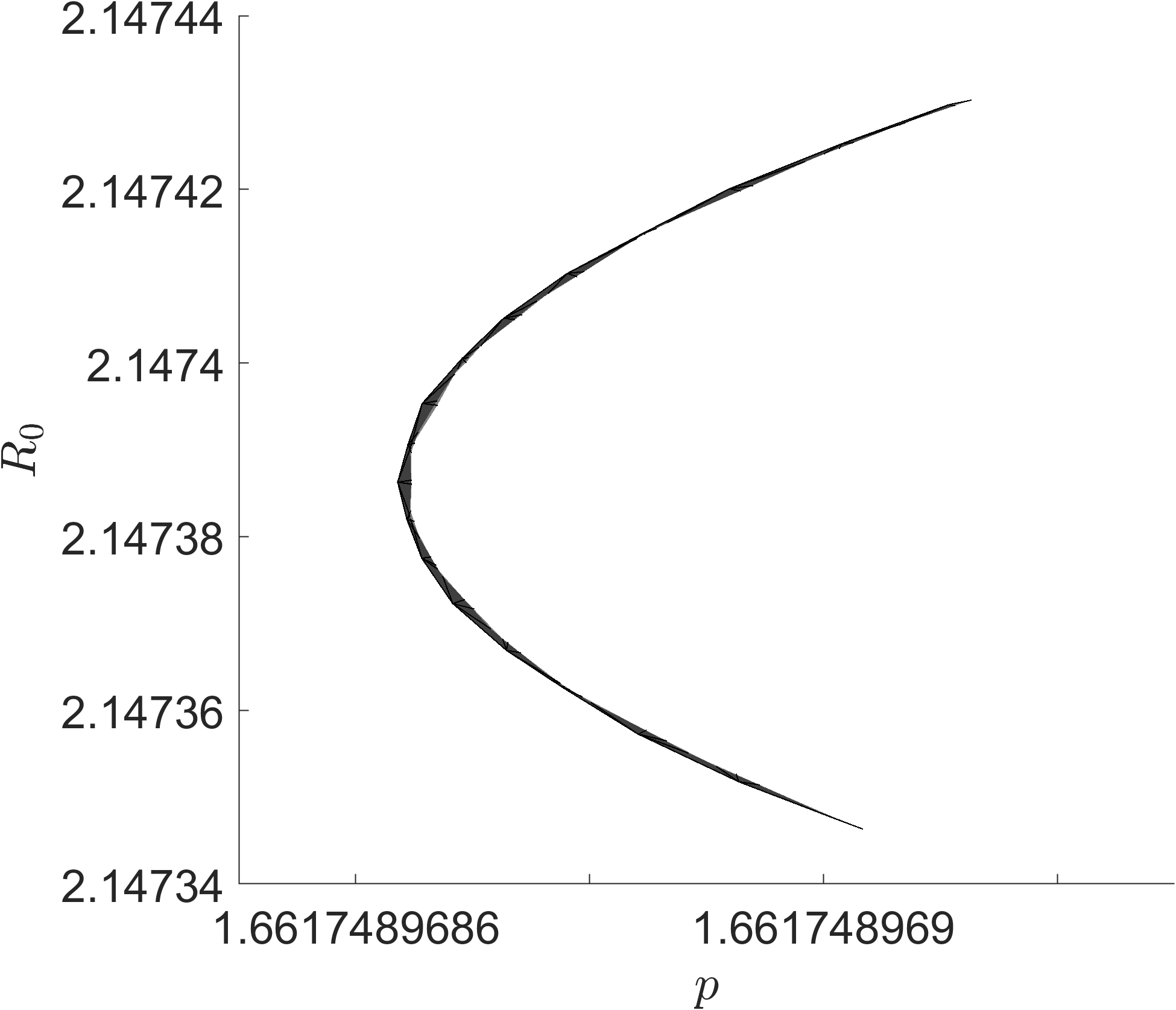}
  \end{minipage}\hfill
  \begin{minipage}[c]{0.4\textwidth}
    \caption{
       Intersection of the surface from Figure \ref{fig:SI-manifold} with the amplitude zero, which corresponds to a Hopf bifurcation curve in the time-delay SI model. The validation radius is $2\times 10^{-5}$ over the entire manifold, so the location $R_0\approx 2.1474$ of the bubble bifurcation is consistent with the analysis of Leblanc. A tighter validation radius could be obtained with a smaller step size. Because of how the curve is plotted, skew simplices from them 2-manifold have the effect of making the curve appear ``thicker" in some parts than others.
    } \label{fig:SI-curve}
  \end{minipage}
\end{figure}

We validated a patch of manifold initially with 406 simplices at a step size of $5\times 10^{-6}$. In the validation of nearly every simplex, three layers of adaptive refinement were needed to keep the $Y_0$ bound under control. We have plotted the manifold in Figure \ref{fig:SI-manifold} without the refinements included. The projection into the $(R_0,p)$ plane is provided in Figure \ref{fig:SI-curve}. The results are consistent with the analysis of Leblanc \cite{LeBlanc2016}.

\subsection{FitzHugh-Nagumo equation}
\begin{figure}[!htb]
\centering\includegraphics[scale=0.43]{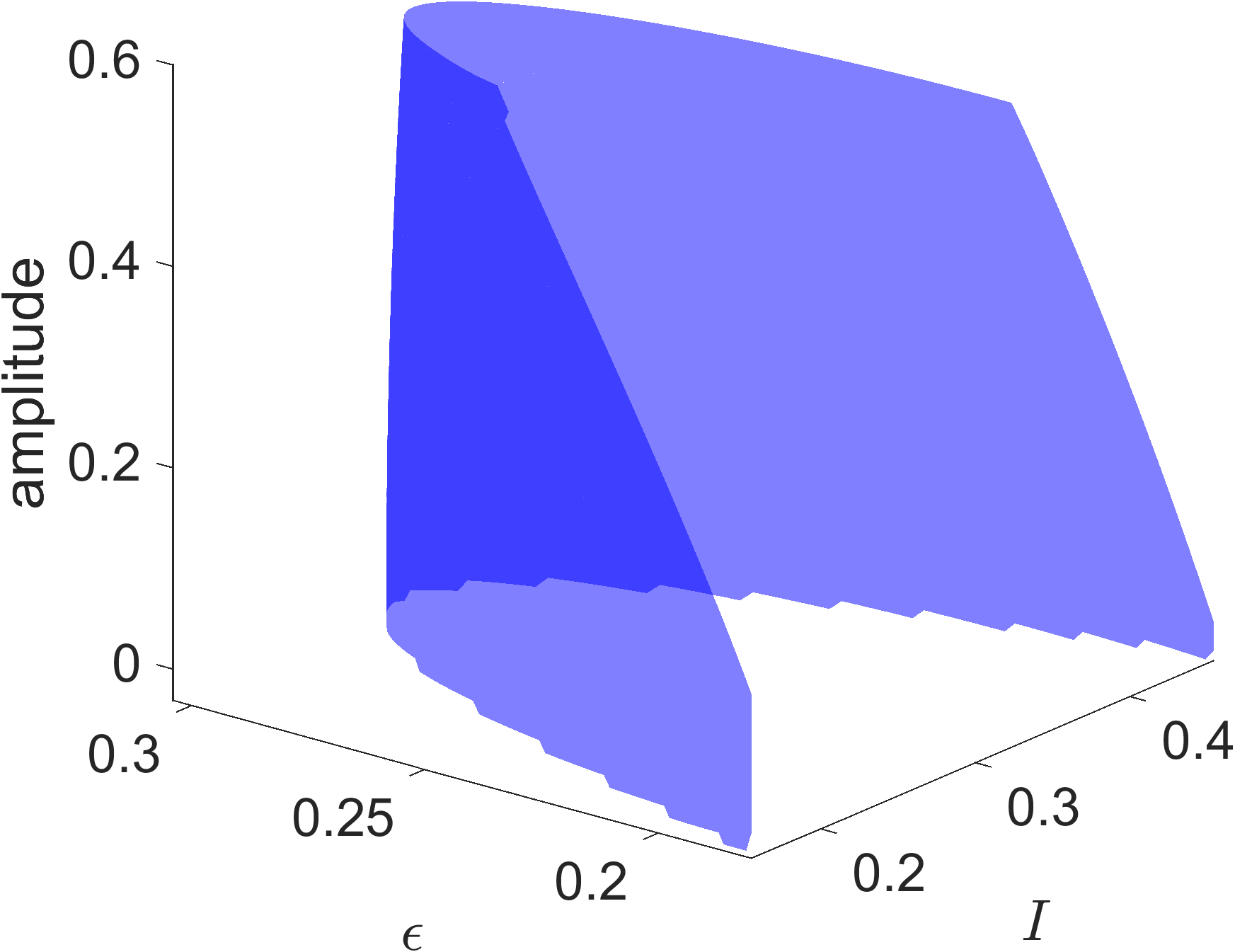} \hspace{3mm} \includegraphics[scale=0.43]{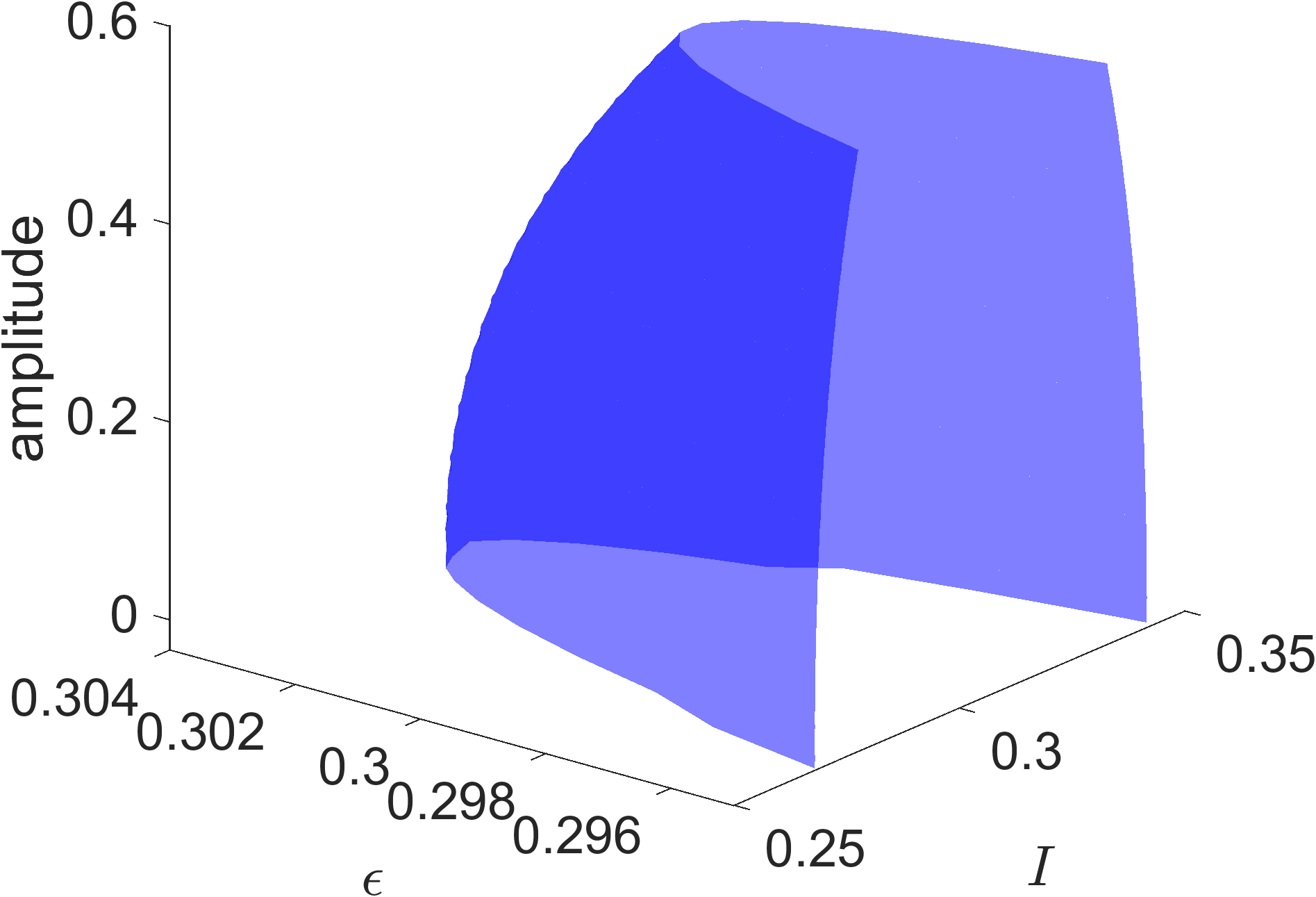}
\caption{Left: Projection of the manifold of periodic orbits for the FitzHugh-Nagumo equation into  $\epsilon\times I\times \mbox{amplitude}$ space. At this level of scaling, the curvature of the manifold is not easily visible. Right: zoomed-in portion near the bubble bifurcation. Here, the curvature is more easily seen. There are over 9000 simplices in the left figure, so we have not plotted the edges.}\label{fig:FHN0}
\end{figure}

\begin{figure}[!htb]
  \begin{minipage}[c]{0.6\textwidth}
    \includegraphics[width=0.9\textwidth]{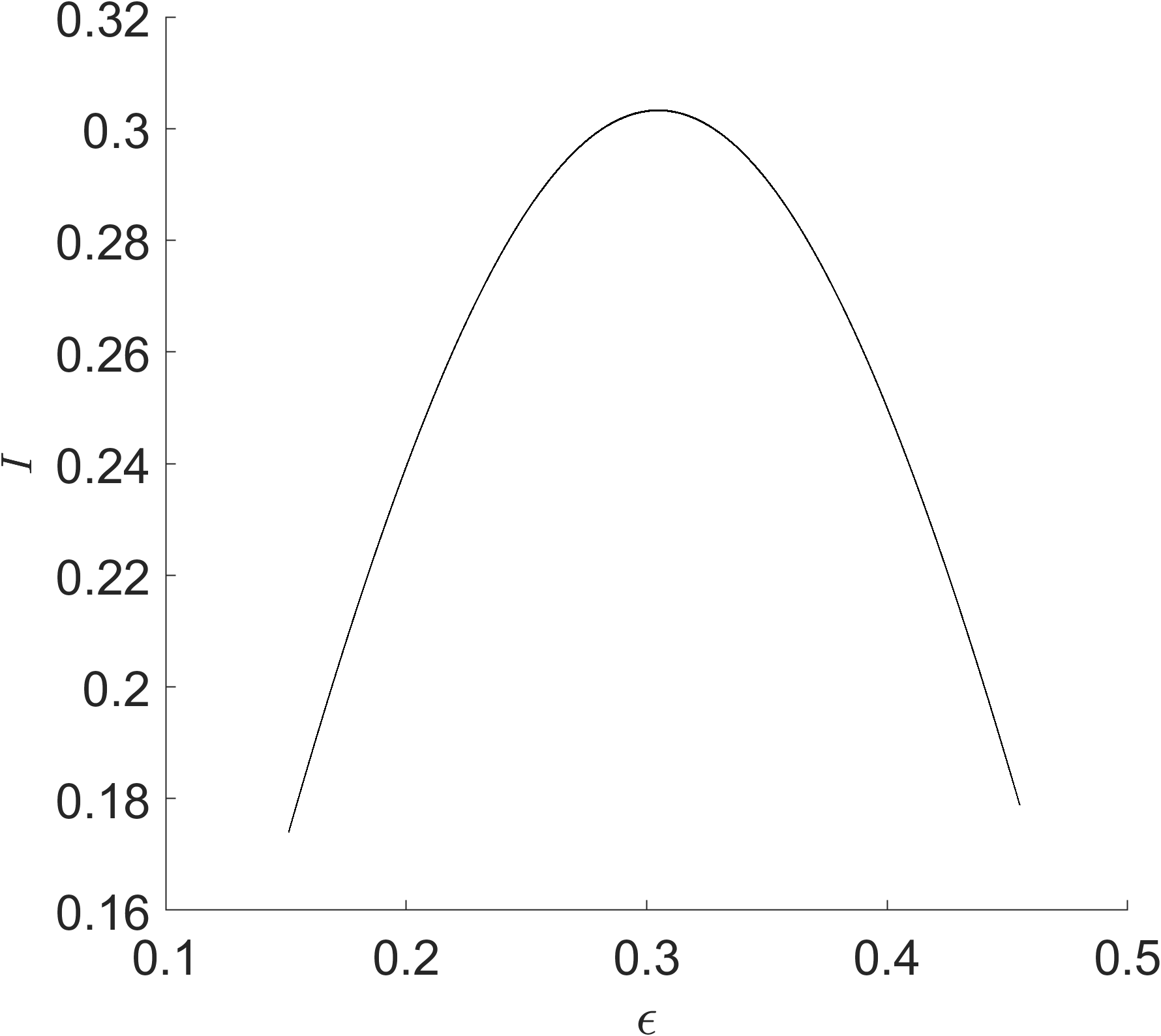}
  \end{minipage}\hfill
  \begin{minipage}[c]{0.4\textwidth}
    \caption{
       Intersection of the surface from the left pane of Figure \ref{fig:FHN0} with the amplitude zero plane. This corresponds to the Hopf bifurcation curve for the Fitzhugh-Nagumo system. The curve resembles a parabola, with a bubble bifurcation at its vertex. 
    } \label{fig:FHN1}
  \end{minipage}
\end{figure}

The FitzHugh-Nagumo ODE is
\begin{align*}
\dot u&=u(u-\alpha)(1-u)-w+I\\
\dot w&=\epsilon(u-\gamma w)
\end{align*}
for scalar parameters $\alpha,\epsilon,\gamma,I$. It is a cubic vector field in the state variables, and numerical simulations suggest the existence of Hopf bubbles (see Section 5.8 of \cite{Dupont2016}). We fix $\alpha = 0.1,\, \gamma = 1$, while leaving $\epsilon$ and $I$ as parameters for the continuation. We start the continuation near $(\epsilon,I)=(0.3,0.3)$ and compute a triangulation of the manifold with 9006 simplices (including those needed for adaptive refinement) at step size 0.01. For this example, we used $N=7$ Fourier modes. A plot of the proven simplices from the manifold is provided in Figure \ref{fig:FHN0}. The Hopf bifurcation curve appears in Figure \ref{fig:FHN1}.

\subsection{An ODE with a periodic orbit 2-manifold resembling a fish}
\begin{figure}
\centering\includegraphics[scale=0.4]{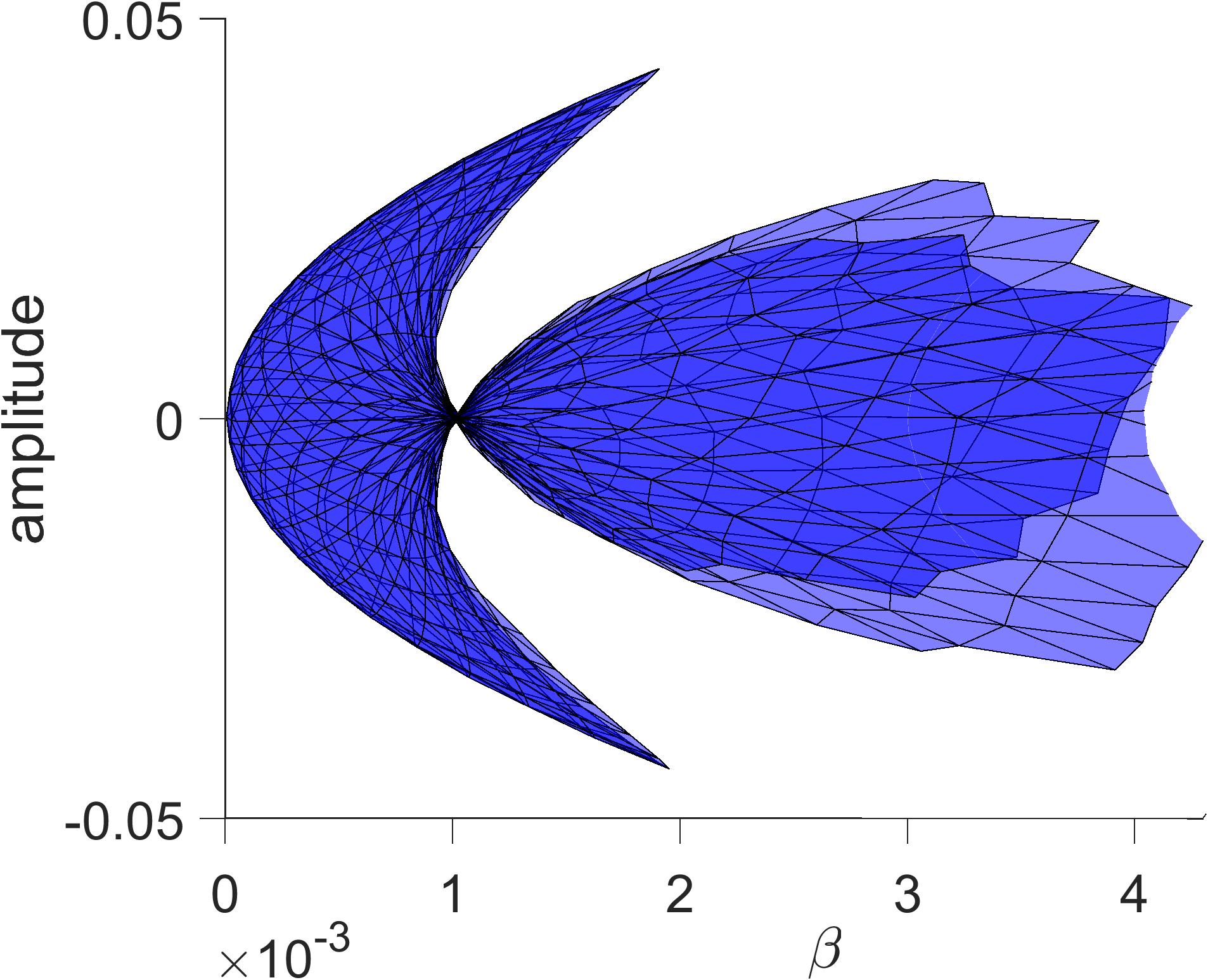} \hspace{8mm} \includegraphics[scale=0.25]{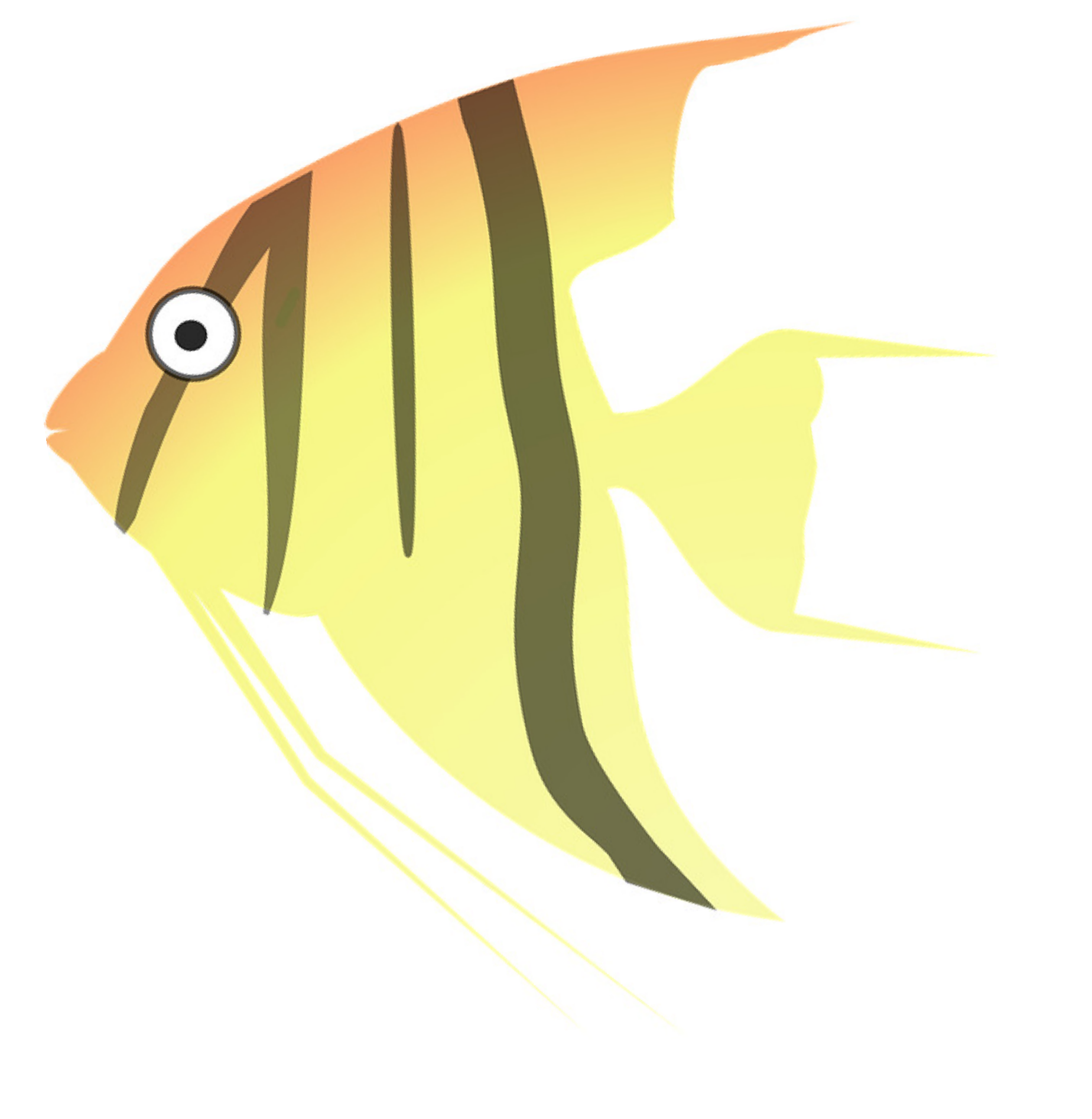}
\caption{Left: Projection of the manifold of periodic orbits for the ODEs \eqref{fish1}--\eqref{fish3} into $\mbox{amplitude}\times\alpha\times\beta$ space, viewed from the $\alpha$ axis. Note the ``pinching" of the manifold, at which this projection becomes singular, near $\beta\approx 10^{-3}$. Right: an illustration of an angelfish, for comparison.}\label{fig:Fish0}
\end{figure}

\begin{figure}[!htb]
\centering\includegraphics[scale=0.4]{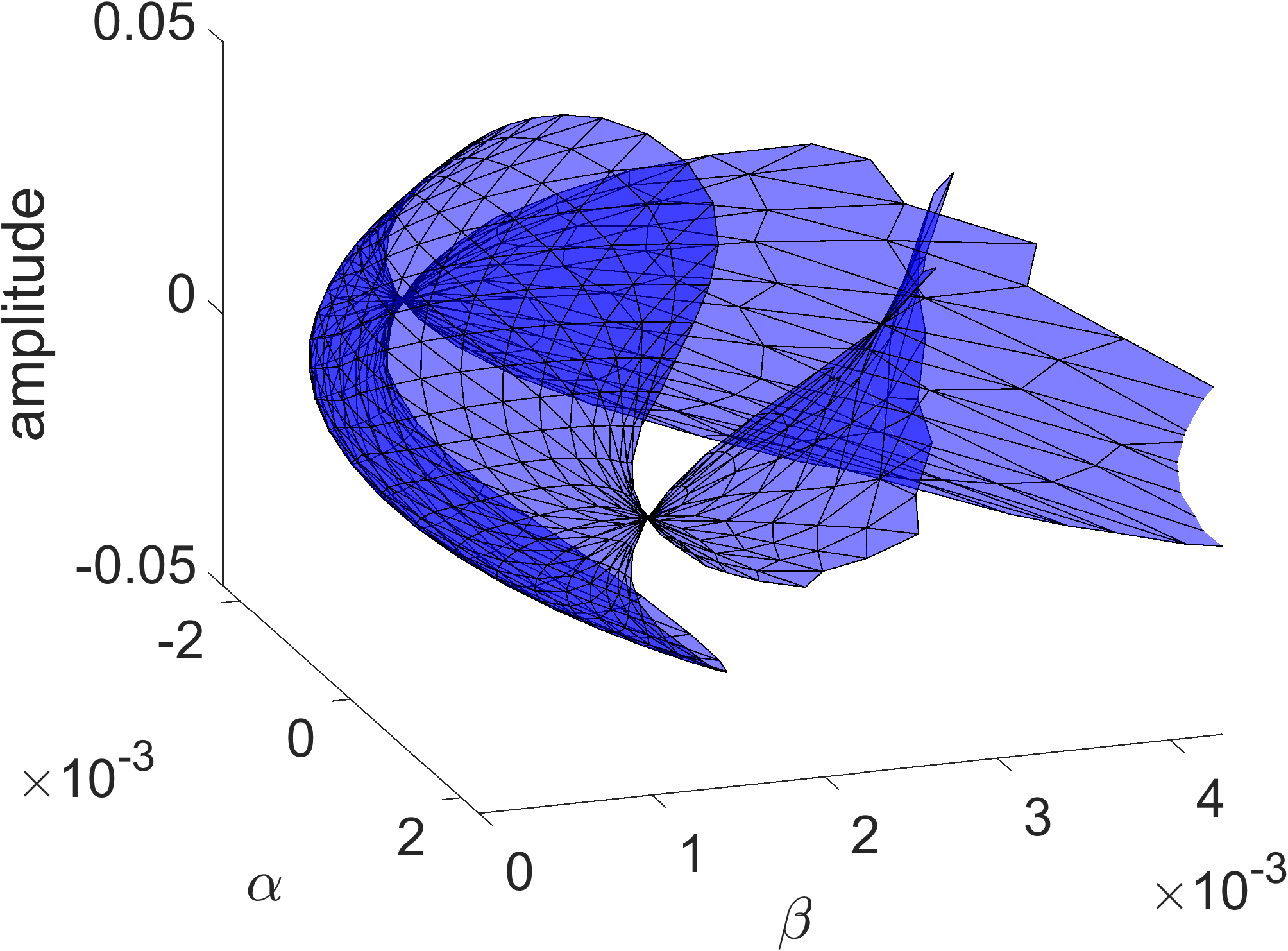} \hspace{4mm}\centering\includegraphics[scale=0.4]{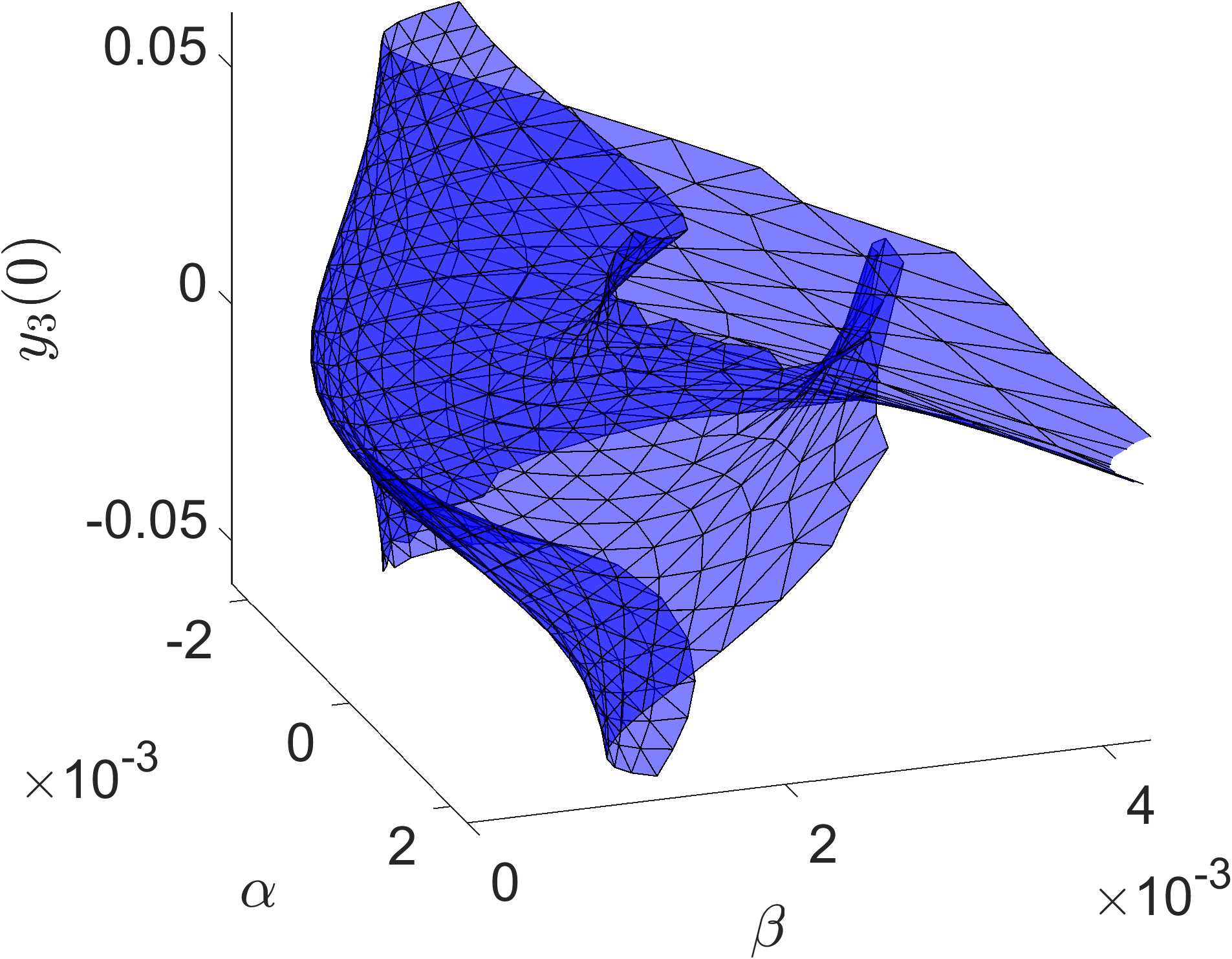} 
\caption{Left: Projection of the manifold of periodic orbits for the ODEs \eqref{fish1}--\eqref{fish3} into $\mbox{amplitude}\times\alpha\times\beta$ space. There are Hopf bubbles prior to the pinching phenomenon that happens near $\beta\approx 10^{-3}$. Right: projection into $\alpha\times\beta\times y_3(0)$ space, where $y_3$ denotes the third component of the periodic orbit in the blown-up coordinates. The pinching point in the left figure projection is caused by a pair of simultaneous folds, clearly visible in the right figure projection. }\label{fig:Fish1}
\end{figure}

\begin{figure}[!htb]
  \begin{minipage}[c]{0.6\textwidth}
   \includegraphics[width=0.8\textwidth]{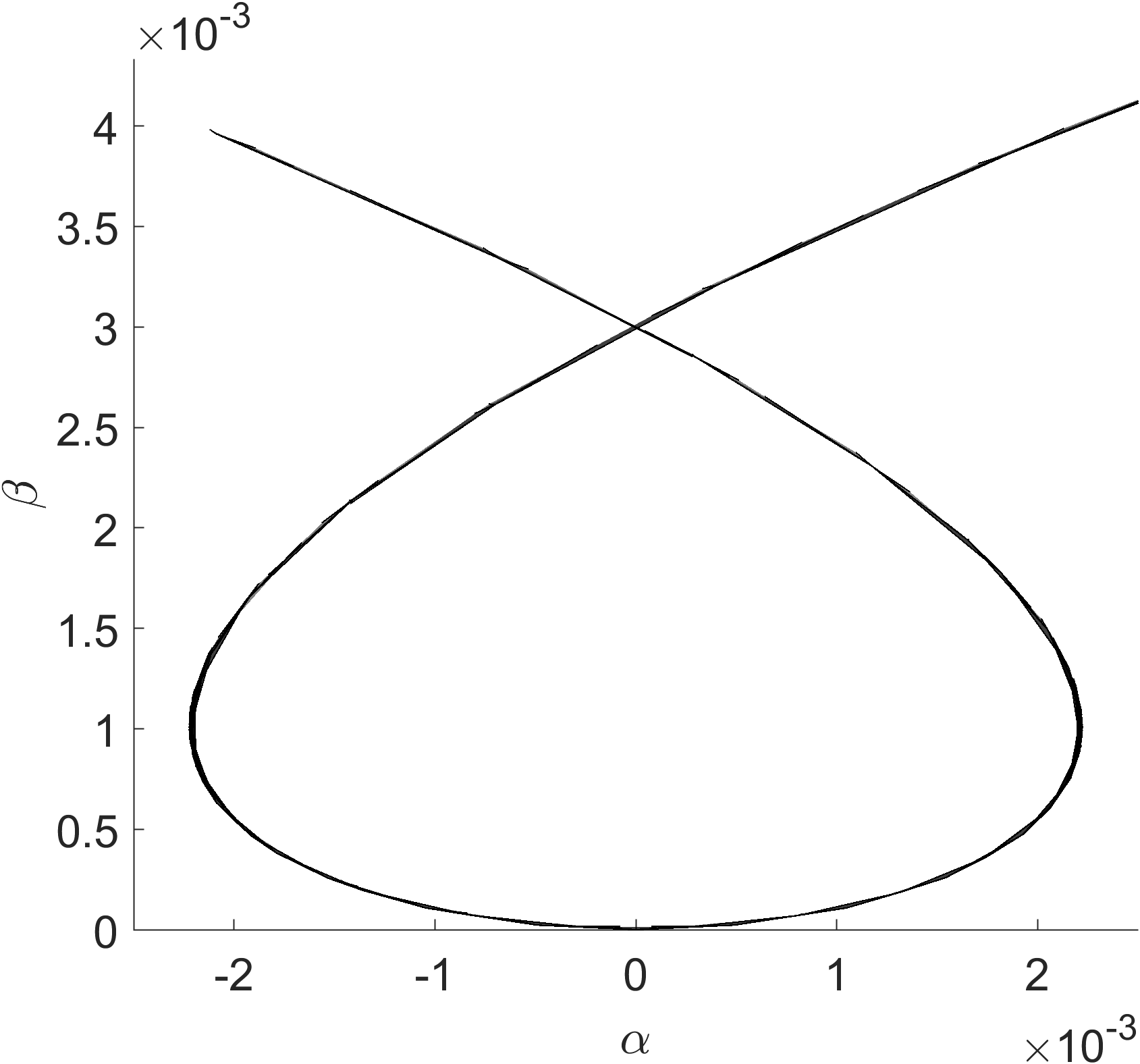}
  \end{minipage}\hfill
  \begin{minipage}[c]{0.4\textwidth}
    \caption{
       Intersection of the surface from the left pane of Figure \ref{fig:Fish1} with the amplitude zero, which corresponds to a Hopf bifurcation curve for the ODE system \eqref{fish1}--\eqref{fish3}. There is a bubble bifurcation at $(\alpha,\beta)=(0,0)$. The apparent self-intersection of the Hopf curve is a consequence of the projection, and does not represent a bifurcation point. Because of how the curve is plotted, skew simplices from them 2-manifold have the effect of making the curve appear ``thicker" in some parts than others.
    } \label{fig:Fish2}
  \end{minipage}
\end{figure}

Consider the three-dimensional ODE system
\begin{align}
\label{fish1}\dot y_1 &= \beta y_1 - y_2 - y_1(y_1^2 + y_2^2 + y_3^2 + \alpha^2)\\
\dot y_2 &= y_1 + \beta y_2 - y_2(y_1^2 + y_2^2 + \epsilon y_3^2 + \alpha^2)\\
\label{fish3}\dot y_3 &= -y_3^5 + 3y_3^3 - 0.01y_3 + 0.1\alpha + 0.01(y_1^2+y_2^2),
\end{align}
for parameters $\alpha,\beta$, and a real control parameter $\epsilon$. When $\epsilon=1$, a change of variables to cylindrical coordinates shows that periodic orbits are in one-to-one correspondence with solutions $(r,z)$ of the set of algebraic equations
\begin{align*}
0&=\beta - r^2 - \alpha^2 - z^2\\
0&=-z^5 + 3z^3 - 0.01z + 0.1\alpha + 0.01r^2.
\end{align*}
When $\epsilon\neq 1$, the radial symmetry in $(y_1,y_2)$ is broken and this change of variables is no longer informative. We set $\epsilon=0.8$ in \eqref{fish1}--\eqref{fish3} and used our validated continuation scheme to rigorously compute a 2-manifold of periodic orbits. In the projection of amplitude and parameters $(\alpha,\beta)$, the result is a figure that qualitatively resembles an angelfish. See Figure \ref{fig:Fish0}. In this projection, the manifold has several folds and appears to exhibit a singularity where it pinches onto a single point. Plotting the manifold in a different projection more clearly allows us to see that this singularity is merely an artifact of the projection; see Figure \ref{fig:Fish1}. The Hopf bifurcation curve is plotted in Figure \ref{fig:Fish2}. For this example, we used $N=9$ Fourier modes and a step size 0.01. We computed a comparatively small portion of the manifold, since the interesting geometry was localized close to $(\alpha,\beta)=(0,0)$. We computed and validated 1007 simplices. This example did not require any adaptive refinement.

\section{Discussion}\label{sec:discussion}
We have proposed validated continuation as an alternative way of exploring degenerate Hopf bifurcations. In combination with rigorous numerics and additional \textit{a posteriori} post-processing, one can prove the existence of Hopf bifurcation curves and bubble bifurcations. The library \textit{BiValVe} is rather flexible, and with the additions of the present paper, can handle multiparameter continuation problems for periodic orbits in ordinary and delay differential equations, near and far from Hopf bifurcations.

Without access to second derivatives of solutions of the zero-finding problem \eqref{Fmap}, it is difficult to prove bubble bifurcations with quadratic folds. That is, we are only able to prove the weaker characterization of Definition \ref{def-generateHopf2}. This is a major barrier in applying the method to delay equations. We believe that a suitable re-formulation of the zero-finding problem, taking into account the additional unbounded operators that result from derivatives with respect to frequency, could resolve the issue.

Another way to prove the ``quadratic fold" part of the bubble bifurcation, would be to compute the second derivatives of Hopf bifurcation curves without using the machinery of sequences spaces. Along these lines, it would be interesting to use pseudo-arclength continuation to continue Hopf bifurcation curves directly. Computer-assisted proofs of isolated Hopf bifurcations in delay differential equations are completed in \cite{Church2022}, and with minimal changes, pseudo-arclength continuation could be used to do continuation of Hopf bifurcations. The map from \cite{Church2022} is finite-dimensional and as smooth as the delay vector field, so the derivatives of the Hopf curve could be rigorously computed that way instead. However, this trick can not be used to prove that $(\alpha,a)\mapsto \pi^\beta u_{h^{-1}(\alpha,a)}$ has a strict local extremum at the bifurcation point. Indeed, the latter map is defined in terms of the periodic orbits themselves, rather than the algebraic properties of the vector field and the eigenvalues of the linearization.

As remarked in \cite{VandenBerg2020a}, the $Z_2$ bound associated to delay periodic orbit validation suffers from a fundamental limitation: it scales linearly with respect to the number of Fourier modes. Therefore, while we have not needed to use many Fourier modes in our examples, it would be very costly (or infeasible) to do continuation of a periodic orbit that required many modes to represent. This is because the $Y_0$ bound is naturally dependent on step size, so even if an isolated solution has an exceptionally good numerical defect, a very small step size might be needed to hedge against a large $Z_2$. In this way, while we can compute manifolds of periodic orbits with delay near (degenerate) Hopf bifurcations, we expect that in large-amplitude regimes or for complicated orbits, completing a validation would be difficult. To compare, the situation is far better for ordinary differential equations. The $Z_2$ bound is generally unharmed by having many Fourier modes, and second derivatives of the solutions can be computed by solving an auxiliary zero-finding problem using similar techniques from rigorous numerics. 

\begin{figure}
  \begin{minipage}[c]{0.6\textwidth}
    \includegraphics[width=0.9\textwidth]{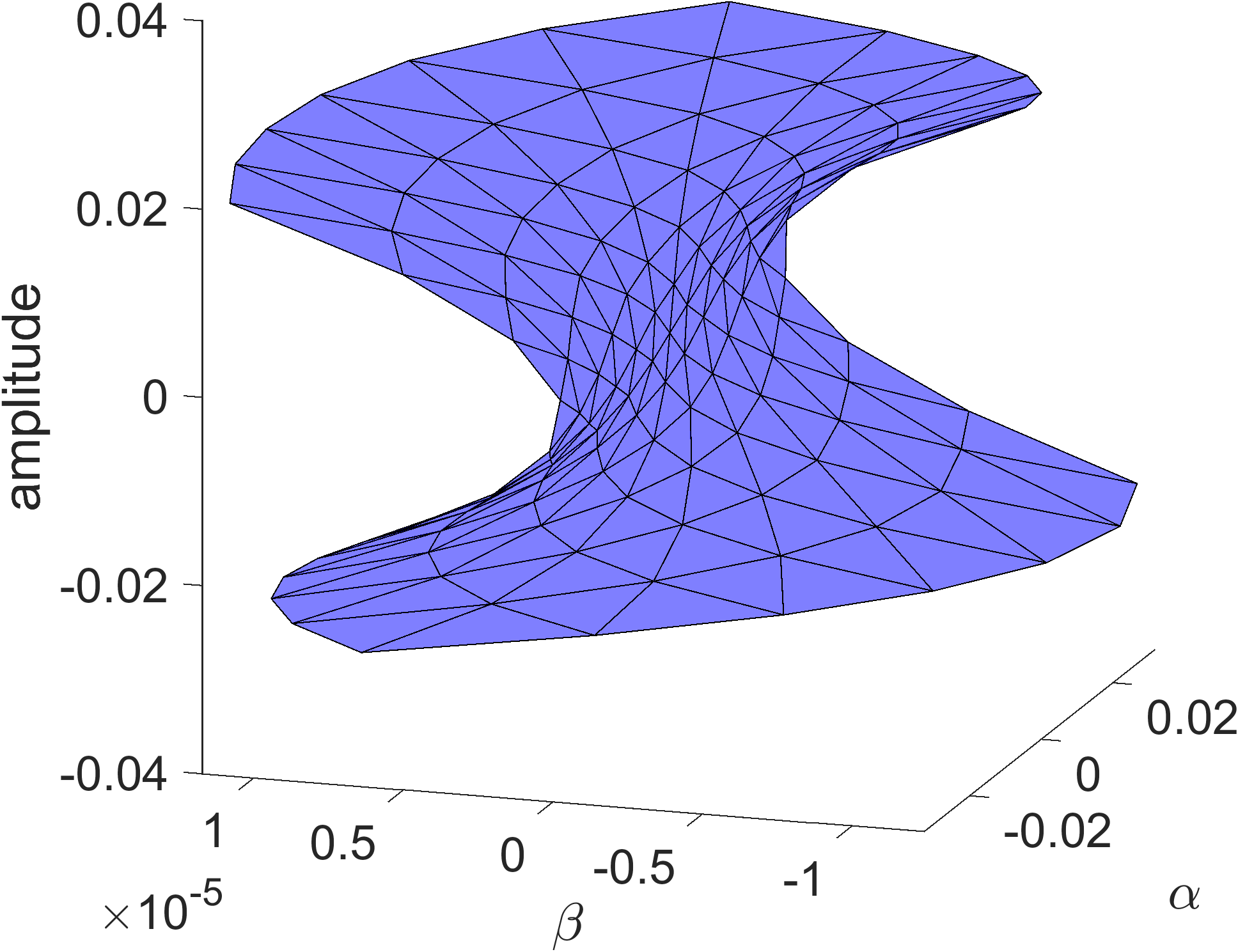}
  \end{minipage}\hfill
  \begin{minipage}[c]{0.4\textwidth}
    \caption{
       Proven section of the manifold of periodic orbits associated to the Bautin normal form. Note that the amplitude is in fact equal to $r$ from equation \eqref{r-Bautin}.
    } \label{fig:Bautin}
  \end{minipage}
\end{figure}

There are other codimension-2+ bifurcations that could be studied from the point of view of validated multi-parameter continuation. For example, the cusp bifurcation should be amenable to this type of analysis, and is simpler than the present work because it involves only bifurcations of fixed points rather than periodic orbits. There is also the Bautin bifurcation, for which the analysis of Section \ref{sec:analytical_proofs} could be replicated. In fact, our continuation scheme is able to validate manifolds of periodic orbits passing through Bautin points. As a very brief final example, recall the normal form \eqref{ex.Bautin1}--\eqref{ex.Bautin2}, which has a Bautin bifurcation at $(x,y)=(0,0)$ at the parameters $(\alpha,\beta)=(0,0)$. Periodic orbits in this ODE are equivalent (by a change of variables to polar coordinates) to scalar solutions $r$ of \begin{align}\label{r-Bautin}0=r(\beta + \alpha r^2 - r^4).\end{align}
Agnostic to this particular representation of the zeroes, our code is able to validate a large section of the manifold of periodic orbits directly from the ODEs. See Figure \ref{fig:Bautin}. As expected, we were able to validated this manifold using very few Fourier modes: three, in this case. 

Hopf bubbles have been observed in the Mackey-Glass equation \cite{Krisztin2011} at the classical parameters, and some of our preliminary investigations suggest that the equation possesses a bubble bifurcation. It would be interesting to use a combination of polynomial embedding and blow-up to investigate this bifurcation. However, the added complexity of using both blow-up and polynomial embedding presents a challenge; the resulting (polynomial) delay vector field ends up being high-order with dozens of distinct nonlinear terms.

\section*{Acknowledgments}
Kevin E.\ M.\ Church is supported by a CRM-Simons Postdoctoral Fellowship. Elena Queirolo is supported by Walter Benjamin Programme DFG-Antrag QU579/1-1.


\begin{thebibliography}{10}

\bibitem{Ashwin2016}
Peter Ashwin, Stephen Coombes, and Rachel Nicks.
\newblock {Mathematical Frameworks for Oscillatory Network Dynamics in
  Neuroscience}.
\newblock {\em The Journal of Mathematical Neuroscience}, 6(1):2, dec 2016.

\bibitem{Bautin}
N. Bautin.
\newblock{Behaviour of dynamical systems near the boundaries of stability regions}.
\newblock{\em OGIZ GOSTEXIZDAT}. Leningrad. In Russian. 1949.

\bibitem{Braza2003}
Peter~A. Braza.
\newblock {The Bifurcation Structure of the Holling--Tanner Model for
  Predator-Prey Interactions Using Two-Timing}.
\newblock {\em SIAM Journal on Applied Mathematics}, 63(3):889--904, jan 2003.

\bibitem{CodeURL}
Kevin~E.M. Church and Elena Queirolo.
\newblock{BiValVe: Bifurcation Validation Venture}.
\url{https://github.com/elenaquei/bubbles/releases/tag/v1}

\bibitem{Church2022}
Kevin~E.M. Church and Jean-Philippe Lessard.
\newblock {Rigorous verification of Hopf bifurcations in functional
  differential equations of mixed type}.
\newblock {\em Physica D: Nonlinear Phenomena}, 429:133072, jan 2022.

\bibitem{Crandall1977}
Michael~G. Crandall and Paul~H. Rabinowitz.
\newblock {The Hopf Bifurcation Theorem in infinite dimensions}.
\newblock {\em Archive for Rational Mechanics and Analysis}, 67(1):53--72,
  1977.

\bibitem{Dupont2016}
Genevi{\`{e}}ve Dupont, Martin Falcke, Vivien Kirk, and James Sneyd.
\newblock {\em {Models of Calcium Signalling}}, volume~43 of {\em
  Interdisciplinary Applied Mathematics}.
\newblock Springer International Publishing, Cham, 2016.

\bibitem{Morshedy2021}
Hassan~A. El-Morshedy and Alfonso Ruiz-Herrera.
\newblock {Asymptotic convergence in delay differential equations arising in
  epidemiology and physiology}.
\newblock {\em SIAM Journal on Applied Mathematics}, 81(4):1781--1798, 2021.

\bibitem{Faria1995}
Teresa Faria and Luis~T. Magalh{\~{a}}es.
\newblock {Normal forms for retarded functional differential equations with
  parameters and applications to hopf bifurcation}, 1995.

\bibitem{Gameiro2016a}
Marcio Gameiro, Jean~Philippe Lessard, and Alessandro Pugliese.
\newblock {Computation of Smooth Manifolds Via Rigorous Multi-parameter
  Continuation in Infinite Dimensions}.
\newblock {\em Foundations of Computational Mathematics}, 16(2), 2016.

\bibitem{Golubitsky1981}
Martin Golubitsky and William~F. Langford.
\newblock {Classification and unfoldings of degenerate Hopf bifurcations}.
\newblock {\em Journal of Differential Equations}, 41(3):375--415, 1981.

\bibitem{Hassard1993}
Brian Hassard and Katie Jiang.
\newblock {Degenerate Hopf bifurcation and isolas of periodic solutions in an
  enzyme-catalyzed reaction model}.
\newblock {\em Journal of Mathematical Analysis and Applications},
  177:170--189, 1993.

\bibitem{Jiang2020}
Weihua Jiang, Qi~An, and Junping Shi.
\newblock {Formulation of the normal form of Turing-Hopf bifurcation in partial
  functional differential equations}.
\newblock {\em Journal of Differential Equations}, 268(10):6067--6102, 2020.

\bibitem{Kuznetsov}
Yuri A Kuznetsov.
\newblock {Elements of Applied Bifurcation Theory}.
\newblock{Springer New York}, 2004.

\bibitem{Krisztin2011}
Tibor Krisztin and Eduardo Liz.
\newblock {Bubbles for a Class of Delay Differential Equations}.
\newblock {\em Qualitative Theory of Dynamical Systems}, 10(2):169--196, oct
  2011.

\bibitem{Latulippe2018}
Joe Latulippe, Derek Lotito, and Donovan Murby.
\newblock {A mathematical model for the effects of amyloid beta on
  intracellular calcium}.
\newblock {\em PLOS ONE}, 13(8):e0202503, aug 2018.

\bibitem{LeBlanc2016}
Victor~G. LeBlanc.
\newblock {A Degenerate Hopf Bifurcation in Retarded Functional Differential
  Equations, and Applications to Endemic Bubbles}.
\newblock {\em Journal of Nonlinear Science}, 26(1):1--25, 2016.

\bibitem{Liu2015d}
Maoxing Liu, Eduardo Liz, and Gergely R{\"{o}}st.
\newblock {Endemic bubbles generated by delayed behavioral response: Global
  stability and bifurcation switches in an sis model}.
\newblock {\em SIAM Journal on Applied Mathematics}, 75(1):75--91, 2015.

\bibitem{Liu2011a}
Zhihua Liu, Pierre Magal, and Shigui Ruan.
\newblock {Hopf bifurcation for non-densely defined Cauchy problems}.
\newblock {\em Zeitschrift fur Angewandte Mathematik und Physik},
  62(2):191--222, 2011.

\bibitem{Margolis1988}
Stephen~B. Margolis and Bernard~J. Matkowsky.
\newblock {New Modes of Quasi-Periodic Combustion Near a Degenerate Hopf
  Bifurcation Point}.
\newblock {\em SIAM Journal on Applied Mathematics}, 48(4):828--853, aug 1988.

\bibitem{Marsden1976}
J.~E. Marsden and M.~McCracken.
\newblock {\em {The Hopf Bifurcation and Its Applications}}, volume~19 of {\em
  Applied Mathematical Sciences}.
\newblock Springer New York, New York, NY, 1976.

\bibitem{Minicucci2021}
Joseph Minicucci, Molly Alfond, Angelo Demuro, David Gerberry, and Joe
  Latulippe.
\newblock {Quantifying the dose-dependent impact of intracellular amyloid beta
  in a mathematical model of calcium regulation in xenopus oocyte}.
\newblock {\em PLOS ONE}, 16(1):e0246116, jan 2021.

\bibitem{Opoku-Sarkodie2022}
Richmond Opoku-sarkodie and Ferenc~A Bartha.
\newblock {Dynamics of an SIRWS model with waning of immunity and varying
  immune boosting period}.
\newblock pages 1--26, 2022.

\bibitem{P.Lessard2020}
J.~{P. Lessard} and J.~{D. Mireles James}.
\newblock {A functional analytic approach to validated numerics for eigenvalues
  of delay equations}.
\newblock {\em Journal of Computational Dynamics}, 7(1):123--158, 2020.

\bibitem{Rustichini1989}
Aldo Rustichini.
\newblock {Hopf bifurcation for functional differential equations of mixed
  type}.
\newblock {\em Journal of Dynamics and Differential Equations}, 1(2):145--177,
  apr 1989.

\bibitem{Sherborne2018}
N.~Sherborne, K.~B. Blyuss, and I.~Z. Kiss.
\newblock {Bursting endemic bubbles in an adaptive network}.
\newblock {\em Physical Review E}, 97(4):042306, apr 2018.

\bibitem{Sneyd2017}
James Sneyd, Shawn Means, Di~Zhu, John Rugis, Jong~Hak Won, and David~I. Yule.
\newblock {Modeling calcium waves in an anatomically accurate three-dimensional
  parotid acinar cell}.
\newblock {\em Journal of Theoretical Biology}, 419:383--393, apr 2017.

\bibitem{VandenBerg2020a}
Jan~Bouwe van~den Berg, Chris Groothedde, and Jean~Philippe Lessard.
\newblock {A General Method for Computer-Assisted Proofs of Periodic Solutions
  in Delay Differential Problems}.
\newblock {\em Journal of Dynamics and Differential Equations}, 2020.

\bibitem{Vandenberg2021_0}
Jan~Bouwe van~den Berg Elena Queirolo.
\newblock {A general framework for validated continuation of periodic orbits in systems of polynomial ODEs.}
\newblock {\em Journal of Computational Dynamics}. 8(1):59, 2021.

\bibitem{VandenBerg2021a}
Jan~Bouwe van~den Berg, Jean-Philippe Lessard, and Elena Queirolo.
\newblock {Rigorous Verification of Hopf Bifurcations via Desingularization and
  Continuation}.
\newblock {\em SIAM Journal on Applied Dynamical Systems}, 20(2):573--607, jan
  2021.

\bibitem{Wang2017a}
Xiaoli Wang, Junping Shi, and Guohong Zhang.
\newblock {Interaction between water and plants: Rich dynamics in a simple
  model}.
\newblock {\em Discrete \& Continuous Dynamical Systems - B}, 22(7):2971--3006,
  2017.

\end{thebibliography}
\end{document}